\documentclass[preprint, nonblindrev]{informs3aa} 

\usepackage[T1]{fontenc}

\usepackage{etex}
\usepackage{dsfont,leftidx, mathrsfs, appendix}
\usepackage{amsmath, amssymb, subfigure, bm, url, graphicx, epstopdf, epsfig,psfrag, xspace}
\usepackage{verbatim} 
\usepackage{subfigure}
\usepackage{algorithmic}
\usepackage{enumerate}
\usepackage{color}
\usepackage{pdfsync}
\usepackage{hyperref}
\usepackage{color}
\def\bydef{\stackrel{\triangle}{=}}

\def\zp{\mathbb{Z}_{+}}
\def\N{\mathbb{N}}
\def\blw{\mathbf{w}}

\def\bla{\mathbf{a}}
\def\blb{\mathbf{b}}
\def\bld{\mathbf{d}}
\def\blx{\mathbf{x}}
\def\bly{\mathbf{y}}
\def\blu{\mathbf{u}}

\def\bua{\mathbf{A}}

\def\buu{\mathbf{U}}
\def\buw{\mathbf{W}}
\def\budta{\mathbf{\Delta}}

\def\R{\mathbb{R}}
\def\pb{\mathbb{P}}

\def\E{\mathbb{E}}
\def\var#1{\mbox{\normalfont Var}\left(#1\right)}

\def\rp{\mathbb{R}_{+}}

\def\calD#1{\mathcal{D}\left(#1\right)}

\def\calE{\mathcal{E}}

\def\lcal{\mathcal{L}}

\def\calX{\mathcal{X}}
\def\calF{\mathcal{F}}
\def\calN{\mathcal{N}}
\def\calH{\mathcal{H}}
\def\calK{\mathcal{K}}
\def\calY{\mathcal{Y}}
\def\calU{\mathcal{U}}

\def\calB{\mathcal{B}}
\def\calC{\mathcal{C}}
\def\calD{\mathcal{D}}
\def\calG{\mathcal{G}}
\def\calI{\mathcal{I}}
\def\calL{\mathcal{L}}
\def\calV{\mathcal{V}}
\def\calW{\mathcal{W}}
\def\calS{\mathcal{S}}

\def\snorm#1{\lt\| #1 \rt\|_T}

\def\uLam{\Lambda}

\def\olz {{\overline{z}}}
\def\old {{\overline{d}}}
\def\dav {d_{\mbox{a}}}
\def\uld {{\underline{d}}}
\def\ulr {{\underline{r}}}

\def\what{\widehat}

\def\nln{\nnb\\}

\def\bbar{\, \big| \, }
\def\hbar{\, \Big| \, }

\def\sk#1#2{\stackrel{(#1)}{#2}}

\def\bpf{\proof{Proof.}}

\def\del#1{}

\newcommand{\szH}{c_\calH}
\newcommand{\szK}{c_\calK}

\def\bslh{\backslash \{h\}}

\newcommand{\flp}{\textsc{FLP}\xspace}
\newcommand{\lptwo}{\textsc{LP2}\xspace}

\newcommand{\beq}{\begin{equation}}
\newcommand{\eeq}{\end{equation}}
\newcommand{\benum}{\begin{enumerate}}
\newcommand{\eenum}{\end{enumerate}}
\newcommand{\btmz}{\begin{itemize}}
\newcommand{\etmz}{\end{itemize}}
\newcommand{\bthm}{\begin{theorem}}
\newcommand{\ethm}{\end{theorem}}
\newcommand{\bdefn}{\begin{defn}}
\newcommand{\edefn}{\end{defn}}
\newcommand{\lt}{\left}
\newcommand{\rt}{\right}
\newcommand{\nnb}{\nonumber}

\theoremstyle{plain}
\newtheorem{theorem}{Theorem}[section]
	
\newtheorem{lemma}[theorem]{Lemma}

\newtheorem{proposition}[theorem]{Proposition}

\newtheorem{definition}[theorem]{Definition}
\newtheorem{defn}[theorem]{Definition}

\usepackage{natbib}
 \bibpunct[, ]{(}{)}{,}{a}{}{,}%
 \def\bibfont{\small}%
 \def\bibsep{\smallskipamount}%
 \def\bibhang{24pt}%
\EquationsNumberedThrough

\begin{document}


\RUNAUTHOR{Massoulie and Xu}

\RUNTITLE{On the Capacity of Information Processing Systems}

\TITLE{On the Capacity of Information Processing Systems}


\ARTICLEAUTHORS{%
\AUTHOR{Laurent Massoulie}
\AFF{Microsoft Research-Inria Joint Centre, 91120 Palaiseau, France} 
\AUTHOR{Kuang Xu}
\AFF{Operations, Information and Technology,  Stanford Graduate School of Business, Stanford, CA 94305, USA} 
}

\vspace{-25pt}

\ABSTRACT{We propose and analyze a family of \emph{information processing systems}, where a finite set of experts or servers are employed to extract information about a stream of incoming jobs. Each job is associated with a hidden label drawn from some prior distribution. An inspection by an expert produces a noisy outcome that depends both on the job's hidden label and the type of the expert, and occupies the expert for a finite time duration. A decision maker's task is to dynamically assign inspections so that the resulting outcomes can be used to accurately recover the labels of all jobs, while keeping the system stable. Among our chief motivations are applications in crowd-sourcing, diagnostics, and experiment designs, where one wishes to efficiently learn the nature of a large number of items, using a finite pool of computational resources or human agents.

We focus on the \emph{capacity} of such an information processing system. Given a level of accuracy guarantee, we ask how many experts are needed in order to stabilize the system, and through what inspection architecture. Our main result provides an adaptive inspection policy that is asymptotically optimal in the following sense: the ratio between the required number of experts under our policy and the theoretical optimal converges to one, as the probability of error in label recovery, $\delta$, tends to zero. \footnote{This version: May 15, 2016. An extended abstract of this paper is accepted for presentation at the Conference on
Learning Theory (COLT) 2016. We are grateful for the  comments from the anonymous COLT referees. }}


\maketitle


\vspace{-20pt}

\section{Introduction}
An increasing number of modern processing systems has been designed and deployed for the purpose of {learning} and {information extraction}. In these applications, which we  refer to broadly as {\bf information processing systems}, a group of experts or servers  is tasked with performing (noisy) inspections on a large collection of jobs, with the objective of uncovering some hidden features associated with each job up to a level of desirable accuracy. Below are some examples: 
\begin{enumerate}
\item \emph{Crowd-sourcing} (\cite{karger2014budget}): a collection of images is dispatched to a group of human agents, where an agent attaches a label to each assigned image based on her own judgment. A decision maker then aggregates agents' responses to produce a ``best'' label for each image. 
\item \emph{Medical diagnostics} (\cite{gerdtz2001triage}):  medical data of patients is reviewed by physicians or nurses with different domains of expertise, with the goal of correctly identifying the patients' diseases. 
\item \emph{Quality management} (\cite{baker1996shmoo}): a set of products undergoes a number of different tests performed by specialized machines, to identify whether a product is faulty and the type of fault it contains. 
\end{enumerate}

The presence of {\bf resource constraints} is a crucial  feature shared across many of these systems:  the amount of processing resources, such as human agents, machines, or computer servers, is finite,  and yet, each inspection or test requires the corresponding resource to commit a non-trivial amount of effort. This raises a natural question: 

\begin{center}
\emph{How much information can we extract using a finite amount of processing resources? }
\end{center}
The main objective of the present paper is to address this question, and gain understanding of the ``capacity'' of an information processing system. We will approach this problem by studying the {minimum required system size}, defined, roughly speaking, as the minimum number of experts needed in order to learn a sufficient amount of information about \emph{every} job in a stream of arrivals, while ensuring system stability. \footnote{Throughout the paper, we will use the term ``expert'' to refer to a single unit of processing resource, with the understanding that it may represent a computer server, testing machine, or human agent, depending on the application.}

We now informally describe our model. The system consists of $m$ experts with different \emph{types} (expertise), where the fraction of experts of type $k$ is $\rho_k$. The system receives a stream of incoming jobs arriving at unit rate, where each job is associated with a \emph{label}, hidden from the decision maker, which is drawn i.i.d.~from some prior distribution, $\pi$. An atomic unit of processing is called an \emph{inspection}: the decision maker may assign an expert to perform an inspection on a job, which occupies the expert for a random period of time, with a mean that depends on the type of the expert. The inspection leads to a (noisy) \emph{outcome}, whose distribution,  $p(k,h,\cdot)$, depends  on both the type of the expert, $k$, and the true label of the job under inspection, $h$. The goal of the decision maker is to assign inspections intelligently, and use the resulting outcomes to produce, for each job, a \emph{classification} (i.e., prediction) of its hidden label. \del{representing what the decision maker believes the true label to be.}We say that the system is \emph{stable} if all jobs will receive a classification in a finite amount of time. 

Note that we cannot have a meaningful discussion on the resource requirement of this system without specifying how accurate the classifications need to be. Indeed, in the absence of any accuracy constraint, the decision maker can simply make up classifications without performing a single inspection, and the system would always be trivially stable. Therefore, we will designate an accuracy parameter, $\delta \in (0,1)$, and demand that the probability of mis-classification for each job be at most $\delta$. Since inspections take up the experts' bandwidth, we expect that a smaller $\delta$ would demand more inspections, which, in turn, require more experts. 

The main goal of this paper is to understand the \emph{minimum number of experts} needed for a given accuracy parameter, $\delta$, and how to achieve it via intelligent policies. This definition of minimum system size is motivated by the intuition that a decision policy that requires the least number of experts to stabilize the system also, in a sense, most efficiently utilizes the processing resources.

{\bf Preview of main result.} The main result of this paper proposes an inspection architecture which, for \emph{any} non-trivial outcome distributions, {asymptotically achieves the minimum system size} in the regime of high accuracy ($\delta\to 0$). Specifically, the ratio between the required number of  experts under our policy, $m$, and that of a theoretical minimum, $m^*$, converges to $1$ uniformly across all prior distributions of job labels, as $\delta \to 0$: 
\begin{equation}
\sup_{\pi}\frac{m}{m^*} = 1+ \mathcal{O}\left(\sqrt{\frac{\ln\ln(1/\delta) }{\ln (1/\delta)}}\right), \quad \mbox{as $\delta \to 0$}. 
\end{equation}
Moreover, the policy does not require knowing the prior distribution, $\pi$, but adapts to it automatically.

We conclude this section by highlighting two main design challenges in creating an efficient decision policy for information processing systems. The first challenge arises from the fact that the processing resources are often heterogeneous, a result of the variations in expertise, machine functionality, or personal trait. In our model, this is captured by the fact that the  outcomes distributions, $\{p(h,k,\cdot)\}_{h\in \calH, k\in \calK}$, may vary significantly depending on the type of the inspecting expert. A key implication of such heterogeneity is that the decision policy must be sufficiently adaptive to a job's past inspection outcomes, because depending on what we ``believe'' the job's label to be, some combinations of expert types may be more efficient than others. For instance, consider a testing system where job labels correspond to three domains that might have caused an error in a product: $\{\mbox{network, hardware, software}\}$, and the experts are of two kinds: 
\begin{enumerate}
\item General experts: they know a little about all domains, and produce noisy inspection outcomes. 
\item Specialists:  who have deep expertise in one of the three domains but know nothing about the other two.  Their inspections contribute strong signals towards confirming an error in their own domain, but are otherwise non-informative if the job's true label lies elsewhere. 
\end{enumerate}
Intuitively, a sensible decision policy for this setting should be adaptive: a job can be first inspected by a few general experts to zone in on a possible domain, and depending on their opinions, the job will then be sent to the corresponding specialists to further ``confirm'' the  {diagnosis}, with some small possibility of back-and-forth if the initial {diagnosis} had been incorrect. In contrast, sending a job directly to a specialist at the beginning would risk the specialist being of the wrong kind, and never sending a job to any specialist would overwhelm the generalists who can only provide limited information per inspection; in both cases, there would be waste of processing resources.

 The second challenge is more nuanced, and stems from the combined effect of the resource constraint and expert heterogeneity: the ``optimal'' course of inspections for an individual job does not only depend on the types of experts available, but also on the \emph{prior distribution} governing the proportions of labels among its fellow jobs. In other words, there will be  \emph{contention} among jobs as they ``compete'' for the same processing resources. For instance, in the above-mentioned testing example, suppose there emerges a {disproportionately} large fraction of jobs with the label ``network'', while the fraction of specialists in ``network'' remains fixed. In such a case, it is plausible that some of the generalists may now need to be enlisted to provide assistance by performing more inspections on these {jobs} than they would before. Therefore, the mixture of inspections that a job receives may shift as the prior distribution changes, which shows that a good decision policy cannot be overly centered around individual jobs, and must be aware of the overall arrival pattern.

\subsection{Related Research}

The present paper intersects with two main areas of research: statistics and dynamic resource allocation. Most related to our work is the literature on sequential hypothesis testing in statistics, dating back to the seminal work of \cite{wald1945sequential}, which studies the problem of distinguishing hypotheses of a distribution by sequentially drawing samples from it, with the objective of minimizing a combined cost involving the number of samples and the resulting probability of error (see also \cite{siegmund2013sequential} and the references therein). Notably, \cite{chernoff1959sequential} considers an important generalization of Wald's problem, where, instead of drawing samples from the same distribution and deciding when to stop, the decision maker has the additional freedom to choose from \emph{multiple} available experiments, and the distribution of the experimental outcome depends on both the true hypothesis and the type of the experiment. Chernoff identifies a dynamic testing policy that asymptotically achieves the minimum sample complexity, and computes explicitly the leading factor via the solution to a zero-sum game. The current paper draws inspiration from this literature, and especially the multi-experiment version of \cite{chernoff1959sequential}. Yet, there are some key differences:  while sequential hypothesis testing aims to reduce the sample complexity associated with testing a \emph{single} distribution, we are interested in performing tests for \emph{multiple} jobs simultaneously using finite resources. The resource constraint creates {contention} and coupling among the jobs, and as was {alluded} to in Introduction, policies designed to minimize the sample complexity associated with testing a single hypothesis do not easily extend to our problem.  

Our work is also related to the literature on dynamic resource allocation, and particularly multi-class, multi-server queueing networks (\cite{harrison1999heavy, TW08, TsiXu15Flex}). In our model the inspection outcomes may have different distributions depending on the expert's type and the job's label, which is roughly analogous to a queueing network where the service rate depends on the class of the server and the job being served. However, our system differs from this literature in two major aspects. First, while in conventional processing systems jobs typically come with an \emph{exogenous} service requirement, also referred to as workload or job size, in our system a job's service requirement is defined \emph{endogenously}  in relation to how much information needs to be gathered in order to uncover its label. Second,  in a queueing network the class of an incoming job is typically known to the decision maker, while in our model the job labels are hidden, and the very objective of processing is to uncover them. 

There have also been a growing interest in the connections between information acquisition and resource allocation \cite{alizamir2013diagnostic,Bimpikis2015learning, jkk2016,  harrison2015investment, levi2015testing}. They are similar to our work in spirit, but differ in models and objectives. The authors of \cite{alizamir2013diagnostic} study a single-server queue where the server decides on how many tests to perform on each job to reach a binary {diagnosis}, while achieving an optimal delay-accuracy trade-off. The  information structure in \cite{alizamir2013diagnostic} is substantially more restricted than in our model: the tests are identical, and outcomes and job types are binary; on the other hand, \cite{alizamir2013diagnostic} analyzes queueing delay which we do not consider. The recent papers \cite{jkk2016,Bimpikis2015learning} study processing systems with heterogeneous job and server types and capacitated processing resources, where, similar to our model, it is important for the decision maker to learn about the job types in order to identify the best processing scheme. These models differ from ours in that learning serves as a means towards another  objective, such as maximizing total rewards in \cite{jkk2016} or improving  throughput or delay in \cite{Bimpikis2015learning, levi2015testing}, whereas information extraction is the intrinsic purpose of processing in our problem. 

On the technical end, we build on several existing techniques and ideas. To establish the stability of our decision policy, which involves three interconnected stages, we will leverage a program pioneered by \cite{rybko1992ergodicity} and \cite{dai1995positive} in the context of queueing networks, which approximates the system dynamics using a certain fluid limit, and subsequently {uses} a contraction property of the fluid limits to derive the stability of the original system. A sub-routine of our inspection policy dynamically creates workload vectors by dynamically solving a linear program to minimize the incremental change to a potential function, which is reminiscent of, and inspired by, the family of max-weight scheduling policies (\cite{tassiulas1992stability}). Finally, to derive upper and lower bounds on the probability of error under an  inspection policy, we make elementary uses of well-known techniques in probability theory and statistics, such as changes of measures and concentration inequalities for martingales. 

\section{Model and Metrics}
\label{sec:setup}

\subsection{The Model} 

\emph{System primitives.} The system evolves in continuous time, indexed by $t\in \rp$. There is a stream of jobs that arrives to the system according to a Poisson process with rate $\lambda_0>0$. Without loss of generality,  we assume that $\lambda_0=1$, because one can simply scale the expressions of system size by $\lambda_0$ so that the results in the paper apply to other values of $\lambda_0$ as well. We index jobs by the order in which they arrive, and refer to the $i$th job that arrives to the system as \emph{job $i$}.  Job $i$ is associated with a \emph{label}, $H_i$, which belongs to a finite set, $\calH$, whose cardinality is $\szH$. The job labels are independent and identically distributed according to a prior distribution, $\pi=\{\pi_h\}_{h \in \calH}$, with $\pb(H_1 = h) = \pi_h$, and the labels are unknown to the decision maker. We assume all entries of $\pi$ to be positive.  

The system is equipped with $m$ \emph{experts}. Let $\calE$ be the set of experts. Each expert is associated with a \emph{type}, $ k $, from a finite set, $ \calK$,  with cardinality $\szK$. The number of type-$k$ experts is $ \rho_k m$, $k\in  \calK$, where $\rho_k$ is the fraction of type-$k$ experts in the system, with $\rho_k > 0$ and $\sum_{ k \in \calK} \rho_k =1$. We will refer to $\{\rho_k\}_{k \in \calK}$ as the \emph{expert mixture}.


\paragraph{Inspections and resource constraints.} An {expert} can be called upon to perform an inspection of any job  in the system. At the end of an inspection, a random \emph{outcome} is produced which takes values in a finite set, $\calX$. We denote by $X_{i,j}$ the outcome of the $j$th inspection performed on job $i$. Suppose that the true label of job $i$ is $h$, and that the $j$th inspection is performed by an expert of type $k$, then $X_{i,j}$ is a random variable  distributed according to the outcome distribution,  $p(h, k, \cdot)$, and is independent from all other parts of the system. Note that this implies also that an expert may inspect a job multiple times, producing i.i.d.~outcomes. The diversity in the set of outcome distributions, $\{p(h,k,\cdot)\}_{h\in \calH, k \in \calK}$, captures the possibility that experts of various types may have different expertise. We assume the set of outcome distributions is known to the decision maker. 

The pool of experts is \emph{resource constrained},  in the sense that each individual expert can only perform,  on average, a finite number of inspections per unit time. Formally, at any time $t$, we assume that each of the $m$ experts can be in one of the two states: IDLE and BUSY, and all experts are initialized in state IDLE. An IDLE expert of type $k$ can be assigned to initiate an inspection of a job currently in the system. Once the inspection starts, the expert enters the state BUSY for a duration that is exponentially distributed with mean $1/\mu_k$, $\mu_k >0$, independent from the rest of the system. We refer to $\{\mu_k\}_{k \in \calK}$ as the \emph{inspection rates}, since $\mu_k$ corresponds to the average number of inspections that a type-$k$ expert can perform in unit time. The expert returns to the IDLE state once the inspection is completed. The inspections are non-preemptive, so that an expert cannot start inspecting a different job before the previous inspection has been finished.  We assume that \del{a single expert cannot inspect multiple jobs simultaneously, but }multiple experts are allowed to inspect the same job at the same time. The latter assumption is motivated by applications where jobs, such as images or data files, can be duplicated at relatively low costs or be accessible to multiple experts concurrently. 

 Note that the inspection rate depends on the expert's type but not the job's hidden label.  Similar to the earlier assumption that the arrival rate $\lambda_0=1$, without loss of generality, we assume that $\{\mu_k\}_{k \in \calK}$ is normalized so that the average inspection rate across different expert types is 1:\footnote{In particular, if $\overline{\mu} \neq 1$, then one can scale our results that concern the minimal system size by $1/\overline{\mu}$. }
\begin{equation}
\overline{\mu} \bydef \sum_{k \in \calK} \rho_k \mu_k= 1. 
\label{eq:avgRate}
\end{equation}

\paragraph{Departure rule.} At any time $t$, the system operator can choose to let a job depart from the system. Upon job $i$'s departure from the system, the operator must produce a \emph{classification}, $\what{H}_i$, representing her belief of job $i$'s true label. We say that there is an \emph{error} if the classification does not match the true label, i.e., if $\what{H}_i \neq H_i$.

\subsection{Conditions on Outcome Distributions}

The set of outcome distributions, $\{p(h,k,\cdot)\}_{h\in \calH, k \in \calK}$, plays a central role in our model, because the job labels can only be differentiated via the inspection outcomes. In the present paper, we allow for essentially any outcome distribution over $\calX$, with the exception of two conditions. Informally, we assume that $(1)$  all outcomes are ``noisy'', so that no single outcome can distinguish between two job labels with certainty, and $(2)$   for any two job types, there exists at least one type of experts who can distinguish them. 	Neither assumption leads to a severe loss of generality, as we explain shortly. 

We now give formal definitions of the two conditions. In our regime of interest, where the target classification error is small, it turns out that an important measure of the informativeness of an inspection outcome is that of KL-divergence, defined as follows. Fix $i,j\in \N$,  $h, l$ in $\calH$,  and $k\in \calK$, denote by $Z_{i,j}(h, l, k)$ the likelihood associated with the $j$th inspection to job $i$ done by a type-$k$ expert, as follows: 
\begin{equation}
Z_{i,j}(h, l, k) = \ln \frac{p(h, k,X_{i,j})}{p(l, k,X_{i,j})}.
\label{eq:zij}
\end{equation}
The KL-divergence from the outcome distribution $p(h, k, \cdot)$ to $p( l, k, \cdot)$, denoted by $D(h, l, k)$, is defined by the expected value of $Z_{i,j}(h, l, k)$ conditional on the true label of job $i$ being $h$: 
\begin{equation}
D(h, l, k)= \E\left(Z_{1,1}(h,l,k) \bbar H_1 = h\right) = \sum_{x \in \calX}  p( h, k, x) \ln \frac{p(h, k, x)}{p(l, k, x)}. 
\label{eq:dklDef}
\end{equation}
Intuitively, the value of $D(h,l,k)$ captures the ability of an expert of type $k$ in telling apart whether a job's label is $h$ versus $l$; a higher value of $D(h,l,k)$ indicates that the expert's inspection, on average, provides stronger evidence that the true label of a job {is} more likely to be $h$ than $l$.

We assume that the outcome distributions $\{p(h,k,\cdot)\}_{h\in \calH, k\in \calK}$ satisfy two conditions as follows, expressed in terms of KL-divergence. 

\begin{enumerate}

\item The outcome distributions should be sufficiently diverse so that accurate classifications are possible.  For any two distinct job labels, $h, l \in \calH$, we assume that there exists at least one expert type, $k$, for whom the  outcome distributions, $p(h,k,\cdot)$ and $p(l,k,\cdot)$, are non-identical. This is equivalent to saying that there exists $\uld>0$, such that 
\begin{equation}
\uld = \min_{h,l\in \calH, h\neq l} \, \max_{k\in \calK} D(h,l,k) >0. 
\label{eq:minimumInfo}
\end{equation} 
\item All KL-divergences should be finite, so that a single inspection cannot distinguish two job labels with absolute certainty: 
\begin{equation}
\old  = \max_{\stackrel{h, l \in \calH,\, h\neq l}{k\in \calK}} D(h,l,k) < \infty. 
\label{eq:dbar}
\end{equation}
\end{enumerate}

We note that the above two conditions do not significantly restrict the outcome distributions. The first condition is in fact necessary for the problem to be non-trivial, for otherwise one would not be able to distinguish between jobs of labels $h$ and $l$. The second condition requires all inspections to be ``noisy'', and it is a natural assumption in applications where the inspections have some intrinsic variability, such as when inspections are performed by human agents, or by machines that are subject to idiosyncratic noise. 


\section{An example of outcome distributions}

For concreteness, we describe here a simple example of a family of outcome distributions. Consider the case where a job is an image, containing one out of three possible animals, with $\calH = \{\mbox{cat}, \mbox{dog}, \mbox{rabbit}\}$. There are three types of experts, with $\calK= \{1,2,3\}$, who are asked to perform inspections leading to binary outcomes (e.g., ``like'' or ``dislike''), with $\calX = \{0,1\}$. Fix $p,q\in (0,1)$, $p \neq q$. Denote by $A$ and $B$ the Bernoulli distribution with mean $p$ and $q$, respectively. The outcome distributions, $p(h,k,\cdot)$, are given by the following matrix, where the column corresponds to the labels and the row to the expert types: $\left( 
\begin{array}{ccc}
A & A & B \\
A & B & A \\
B & A & A
\end{array} 
\right)$. 
For instance, the entry $(2,1)$ indicates that a type-$2$ expert's inspection of an image containing a `cat' is distributed according to $A$. Note that the outcomes are statistically identical when a type-$1$ expert inspects an image with a cat versus a dog, but are different from the outcome when the image contains a rabbit. However, it is not difficult to  see that if we were to obtain many inspections from any two expert types, one could eventually uncover the true label of an image. This example illustrates that  for the overall processing system to be effective, it is not necessary that a single expert be able to distinguish all job labels.

\subsection{Inspection Policies and Performance Metrics}


\paragraph{Inspection policies.} To facilitate our discussion, we now introduce the concept of an inspection policy. An inspection {policy}, $\psi$, has access to the entire system state and all past history. At any time $t$, it has the ability to: 
(1) let an IDLE expert initiate an inspection on a job; 
(2) let a job $i$ depart from the system, in which case the policy will have to produce a classification, $\what{H}_i$, for the job's label, ${H}_i$.  In addition, the inspection {policy} can take as input the following parameters: $(1)$
 the number of experts, $m$, the expert mixture, $\{\rho_k\}_{k \in \calK}$, and the inspection rates, $\{\mu_k\}_{k \in \calK}$;  $(2)$ the outcome distributions, $\{p(h,k, \cdot)\}_{h \in \calH, k \in \calK}$; 
$(3)$ an accuracy parameter, $\delta$; $(4)$ (potentially) the prior distribution, $\pi$, of the job labels. An inspection policy that does not require the knowledge of the prior distribution, $\pi$, is said to be \emph{prior-oblivious}. 

In our system, there are two main performance criteria for an inspection policy that are of interest. First, we would like an inspection policy to accurately recover the true labels of all jobs, quantified in the following definition. 
\begin{definition}[$\delta$-accuracy]
A policy is {\bf $\delta$-accurate}, if for all $h\in \calH$ and $i \in \N$,  we have that $\pb(\what{H}_i\neq h \bbar H_i = h) \leq \delta$. \end{definition} 
That is, the resulting probability of misclassification on any job of type $h$ is at most $\delta$ under a $\delta$-accurate policy. Since the fundamental task of our processing system is to recover job labels reliably, we assume throughout the paper that an inspection policy should always be $\delta$-accurate, where $\delta$ is the accuracy parameter set by the decision maker.

In addition to accuracy, a second important benchmark is stability. That is, every job should be able to depart from the system in a finite amount of time. Formally, denote by $Q(t)$ the total number of jobs in the system at time $t$. We say that the system is {\bf stable} under an inspection policy, if the resulting process $\{Q(t)\}_{t\in \rp}$ is positive recurrent.  Compared to the definition of accuracy, the notion of stability is more delicate, because whether an inspection policy can stabilize a system depends on the \emph{relative} magnitudes between the number  of experts, $m$, and the accuracy requirement, $\delta$: because the inspections are noisy, as $\delta$ decreases, each job will require a larger number of inspections to achieve a desired classification accuracy. Since the arrival rate of jobs is assumed to be fixed, this further implies that the number of experts, $m$, must also grow accordingly. 

Therefore, as {alluded} to in the Introduction, a natural way of assessing how efficient an inspection policy is at stabilizing the system is to measure the minimum system size (i.e., $m$) needed in order for the system to be stable for a given accuracy target, $\delta$. This inspires the following notion of {\bf resource efficiency}, where we compare the minimum system size required by an inspection policy against that of a theoretical optimal, as follows. Fix an expert mixture, $\{\rho_k\}_{k \in \calK}$, inspection rates, $\{\mu_k\}_{k \in \calK}$, and outcome distributions, $\{p(h,k, \cdot)\}_{h\in \calH, k\in \calK}$. For a policy, $\psi$, define $m_\psi(\delta, \pi)$ as the smallest system size required under $\psi$ in order to ensure stability: 
\begin{equation}
m_\psi(\delta, \pi) = \min\{m \in \N: \mbox{given $\delta$ and $\pi$, a system with $m$ experts is stable under  $\psi$}\}. 
\end{equation}
Given prior distribution $\pi$ and $\delta>0$, we define $m^*(\delta, \pi)$ as the smallest system size for which there exists  a $\delta$-accurate inspection policy that stabilizes the system. That is, $m^*(\delta, \pi)$ represents the minimal amount of processing resources required to ensure stability under an ``optimal'' inspection policy. The following definition serves as the main performance metric of this paper. 

\begin{definition} We say that an inspection policy, $\psi$, is {\bf resource efficient}, if 
\begin{equation}
\limsup_{\delta \to 0} \frac{{m}_\psi(\delta, \pi)}{m^*(\delta, \pi) } =1, \quad \mbox{for all prior distribution, $\pi$.}
\end{equation}
We say that $\psi$ is {\bf strongly resource efficient}, if the above convergence occurs uniformly over all prior distributions: 
\begin{equation}
\limsup_{\delta \to 0} \, \sup_{\pi}\frac{{m}_\psi(\delta, \pi)}{m^*(\delta, \pi) } =1. 
\end{equation}
\end{definition}

\section{Main Result}

The main result of the current paper is the following theorem. 

\begin{theorem}
\label{thm:main}  Fix an expert mixture, $\{\rho_k\}_{k \in \calK}$, inspection rates, $\{\mu_k\}_{k \in \calK}$, and outcome distributions, $\{p(h,k, \cdot)\}_{h\in \calH, k\in \calK}$. There exists a prior-oblivious, strongly resource efficient inspection policy, $\psi$. In particular, there exist $c_0, \delta_0>0$, such that
\begin{equation}
\sup_{\pi}\frac{{m}_\psi(\delta, \pi)}{m^*(\delta, \pi)}  \leq 1+ c_0 \sqrt{\frac{\ln\ln (1/\delta) }{\ln (1/\delta)}}, \quad \forall \delta \in (0,\delta_0). 
\label{eq:mainIneq}
\end{equation}
\end{theorem}

%

We highlight two important features of the theorem. First, the inspection policy is strongly resource efficient, which implies that its performance guarantee in comparison to the theoretical optimal holds \emph{independently} of the prior distribution.  Second, the inspection policy is \emph{prior-oblivious} so that it can operate without any knowledge of the prior distribution of the job labels. This feature is especially important for our problem, since knowing the prior distribution would have likely required first learning the labels of the incoming jobs, which is the very task that we are trying to solve! Moreover, because a prior-oblivious policy automatically adapts to any prior distribution, it is also more robust if the prior distribution were to shift over time, a likely scenario for many applications.  We will provide in Section \ref{sec:DefOfPolicy} a complete description of the strongly resource-efficient inspection policy in Theorem \ref{thm:main}, which is based on a three-stage architecture that first generates a \emph{coarse} label estimate for each job, and subsequently uses the majority of the processing resources to \emph{verify} the validity of the coarse estimates, in an adaptive manner. This policy also inspires a simple heuristic algorithm, discussed in Appendix \ref{sec:heuristicPolicy}, which can be much easier to implement in practice.

\section{Proof Overview and Preliminaries}
\label{sec:profideas}

\subsection{The Main Ideas}
\label{sec:mainIdeas}

The remainder of the paper is devoted to the proof of Theorem \ref{thm:main}. Before delving into the details, we will start by illustrating the main ideas of the proof.  Our main goal is to design an inspection architecture to extract information efficiently using a finite number of experts. We will break this general problem further into three sub-problems, in the following order: 
\begin{enumerate}[(a)]
\item What type  of information is sufficient for making accurate classification decisions? 
\item How {much} information do we need to gather for an \emph{individual} job, via inspections, in order to produce an accurate classification of its label?
\item How can we gather such information for \emph{all} jobs simultaneously in a resource-efficient manner? 
\end{enumerate}

We now address each of the three points in order. For point $(a)$, the following notion of cumulative log-likelihood ratio, a concept widely used in statistics, will be central in quantifying information in our problem. Fix $i\in \N$ and $t\in \rp$. Denote by $N_{i,t}$ the total number of inspections received by job $i$ by time $t$, and by $K_{i, j}$ the type of the expert who performed the $j$th inspection on job $i$. For $h$ and $l$ in $\calH$, we define the cumulative log-likelihood ratios for job $i$ at time $t$ as
\begin{equation}
S_{i, t}(h, l) = \sum_{j=1}^{N_{i,t}}Z_{i,j}(h, l, K_{i,j}) =  \ln \lt( \prod_{j=1}^{N_{i,t}}\frac{p(h, K_{i,j},X_{i,j})}{p(l, K_{i,j},X_{i,j})} \rt), \quad h, l \in \calH, 
\label{eq:Sinsample}
\end{equation}
where $Z_{i,j}(h,l,k)$ is the log-likelihood ratio defined in Eq.~\eqref{eq:zij}. Intuitively, the fact that $S_{i, t}(h, l)>0$ implies that given the inspections and their outcomes up till time $t$, the label of job $i$ is \emph{more likely} to be $h$ than $l$, and such likelihood intensifies as the value of $S_{i, t}(h, l)$ increases. As we will see in a moment, the set  $\{S_{i, t}(h, l)\}_{h, l \in \calH}$ serves as a summary statistic that is sufficient for producing classifications for job $i$'s label. 

Point $(b)$ concerns the quantity of information needed to make an \emph{accurate} classification. In light of the preceding discussion, we could equally ask: at the time when a classification has to be made about job $i$'s label, what conditions should  $\{S_{i, t}(h, l)\}_{h, l \in \calH}$ satisfy in order for the classification error to be small?  The next lemma provides such a sufficient condition.  Define $\what{H}_{i,t}$ as the \emph{maximum-likelihood (ML) estimator} for the true label of job $i$, $H_i$, given the inspections performed on job $i$ up till time $t$, i.e.,
\begin{equation}
	\what{H}_{i,t} \in \{h \in \calH: S_{i, t}(h, l) \geq 0, \, \forall l \neq h\}=  \arg\max_{h \in \calH}   \prod_{j=1}^{N_{i,t}}p(h, K_{i,j},X_{i,j}), 
	\label{eq:MLdef}
\end{equation}
with ties broken arbitrarily\footnote{Note that the equality in the above equation follows from the definition: the most likely label is also the one that is no less likely than any other labels.}. We will denote by $S^F_{i}(h,l)$ and $\what{H}_{i}^F$ the value of $S_{i,t}(h,l)$ and $\what{H}_{i,t}$ at the time when job $i$ departs from the system, respectively.  We have the following lemma, which is a special case of a more general result, Lemma \ref{lem:SijSufficient2}, in Appendix \ref{app:lem:SijSufficient2}.   
\begin{lemma}
\label{lem:SijSufficient}
Fix $i \in \N$ and $x>0$. Denote by $\calG_x$ the event:
\begin{equation}
\calG_x = \{\exists h'\in \calH, \mbox{ s.t. } S_{i}^F(h', l) \geq x, \quad \forall l\in \calH, l \neq h' \}.
\end{equation}
We have that
\begin{equation}
\pb(\what{H}^F_{i} \neq h \, , \, \calG_x\bbar H_i=h)\leq \szH\exp(-x), \quad \forall h \in \calH. 
\end{equation}
\end{lemma}

Lemma \ref{lem:SijSufficient} shows that if, for a large value of $x$, the event $\calG_x$ occurs with high probability under any job label, then the resulting probability of mis-classification must be small. This achievability result is complemented by the following converse, which states that in order for any inspection policy to be $\delta$-accurate, the expected value of the cumulative log-likelihood ratio, $\E(S_{i}(h,l)) $, must satisfy a lower bound that is approximately $\ln(1/\delta)$. The proof of the lemma utilizes a coupling argument similar to that in \cite{wald1945sequential} for establishing a lower bound on sample complexity in sequential hypothesis testing, and is given in Appendix \ref{app:lem:EofShl}.

\begin{lemma} 
\label{lem:EofShl}
Fix $\delta \in (0,1)$. If an inspection policy is $\delta$-accurate, then for all $i \in \mathbb{N}$ and $h\in \calH$, 
\begin{equation}
\E \lt(S^F_{i}(h,l) \bbar H_i = h\rt) \geq (1-\delta)\ln\frac{1-\delta}{\delta} - e^{-1}, \quad \forall l \in \calH \bslh. 
\end{equation}
\end{lemma}
The preceding lemmas combined hence give us  a more complete picture of the \emph{information requirement} for accurate classification: it suffices that by the time a job $i$ departs from the system, there exists one label, $h$, whose cumulative log-likelihood ratio when compared against any other alternative label, $l$, is sufficiently large, i.e., $S_{i,t}(h,l)$ is at least $\ln(\szH/\delta)$ for all $l\neq h$ (Lemma \ref{lem:SijSufficient}), and this is essentially necessary (Lemma \ref{lem:EofShl}). 

More importantly, the above discussion reveals  a natural link from \emph{information need} to \emph{service requirement}. While in a traditional processing system a job may come with a certain \emph{size}, we can think of a job  in the information processing system as having a \emph{vector-valued} service requirement: for a job $i$ with true label $h$, if we interpret the quantity $S_{i,t}(h,l)$ as the amount of ``work'' already performed along the $l$th coordinate, then the job's service requirement would be, roughly speaking, that the work performed along \emph{all} coordinates, $l \in  \calH\backslash \{h\}$,  should surpass $\ln(\szH/\delta)$. 

This brings us to the last, and arguably most complex, sub-problem, $(c)$:  \emph{how} do we satisfy these service requirements in an efficient manner, and simultaneously for multiple jobs? Going back to the definition of $S_{i,t}(h,l)$ in Eq.~\eqref{eq:Sinsample}, we see that it can be written as the summation of the $Z_{i,j}(h,l,K_{i,j})$'s. Notably, if job $i$ has true label $h$, then $\E(Z_{i,1}(h,l,k)) = D(h,l,k)$, where $D(h,l,k)$ is the KL-divergence defined in Eq.~\eqref{eq:dklDef}. In other words, one inspection performed by an expert of type $k$ contributes, in expectation, $D(h,l,k)$ amount of ``work'' to $S_{i,t}(h,l)$. Viewing our processing task from this angle reveals a resemblance with a certain multi-class multi-server queueing system, where jobs come in different classes (and in this case, labels), and the amount of work a server can contribute to a job's service requirement during a unit time period depends both on the server's type and the type of the job being treated. 

Unfortunately, there remains a difficult yet fundamental obstacle that prevents us from directly applying our understanding of multi-class queueing systems to designing inspection policies. The above analogy makes it clear that some expert types may be more informative for a certain job label than others, because the values of $D(h,l,k)$ can vary across $h,l$, and $k$. Hence, to best harness the processing power of the experts and satisfy the service requirements across all jobs, the inspections should be arranged in a way that takes into account the types of experts performing the inspections \emph{and} the labels of jobs being inspected, for otherwise an expert could end up wasting her time working  on jobs that she has little expertise on. However, this brings us to a circular argument: while efficient inspection beckons a policy to be aware of job labels, we simply do \emph{not} know the job labels, for otherwise there would have not been a need to perform any inspection to begin with! 

\subsubsection{Overview of the Inspection Policy} 
\begin{figure}
\centering
\includegraphics[scale=.7]{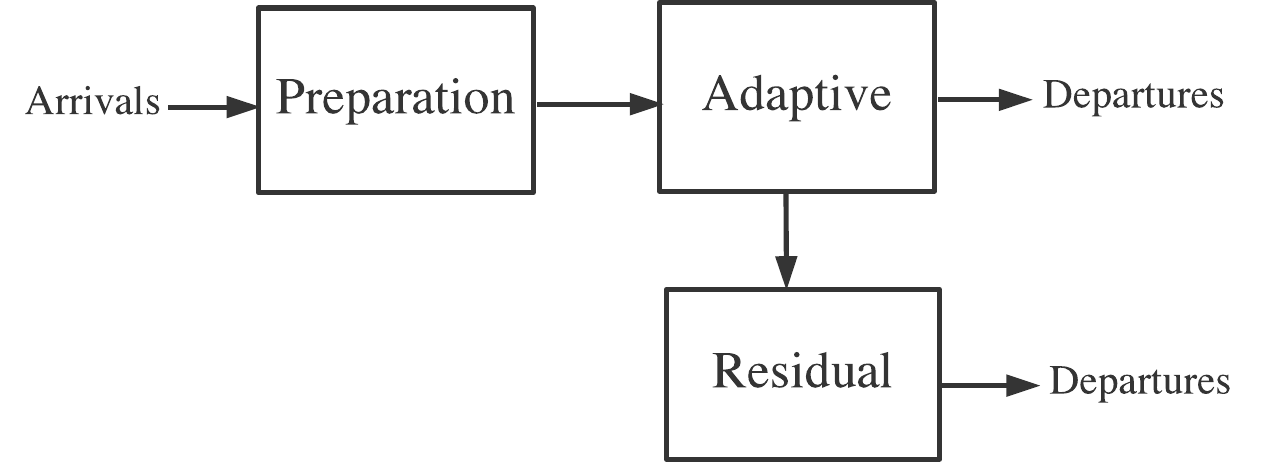}
\caption{Overall architecture associated with the proposed inspection policy.}
\label{fig:blockDiag}
\vspace{-10pt}
\end{figure}

Our inspection policy will make use of a three-stage architecture to circumvent the above-mentioned ``circular'' logic, illustrated in Figure \ref{fig:blockDiag}. 
\begin{enumerate}
\item In the first stage (Preparation), the policy ``boot-straps'' each incoming job, by inspecting it using randomly chosen experts with the goal of generating a \emph{coarse estimate} of its  true label. 

\item In the second stage (Adaptive), the policy performs the majority of the inspections and in an adaptive manner, with the main goal of {verifying} whether the coarse estimates are correct. Most of the jobs with a correct coarse estimate will be able {to} obtain an accurate classification by the end of the Adaptive stage and depart from the system. 

\item The third stage (Residual) treats those jobs whose coarse estimates were erroneous to ensure that they, too, will  receive an accurate classification.
\end{enumerate} 
We will show that $(1)$ the coarse estimates in the Preparation stage are sufficiently accurate so that little resource is wasted in the Adaptive stage, and $(2)$ the processing resources required in the Preparation and Residual stages amount to only a small fraction of the total resources. Together, they lead to the  resource efficiency of our inspection policy. 

We can also interpret the high-level structure of this three-stage architecture through a \emph{learning} versus \emph{verification} dichotomy: all jobs are first inspected by some ``generalists'' (i.e., random experts) to \emph{learn} a coarse label estimate. The system then enlists the ``specialists'' to \emph{verify} the validity of these estimates to a high accuracy. If a coarse estimate is deemed incorrect by the ``specialists'', the job is then sent \emph{back} to the ``generalists'' to perform learning thoroughly to reach an accurate estimate, albeit in a less efficient manner.

\subsection{Proof Outline}

We now provide a brief outline of the main steps of the proof. Expanding upon the informal discussion in the previous subsection, we formally describe in Section \ref{sec:DefOfPolicy} a prior-oblivious inspection policy that will be used to achieve the scaling in Theorem \ref{thm:main}.  In Section \ref{sec:lowerbound}, we build on Lemma \ref{lem:EofShl} and establish a lower bound on the number of experts that \emph{any} $\delta$-accurate policy must satisfy, which is expressed in terms of a solution to a certain linear optimization problem. This lower bound will serve as our performance benchmark of what an ``optimal'' inspection policy could achieve in terms of minimum system size.  In Sections \ref{sec:StabilityPrepAdaptive} and \ref{sec:StabilityResidual}, we establish a sufficient condition for the number of experts under which the proposed policy would stabilize all three stages. In particular, Section \ref{sec:StabilityPrepAdaptive}  contains the most technical portion of our proof, which relies on a fluid model to analyze the joint dynamics of the Preparation and Adaptive stages.  We complete the proof of Theorem \ref{thm:main} in Section \ref{sec:CompleteProof}, in two steps. We first show that the proposed policy is $\delta$-accurate, which, in light of Lemma \ref{lem:SijSufficient}, follows by construction in a straightforward manner. We then compare the sufficient condition on the number of experts, established in Sections \ref{sec:StabilityPrepAdaptive} and \ref{sec:StabilityResidual},  to the lower bound in Section \ref{sec:lowerbound}, and demonstrate that the ratio between the two converges to $1$ uniformly over all prior distributions for the job labels. This completes the proof of Theorem \ref{thm:main}. 


\section{Design of the Inspection Policy}
\label{sec:DefOfPolicy}

We present in this section the inspection policy that we will use to prove Theorem \ref{thm:main}. The main job-flow of the policy consists of three stages: Preparation, Adaptive, and Residual, as is illustrated in Figure \ref{fig:blockDiag}. We begin by explaining some basic actions of the experts. 

\subsection{Randomized Expert Visits}
\label{sec:expertVis}
We first introduce some notation. Define \vspace{-5pt}
\begin{equation}
r_k = \mu_k  \rho_k, \quad  k \in \calK,
\label{eq:rkdef}
\end{equation} 
and from the assumption in Eq.~\eqref{eq:avgRate}, we have that  $\sum_{k\in \calK} r_k = \overline{\mu} = 1. $. Because the lengths of inspections are exponentially distributed, ${r}_k$ is the probability that the next available expert is of type $k$ assuming that all experts are in BUSY in the present moment, and $mr_k$ is the average number of inspections that the pool of type-$k$  experts can complete in unit time. Let $d(h,l)$ be the average KL-divergence for when a job is inspected by an expert whose type is randomly drawn according to the distribution $\{{r}_k\}_{k\in \calK}$:
\begin{equation}
d(h,l) = \sum_{k \in \calK} D(h, l, k) r_k.
\end{equation}
 Denote by $\dav$ the minimum value among the $d(h,l)$,  $\dav = \min_{h, l \in \calH, h\neq l}d(h,l) $, and by $\olz$ the maximum log-likelihood ratio\footnote{The set $\calX$ being finite ensures that $\olz<\infty$. }
\begin{equation}
\olz  = \max_{{h, l \in \calH, \, k\in \calK}}\  \max_{x: p(l,k,x)\neq 0} \lt|\ln \frac{p(h,k,x)}{p(l,k,x)} \rt|. 
\label{eq:overlinez}
\end{equation}
Finally, define the constant
\begin{equation}
\zeta_0 =  \frac{8\olz^2 + 2 \dav }{\dav^2}. 
\end{equation}

\paragraph{Expert visits.} We say that an expert of type $k$ \emph{goes on a vacation} to mean that she starts processing a ``dummy'' job and remains in state  BUSY for a period of time that is exponentially distributed with mean $1/\mu_k$.   Suppose that an expert completes inspecting a previous job at time $t$, then she will choose to {\bf visit} one of the three stages, which means that the expert will either initiate an inspection for a job in that stage, or go on a vacation, depending on the inspection rules which will be specified in the next subsection. The choice of which stage to visit will be made by a simple randomized rule, independent of the rest of the system: the expert visits the Preparation, Adaptive, and Residual stage, with probability $q^P$, $q^A$, and $q^R$, respectively, where 
\begin{align}
\label{eq:qPdef}
q^P = & \frac{\zeta_0 \ln\ln(1/\delta)+\ln^{-1}(1/\delta)}{m}, \quad \quad q^R = \frac{3\szH\zeta_0(1+\ln(4\szH) \ln^{-1}(1/\delta)) +1 }{m},\\
q^A = & 1-q^P-q^R, 
\label{eq:qAdef}
\end{align}
and we assume that the policy will be applied under a range of parameters where all expressions above lie in the interval $(0,1)$.  

\subsection{Multi-stage Inspection Policy}

We now describe our inspection policy in detail, where the exposition for each stage is broken down into three parts: $(1)$ Workload: how many inspections need be completed on a job in each stage;  $(2)$ Departure rules: where the job goes next; $(3)$ Expert actions: how experts perform inspections. 


\subsubsection{Preparation Stage} Every job that arrives to the system will first be processed in the Preparation stage. The objective is to produce a coarse estimate of the job's label using only a \emph{small} number of inspections, performed by random experts whose types are drawn according to $\{{r}_k\}_{k\in \calK}$. The randomization in the expert types ensures that information is gained about the job's true label at a non-zero rate.  The coarse label estimate will then be used to ``bootstrap'' processing in the Adaptive stage to further enhance the classification accuracy. 

\paragraph{Workload.} 

Every job will receive $n^P$ inspections in the Preparation stage, where
\begin{equation}
n^P =  \zeta_0 \ln\ln(1/\delta). 
\end{equation}

\paragraph{Departure rules.} When a job has received the \emph{outcomes} from all $n^P$ inspections, it departs from the Preparation stage and enters the Adaptive stage. 

\paragraph{Expert actions.} Denote by $\buw_0(t)$ the total number of \emph{uninitiated} inspections in the Preparation stage at time $t$. An expert who visits the Preparation stage at time $t$ will attempt to initiate an inspection for a job in the Preparation stage in a first-come-first-serve fashion. If $\buw_0(t)=0$, then the expert goes on a vacation.

\subsubsection{Adaptive Stage} 

The Adaptive stage is the ``power-house'' of the system that performs the majority of all inspections. Its defining feature is that the number of inspections that a job receives from each expert type will be decided \emph{adaptively} depending both on the coarse label estimate from the Preparation stage, and on the existing aggregate workload in the Adaptive stage. 
The main objective of this stage is to \emph{verify} the correctness of the coarse label estimates: \emph{most} jobs with a correct coarse estimate depart from the system after the Adaptive stage, while those with incorrect coarse estimates are likely to be sent to the Residual stage for further processing. 

\paragraph{Workload.} The workload generation in this stage is more complex than that of the Preparation stage. Upon arriving to the Adaptive stage, a job, $i$, is assigned a \emph{workload vector} $\{\uLam_{i,k}\}_{k\in \calK}$, where $\uLam_{i,k}$ is the number of inspections to be performed by experts of type $k$ on job $i$ during its stay in the Adaptive stage. We will denote by $\overline{\uLam}_{i,k}(t)$ the remaining number of inspections by experts of type $k$ at time $t$, defined by the difference between $\uLam_{i,k}$ and the number of inspections already \emph{initiated} by experts of type $k$ for job $i$ by time $t$. Denote by $Q^A(t)$ the set of jobs in the Adaptive stage at time $t$. We define the \emph{workload at expert pool $k$} as 
\begin{equation}
\buw_k(t) = \sum_{i \in Q^A(t)} \overline{\uLam}_{i,k}(t), \quad t\in \rp.
\label{eq:workKthexpool} 
\end{equation}
We will refer to $\buw(\cdot) = \{\buw_k(\cdot)\}_{k\in \calK}$ as the workload process. 

We now explain how the workload vectors, $\{\uLam_{i,k}\}_{k\in \calK}$, are generated. Denote by $\what{H}_i^P$ the maximum likelihood estimator of job $i$, $\what{H}_{i,t}$ (Eq.~\eqref{eq:MLdef}), at the time it exits the  Preparation stage, and suppose that $\what{H}_{i}^P = h$. Let $\{n_{h,k}\}_{k\in \calK}$ be an optimal solution to the following linear optimization problem:
\begin{align}
\label{eq:argminLinearWeight}
\mbox{minimize} \quad & \sum_{k \in \calK} n_k\buw_k(t),  \\
\label{eq:calNh1}
\mbox{s.t.} \quad \sum_{k\in \calK}D(h,l,k)n_k \geq&  \ln(2\szH/\delta)+g_\delta, \quad \forall l \in \calH\bslh, \\
n_k \geq&  0 , \quad \forall k \in \calK, \\
\sum_{k \in \calK} n_k \leq& v_\delta,  \label{eq:calNh3}
\end{align}
with ties broken arbitrarily, where $g_\delta$ and $v_\delta$ are two auxiliary constants that do not depend on $h$: 
\begin{align}
\label{eq:gdeltadef}
g_\delta  =& 3\olz\uld^{-1/2} \sqrt{ \ln(1/\delta)\ln\ln(1/\delta)} \\
v_\delta =& 2\uld^{-1}\ln(1/\delta) \lt[1+({\ln(2\szH) + g_\delta}) \ln^{-1}(1/\delta) \rt]. 
\label{eq:vdeltadef}
\end{align}
We will denote by $\calN_{h}$ the set of all vectors $\{n_k \}_{k\in \calK}$ which satisfy the constraints in Eqs.~\eqref{eq:calNh1} through \eqref{eq:calNh3}. One can verify that the above optimization problem always admits a feasible solution. 
Finally, the workload vector for job $i$ will be obtained by rounding down $\{n_{h,k}\}_{k\in \calK}$: 
\begin{equation}
\uLam_{i,k} = \lfloor n_{h,k} \rfloor , \quad \forall k \in \calK. 
\label{eq:Lamb_i_k}
\end{equation}

\emph{Interpretation of workload:} The workload vector captures the combinations of inspections job $i$ should receive \emph{assuming} that its true label is indeed $\what{H}^P_i$, in which case Eq.~\eqref{eq:calNh1} ensures that the cumulative log-likelihoood ratios are sufficiently large to make an accurate classification. As mentioned earlier, the optimization in Eq.~\eqref{eq:argminLinearWeight} is reminiscent of the family of max-weight scheduling policies (cf.~\cite{tassiulas1992stability}): our policy aims to create workload in order to minimize an inner-product between the new workload and the existing inspections. The proof in subsequent sections will demonstrate that this adaptive procedure allows the inspection policy to work well with any prior distribution, hence making the policy prior-oblivious. 

\paragraph{Departure rules.} A job $i$ departs from the Adaptive stage as soon as it has received the outcomes of all $\uLam_{i,k}$ inspections from experts of type $k$, for all $k\in \calK$. Suppose that the departure occurs at time $t$. The policy then executes the following decision: 
\begin{enumerate}
\item If there exists $h\in \calH$, such that 
\begin{equation}
S_{i,t}(h,l) \geq \ln(2\szH/\delta), \quad \forall l \in \calH, l \neq h, 
\label{eq:adaptiveExitCriteria}
\end{equation}
then job $i$ departs from the system, and a classification is produced by setting $\what{H}_i = \what{H}_{i, t} = h.$
\item Otherwise, job $i$ enters the Residual stage. 
\end{enumerate}

\paragraph{Expert actions.} Suppose that an expert of type $k$ visits the Adaptive stage at time $t \in \rp$. If the workload for the $k$th expert pool, $\buw_k(t)$, is non-zero (Eq.~\eqref{eq:workKthexpool}), then the expert initiates an inspection for a job associated with one unit of work in $\buw_k(t)$, in a first-come-first-serve fashion. If $\buw_k(t)=0$, then the expert goes on a vacation. 

\subsubsection{Residual Stage} The Residual stage acts as a ``clearing house" that treats those jobs whose inspections in the Adaptive stage failed to produce an accurate classification. Similar to the Preparation stage, jobs are inspected by random experts, but they receive significantly more inspections in the Residual stage in order to produce a highly accurate label classification.

\paragraph{Workload.}  The moment a job enters the Residual stage, all of its previous inspections and outcomes are \emph{discarded}. Similar to the Preparation stage, each job will receive a fixed number of inspections,  $n^R$, where 
\begin{equation}
n^R = \zeta_0 \ln(4\szH/\delta).
\end{equation}

\paragraph{Departure rules.} A job $i$ in the Residual stage departs from the system  as soon as it has received the results from all $n^R$ inspections, and a classification of job $i$'s type is produced by setting $ \what{H}_i = \what{H}_{i, t}$.  Note that the classification is made solely based on the inspections in the Residual stage, as we have discarded all inspections from the earlier stages. 

\paragraph{Expert actions.} An expert who visits the Residual stage will attempt to initiate an inspection for a job in the Residual stage in a first-come-first-serve fashion. If there is no job currently in the Residual stage, or if all jobs in the Residual stage  have all of their $n^R$ inspections already initiated, then the expert goes on a vacation.  This concludes the description of our inspection policy. 

\section{Lower Bound on Optimal System Size}
\label{sec:lowerbound}

We establish in this section a fundamental lower bound on the minimum  number of experts required in order for the system to be stable that holds for \emph{any} $\delta$-accurate policy. The following Fundamental Linear Program is central to this lower bound as well as our subsequent analysis.

\begin{definition}
\label{def:FLP}
The {\bf Fundamental Linear Program}, denoted by \flp, is defined as follows.
\begin{align}
\mbox{minimize}  &   \quad\quad\quad m \\
\label{eq:Lp1repconstr1}
\mbox{s.t.} & \quad   \sum_{h\in \calH} n_{h,k}\pi_h \leq r_k m, \quad k \in \calK, \\
\label{eq:Lp1repconstr2}
& \quad \sum_{k \in \calK}n_{h,k} D(h,l,k) \geq \ln (1/\delta),\quad \forall h, l \in\calH, h\neq l,\\
& \quad n_{h,k} \geq 0, \quad \forall h\in \calH, k\in \calK,
\label{eq:FLP}
\end{align}
\end{definition}
where $D(h,l,k)$ and $r_k$ are defined in Eqs.~\eqref{eq:dklDef} and \eqref{eq:rkdef}, respectively. 

We provide some intuition to motivate the above definition. Recall from Lemma \ref{lem:EofShl} that in order for any policy to be $\delta$-accurate, for a job $i$ with label $h$, the expected value of the cumulative log-likelihood ratio $S_{i,t}(h,l)$ should be at least $\ln(1/\delta)$ by the time  job  $i$ departs from the system. Furthermore, an inspection by an expert of type $k$ increases the value of  $S_{i,t}(h,l)$ by $D(h,l,k)$ in expectation. If we interpret the variables, $n_{h,k}$, as the number of inspections a job with true label $h$ should receive from an expert of type $k$, then Eq.~\eqref{eq:Lp1repconstr1}  in \flp corresponds to the resource constraint of having $m \rho_k$ type-$k$ experts, each of whom can  perform $\mu_k$ inspections per unit time, and  Eq.~\eqref{eq:Lp1repconstr2} to the above-mentioned constraint on $\E(S_{i,j}(h,l))$ imposed by Lemma \ref{lem:EofShl}. Therefore, \flp captures the problem of finding minimal system size faced by a decision maker who already \emph{knows} the true labels of the jobs {and is only interested in ``verifying'' them in order to satisfy the condition of Lemma \ref{lem:EofShl}}, and we would expect the optimal value of \flp to be a lower bound for what is achievable in our problem, where the job labels are unknown.


The following proposition is the main result of this subsection, which states that the minimum system size under any $\delta$-accurate policy is essentially no smaller than the solution to \flp, as $\delta\to 0$. The proof builds upon the lower bound in Lemma \ref{lem:EofShl}, and is given in Appendix \ref{app:prop:mstarLWB}. 

\begin{proposition} 
\label{prop:mstarLWB}
Fix $\delta \in (0,1)$ and $\pi$. Denote by $m^*_F$ the optimal value of \flp. We have that
\begin{equation}
m^*(\delta, \pi) \geq b_\delta  \cdot m^*_F,
\end{equation}
where $b_\delta = (1 - \delta)\lt[ 1 - \lt(\ln\frac{1}{1-\delta}+e^{-1} \rt) \ln^{-1}(1/\delta) \rt]$. In particular, $b_\delta \uparrow 1$ as $\delta\downarrow 0$. 
\end{proposition}

The next lemma states some useful properties of \flp. The proof is given in Appendix \ref{app:lem:lp1structure}. 
\begin{lemma}
\label{lem:lp1structure}  
There exists an optimal solution of \flp, $(m^*_F, \{n^*_{h,k}\})$, such that the following holds: 
\begin{align}
\label{eq:sumnhk1}
\old^{-1}\leq & \frac{1}{\ln(1/\delta) }\sum_{k \in \calK}n^*_{h,k} \leq \,  \uld^{-1} , \quad \forall h \in \calH, \\
m^*_F \geq&  \, \old^{-1}   \ln (1/\delta). 
\label{eq:mstar}
\end{align}
\end{lemma}

\section{Stability of Preparation and Adaptive Stages}
\label{sec:StabilityPrepAdaptive}

We  establish in this section a sufficient condition on the number of experts, $m$, in order for the Preparation and Adaptive stages to be jointly stable under the proposed inspection policy. An analogous result for the Residual stage will be established in a subsequent section. We denote by $Q^P(t)$, $Q^A(t)$ and $Q^R(t)$ the number of jobs in the Preparation, Adaptive and Residual stages, respectively, at time $t$. A stage being stable  means that the process $Q^\cdot(\cdot)$ for that stage is positive recurrent. The main result of this section is the following theorem.
\begin{theorem}
\label{thm:prepAndAdapStable}
Define $\ulr = \min_{k \in \calK}r_k = \min_{k \in \calK} \rho_k\mu_k$. The Preparation and Adaptive stages are stable whenever $ m q^P > n^P = \zeta_0 \ln\ln(1/\delta)$, and 
\begin{equation}
 m q^A > \lt(1+\frac{\ln(2\szH)+ g_\delta}{\ln(1/\delta)}\rt)\lt(1+  \frac{2\szH^2  \old }{\uld\, \ulr \ln(1/\delta) }\rt) m_F^*, 
 \end{equation} 
where $m_F^*$ is the optimal value of the Fundamental Linear Program in Definition \ref{def:FLP}. 
\end{theorem}

{\bf Proof Overview for Theorem \ref{thm:prepAndAdapStable}.} The remainder of this section is devoted to the proof of Theorem \ref{thm:prepAndAdapStable}, and, unless stated otherwise, we will use the word ``system''  to refer to the Preparation and Adaptive stages only. We begin by giving an overview of the proof and highlighting some of the main technical challenges that motivate our approach. 

Let us first recall some high-level features of the system dynamics. The expert actions are fairly simple in both stages by simply trying to initiate an inspection in a non-adaptive manner. Creating inspection workloads is also straightforward for the Preparation stage where each job has a fixed number of inspections. The main complexity therefore lies in how the vector-valued workloads are created in the Adaptive stage, which depends both on the job's type estimate from the Preparation stage, and the {aggregate} workloads in the Adaptive stage.  Given the disparity of complexity, a natural approach would be to treat the two stages separately:  the Preparation stage admits simpler dynamics and is easy to analyze while the Adaptive could be tackled using the Foster-Lyapunov criterion. Unfortunately, this approach falls short  because the processing in the Preparation stage destroys the memoryless property of the initial arrival process, rendering its output process non-Markovian. Therefore, the Adaptive stage cannot be treated as an isolated Markov process without taking into account the state of the Preparation stage as well.

 To overcome this problem, we will model the dynamics in both stages \emph{jointly}, and formulate a set of \emph{fluid solutions}, expressed as solutions to a system of ordinary differential equations (ODE), to capture the essential dynamics of this joint process. Specifically, following a general program developed by \cite{rybko1992ergodicity,dai1995positive}, we will show that, under proper scaling, the process of system workloads converges almost surely to a set of {fluid solutions}. We then show that if the number of experts satisfies the conditions stated in Theorem \ref{thm:prepAndAdapStable}, then the fluid solutions exhibit a certain global contraction property with respect to a one-homogeneous Lyapunov function. These two properties will then be used to show that the original workload process is positive recurrent, which in turn implies the stability of the system.  The proof will be carried out in the following main steps: 
\begin{enumerate}
\item We define in Section \ref{sec:stateRep} the Markov process that captures the system dynamics, as well as the fluid solutions. 
\item We show in Section \ref{sec:convgToFluid} (Proposition \ref{prop:convrgToFluid}) that the workload process converges to a fluid solution almost surely  over an appropriately defined probability space. 
\item We establish the global contraction property of the fluid solutions in Section \ref{sec:driftFluid} (Proposition \ref{prop:drainFluid}) by showing a multiplicative decay in a quadratic Lyapunov function along the trajectory of \emph{any} fluid solution. 
\item Finally, we complete the proof in Section \ref{sec:completePStabiPrepandAdap} by combining Propositions \ref{prop:convrgToFluid} (convergence) and \ref{prop:drainFluid} (contraction) with a variant of the Foster-Lyapunov criterion (Theorem 8.13, \cite{robert2003}). 
\end{enumerate}

There are two main technical challenges which we develop novel tools to overcome: $(1)$ the proof requires characterizing the limit points of the solutions to the  linear optimization sub-routine  in Eq.~\eqref{eq:argminLinearWeight} under fluid scaling, which is difficult because the parameters in the objective function, $\{\buw_k(\cdot)\}_{k \in \calK}$, are themselves stochastic variables. We will employ a careful analysis of the (semi-)continuity properties of the optimization problem to study these limit points; $(2)$ the jobs' transitions from the Preparation stage to the Adaptive stage are complex because a job can depart from the Preparation stage only after all of its inspection have {completed} (for otherwise the inspection outcomes would not have  been available). This in turn causes the order of departures from a stage to deviate {from} that of the arrivals, making standard techniques ineffective in showing the convergence of the system's stochastic trajectory to a fluid limit. We will prove convergence by developing explicit bounds on the degree of ``shuffling'', which allows us to conclude that the deviation is not too substantial to invalidate convergence.

\subsubsection{Additional Notation} \label{sec:addiNot}
For a vector $x=(x_1, \ldots, x_n)$, we will denote by $\|x\|_2$ the $l_2$ norm of $x$: $\|x\|_2 = \sqrt{\sum_{i=1}^nx_i^2}$. We will denote by $\snorm{\cdot}$ the maximal norm of a function over the interval $[0,T]$: $\snorm{f(t)} = \sup_{t\in [0,T]} |f(t)|.$  For $x,y\in \R$, we will use $x \vee y$ and $x \wedge y$ to denote $\min\{x,y\}$ and $\max\{x , y\}$, respectively.  For two vectors of the same dimension, $\blx$ and $\bly$, we write $\blx \leq \bly$ if all coordinates of $\blx$ are dominated by those of $\bly$. Similarly, for a set of vectors, $\calY$, we write $\blx \leq \calY$ if $\blx \leq \bly$ for all $\bly \in \calY$. The addition of two sets $\calX+\calY$  is defined to be the set $\{z: \, z = x+y, \, x\in \calX, y\in \calY\}$. 

\subsection{Classification Accuracy of Preparation Stage}
\label{sec:accuPrepStage}

We begin in this subsection with a result that  bounds the classification error of the coarse label estimate from the Preparation stage. Denote by $\what{H}_i^P$ the maximum likelihood estimator, $\what{H}_{i,t}$, when job $i$ exits the Preparation stage. Recall from the construction of our policy that each job will be inspected for the same, deterministic number of times in the Preparation stage, and it is not difficult to verify that $\{\what{H}_i^P\}_{i \in \N}$ are i.i.d. We will denote by $\pi^P$ the distribution of the estimator $\what{H}_1^P$, 
\begin{equation}
\pi^P_h \bydef \pb\lt( \what{H}_1^P  = h\rt), \quad h \in \calH. 
\label{eq:piPdef}
\end{equation}
and by $\epsilon^P$ as the error probability $\epsilon^P  \bydef  \max_{ h \in \calH} \pb(\what{H}_1^P \neq h \bbar  {H}_1 =h)$. 
The following proposition provides an upper bound on $\epsilon^P$, which in turn upper-bounds the distance from $\pi^P$ to the original prior distribution, $\pi$. The proof is given in Appendix \ref{app:prop:erroAfterPrep}, which relies on a generalization of Lemma \ref{lem:SijSufficient} in combination with fact that $S_{i,t}(h,l)$ can be viewed as a martingale under proper conditioning. 

\begin{proposition}
\label{prop:erroAfterPrep} We have that $\epsilon^P  \leq 2  \szH\ln^{-1}(1/\delta)$, and
\begin{equation}
\pi^P_h \leq \pi_h + \epsilon^P \sum_{h' \neq h} \pi_{h'} \leq \pi_h + 2  \szH\ln^{-1}(1/\delta), \quad \forall h \in \calH. 
\label{eq:piP}
\end{equation}
\end{proposition}

\subsection{State Representation and Fluid Solutions}
\label{sec:stateRep}
We describe in this subsection a Markovian representation of the Preparation and Adaptive stages as well as the notion of fluid solutions. We will index jobs in the two stages according to the order in which they arrive. Denote by $I(t)$ the total number of jobs in the system at time $t$, and let $\calI(t) = \{1,2,\ldots, I(t)\}$. For instance, job $1$ corresponds to the oldest job in system at time $t$, and job $I(t)$ corresponds to the youngest. The indices will be updated accordingly in the event of the departure of job $i$, where all jobs with an index greater than $i$ will have their index reduced by $1$.  For each $i\in \calI(t)$, we define a \emph{job state}, $Y_i(t)$, which consists of the following variables: 

\begin{enumerate}
\item $Y^S_i(t) \in \{1,2\}$ represents the stage the job is currently, with $1$ and $2$ corresponds to job $i$ being in the Adaptive and Preparation stage, respectively. 
\item {$\{(X_{i,s},E_{i,s} )\}_{s=1}^{N_{i,t}}$ }contains all the past inspection responses of job $i$, along with the corresponding expert types, where $N_{i,t}$ the number of inspections received by job $i$ by time $t$. 
\item $Y^E_i(t) \subset \calE$ is the set of experts who are in the process of inspecting job $i$ at time $t$. 
\item $N_{i}^P(t)$ is the number of remaining inspections that job $i$ has left in the Preparation stage. 
\item $\overline{\uLam}_{i,k}(t)$ is job $i$'s remaining number of inspections to be completed by experts of type $k$ in the Adaptive stage. 
\end{enumerate}

{The} above variables completely {specify} the state of the Preparation and Adaptive stages at time $t$, and it is not difficult to verify that $\{\calI(t),\{Y_i(t)\}_{i\in \calI(t)}\}_{t \in \rp}$ is a countable-state Markov process. 
We now formally define the \emph{workload process} using the above state representation. Denote by $\buw_0(t)$  the total number of remaining uninitiated  inspections in the Preparation stage
\begin{equation}
\buw_0(t) = \sum_{i \in \calI(t), Y_i^S(t)=1} N^P_i(t). 
\end{equation} 
Similarly, the \emph{workload in the $k$th expert pool} in the Adaptive stage, is defined by: 
\begin{equation}
\buw_k(t) = \sum_{i \in \calI(t), Y_i^S(t)=2} \overline{\uLam}_{i,k}(t). 
\end{equation}

The following notion of fluid solutions will serve as an approximation to the workload process. 
\begin{definition}[Fluid Solutions]
\label{def:fluidSol}
The functions $\blw_k: \rp\to \rp, k\in \{0\}\cup \calK$ are called a \emph{fluid solution} if there exist {Lipschitz-continuous functions} $\bla_k, \bld_k: \rp\to \rp$, $k\in \{0\}\cup \calK$, with $\bla_k(0)=\bld_k(0)=0$ {with} Lipschitz constant $c_\lcal>0$ such that 
\begin{equation}
\blw_k (t) = \blw_k(0)+ \bla_k(t)-\bld_k(t), 
\end{equation}
where, for almost all $t\in \rp$, 
\begin{align}
& \dot{\bla}_0(t) = n^p, \quad \quad \quad  \dot \bld_0(t)  =  \left\{ \begin{array}{ll}
          m q^P, & \quad \mbox{if } \blw_0(t)>0,\\
          n^P, & \quad  \mbox{if } \blw_0(t)=0,\\
         \end{array} \right. \nln 
         \nln
& \lt\{ \dot{\bla}_k(t) \rt\}_{k \in \calK} \leq  (1+\ln^{-1}(1/\delta)) \sum_{h \in \calH}  \pi_h^P \calN^*_h (\blw(t)), \nln
& \dot \bld_k(t)  =  \left\{ \begin{array}{ll}
           r_k m q^A, & \quad \mbox{if } \blw_k(t)>0,\\
          \dot{\bla}_k(t) , & \quad  \mbox{if } \blw_k(t)=0,\\
         \end{array} \right.  \quad \forall k \in 1, \ldots,  \szK,
\end{align}
where $q^P, q^A$ and $n^P$ were defined in Section \ref{sec:DefOfPolicy}, and $\pi^P_h$ in Eq.~\eqref{eq:piPdef}.   $\calN^*_{h}(\blw(t))$ is the set of  optimal solutions for the optimization problem 
\begin{equation}
\min_{\{n_k\}\in \calN_{h} } \, \sum_{k \in \calK} n_k \blw_k(t), 
\label{eq:Nstarw(t)}
\end{equation}
and the set $\calN_h$ was defined in Eqs.~\eqref{eq:calNh1} through \eqref{eq:calNh3}.  Fix $\blw^0\in \rp^{ \szK+1}$. We denote by $\calW(\blw^0)$  the set of all fluid solutions with the initial condition $\blw(0)=\blw^0$. 
\end{definition}

\subsection{Convergence of Stochastic Sample Paths to the Fluid Solutions}
\label{sec:convgToFluid}

We show in this section that the workload process converges to a fluid solution under proper scaling. We will consider a sequence of systems, indexed by $s\in \N$, which have different initial conditions for the workload process at $t=0$ but are otherwise identical. We will denote by $\buw^s(\cdot)$ the workload process $\buw(\cdot)$ in the $s$th system.  For $n\in \N$, and a process $\{X(t)\}_{t\in \rp}$, we will use $\{X(n,t)\}_{t\in \rp}$ to denote the \emph{normalized process}:
\begin{equation}
X(n,t) = \frac{1}{n}X(nt), \quad t\in \rp. 
\label{eq:normalizedprocess}
\end{equation}
The following proposition is the main result of this subsection. 

\begin{proposition}
\label{prop:convrgToFluid}
Fix $\blw^0\in \rp^{ \szK+1}$. Consider a sequence of initial conditions $\{\blw^{(s)}\}_{s \in \N}$, such that for a sequence of positive numbers $\{z_s\}_{s \in \N}$ with $\lim_{s \to \infty} z_s = \infty$, we have that $\lim_{s \to \infty} z_s^{-1}\blw^{(s)} = \blw^0$. 

Suppose $\buw^s(0) = \blw^{(s)}$, for all $s\in \N$. Then,  for all $T >0$,  the following convergence takes place: 
\begin{equation}
\lim_{s \to \infty} \inf_{\blw \in \calW(\blw^0)} \snorm{\buw^{s}(z_s, t) - \blw(t)} = 0,\quad \mbox{almost surely}. 
\end{equation}
\end{proposition}

The remainder of this sub-section is devoted to the proof Proposition \ref{prop:convrgToFluid}. Because the system dynamics is quite complex, we will prove convergence using a sample-path-based approach that helps us isolate the probabilistic aspect of the dynamics from its deterministic counterpart (Similar sample-path-based approaches have also been used in \cite{bramson1998state,tsitsiklis2012power}). In particular, we first identify a large subset of the sample space,  called the \emph{regular set}, which contains the sample paths that exhibit certain typical behaviors. We  then show that convergence to the fluid solutions occurs over \emph{every} sample path in the regular set.  

\subsubsection{Regular Set} We now construct the regular set. We will define all random quantities of the system on the same probability space, $(\Omega, \calF, \pb)$. 

\begin{definition}
\label{def:regEvents} 
Fix $T>0$. We define the following elements of $\calF$. 
\begin{enumerate}
\item $\calC_0$:  Denote by $\Xi_0(t)$ the number of jobs that have arrived to the Preparation stage in the interval $[0,t]$. Define $\calC_0$ as the event where
\begin{equation}
\lim_{z \to \infty} \snorm{\Xi_0(z,t)  -  t } = 0, 
\label{eq:regXi}
\end{equation}

\item$\calC_\calE$: Denote by $R_{e, P}(t)$ and $R_{e, A}(t)$ the total number of times an expert $e\in \calE$ visits the Preparation and Adaptive stage, respectively,  during the interval $[0,t]$. Let $k_e$ be the type of expert $e$. Define $\calC_\calE$ as the event where 
\begin{align}
\lim_{z \to \infty} \max_{e\in \calE} \snorm{R_{e,P}(z,t)  -  q^P \mu_{k_e}t } = 0, \quad   \mbox{ and } 
\quad  \lim_{z \to \infty} \max_{e\in \calE} \snorm{R_{e,A}(z,t)  -  q^A \mu_{k_e} t } =& 0. 
\label{eq:regRAP}
 \end{align}

\item $\calC_H$: Denote by $\what{H}^P_i$ the ML estimator of the type of the $i$th job upon leaving the Preparation stage. Let $g_h (t)  = \sum_{i=1}^{\lfloor t \rfloor} \mathbb{I}(\what{H}^P_i = h)$.  Define $\calC_H$ as the event where
\begin{equation}
\lim_{z \to \infty} \max_{h \in \calH} \snorm{g_h(z,t) - \pi^P_h t} = 0. 
\label{eq:regHhatP}
\end{equation}

\end{enumerate}
\end{definition}

Define the \emph{regular set}, $\calC$, as the intersection of all three events in Definition \ref{def:regEvents}:  $\calC = \calC_0\cap \calC_\calE\cap \calC_H$.  We have the following useful lemma. The proof is a direct consequence of the (functional) law of large numbers applied to each of the three events, and is omitted. 
\begin{lemma}
\label{lem:CisProb1} 
Fix $T>0$. We have that $\pb(\calC_0)= \pb(\calC_\calE) = \pb(\calC_H)=1$, and, consequently, $\pb(\calC)=1. $ \end{lemma}

\subsubsection{Proof of Proposition \ref{prop:convrgToFluid}} We return to the proof Proposition \ref{prop:convrgToFluid}, which will be completed  in two parts. Fix any sample path $\omega\in \calC$. In the first part, we will show that over any finite interval, any sub-sequence of $\{\buw^{z_s}(\cdot)\}_{s\in \N}$ admits a further converging sub-sequence that  converges coordinate-wise to a  Lipschitz-continuous function. We will then show, in the second part, that all such limiting functions are in fact fluid solutions, i.e., for almost all $t$, their derivatives satisfy the conditions given in Definition \ref{def:fluidSol}.  The next result summarizes the first part of the proof. Fix $k \in \{0, 1, \ldots,  \szK\}$. We will write
\begin{equation}
\buw^s_k(t) = \buw^s_k(0)+ \bua^s_k(t)- \budta^s_k(t), \quad  t >0, 
\label{eq:Ws0v1}
\end{equation}
where $\bua^s_k(t)$ and $\budta^s_k(t)$ denote the total number of inspections associated with the arrivals, and the number of inspections that have been initiated, during $[0,t]$, respectively, for jobs associated with the workload $\buw^s_k(\cdot)$. The proof of the following proposition is given in Appendix \ref{app:prop:subseqConvrg}.

\begin{proposition}
\label{prop:subseqConvrg}
Fix $T >0$. Denote by $\calL_c$ the set of coordinate-wise $c$-Lipschitz functions from $[0,T]$ to $\rp^{3( \szK+1)}$. Fix the sample, $\omega\in \calC$, and an increasing sequence, $\{s_i\}_{i\in \N}$, and let $c=[q^A+(v_\delta+1)q^P]m+ n^P$, where $v_\delta$ was defined in Eq.~\eqref{eq:vdeltadef}. Then, there exists  $(\blw, \bla, \bld) \in \calL_c$, and  an increasing sequence, $\{i_j\}_{j\in \N}\subset \N$, such that
\begin{equation}
\lim_{j\to \infty}  \snorm{(\buw, \bua, \budta) ^{s_{i_j}}(z_{s_{i_j}}, t) - (\blw, \bla, \bld)(t)}= 0. 
\end{equation}
We will refer to these $\blw(\cdot)$'s as limit points of $\{\buw^{s}(z_{s_i}, \cdot)\}_{i \in \N}$. 
\end{proposition}

The next result states that \emph{all} of the  limit points in Proposition \ref{prop:subseqConvrg} are in fact fluid solutions. 

\begin{proposition}
\label{prop:limpointarefluid} Let $\blw(\cdot)$ be a limit point as defined in Proposition \ref{prop:subseqConvrg}. Then $\blw(\cdot)$ is also a fluid solution, as defined by Definition \ref{def:fluidSol}. 
\end{proposition}

\bpf Fix $\omega\in \calC$, and a limit point, $(\blw, \bla, \bld)$, with the corresponding sequence $\{i_j\}_{j\in \N}$, as defined in Proposition \ref{prop:subseqConvrg}. To avoid excessive use of subscripts, we will use $\bar{s}_j$ and $\bar{z}_j$ in place of $s_{i_j}$ and $z_{s_{i_j}}$, respectively.  Fix $t\in (0,T)$ to be a time where all coordinates of the limit point are differentiable.  We begin with the Preparation stage, with $k=0$. As was mentioned in the proof of Proposition \ref{prop:subseqConvrg}, each new job that arrives to the Preparation stage creates $n^P$ inspections. We have that
\begin{equation}
\lim_{z\to \infty} \snorm{\bua_0(z,t)-n^Pt} = \lim_{z\to \infty} n^P \snorm {\Xi_0(z,t)-t} = 0,
\end{equation}
where the second equality follows from Eq.~\eqref{eq:regXi}. This shows that  $\dot{\bla}_0(t) = n^P$. 

For $\bld_0(t)$, we consider two cases depending on the value of $\blw_0(t)$. First, suppose that $\blw_0(t)>0$. Then there exists $\overline{\epsilon}>0$, such that for all $\epsilon \in (0,\overline{\epsilon})$, there exists $N_\epsilon>0$ such that, for all $j\geq N_\epsilon$, 
\begin{equation}
\buw^{\bar{s}_j}(t' \bar{z}_j) >0, \quad \forall t' \in [t, t+\epsilon]. 
\label{eq:buwpositive}
\end{equation}
Fix $j$ and $t'$ so that Eq.~\eqref{eq:buwpositive} is true. Because the workload is non-zero, any expert who visits the Preparation stage at time $t'\bar{z}_j$ will necessarily lead to a unit increment in the process $\budta^{\bar{s}_j}_0(\cdot)$. We thus have that, for all $j\geq N_\epsilon$, 
 \begin{equation}
\budta^{\bar{s}_j}_0((t + \epsilon) \bar{z}_j)-\budta^{\bar{s}_j}_0(t \bar{z}_j) = \sum_{e\in \calE}  R_{e,P}( (t+ \epsilon) \bar{z}_j)-R_{e,P}(t  \bar{z}_j), 
\end{equation}
which, with scaling, implies that
\begin{equation}
\lim_{j\to \infty} (\budta^{\bar{s}_j}_0(\bar{z}_j, t + \epsilon)-\budta^{\bar{s}_j}_0(\bar{z}_j, t) ) = \lim_{j\to \infty} \lt( \sum_{e\in \calE}  R_{e,P}( \bar{z}_j, t+ \epsilon)-R_{e,P}( \bar{z}_j, t )  \rt) = q^P m \epsilon, 
\end{equation}
where the last step follows from Eq.~\eqref{eq:regRAP}, and the fact that $\sum_k {\rho_k \mu_k} = \overline{\mu} = 1$ (Eq.~\eqref{eq:avgRate}). Taking the limit as $\epsilon\to 0$, we obtain
\begin{align}
\dot \bld_0(t) = & \lim_{\epsilon \to 0} \epsilon^{-1} \lt[ \lim_{j\to \infty} (\budta^{\bar{s}_j}_0(\bar{z}_j, t + \epsilon)-\budta^{\bar{s}_j}_0(\bar{z}_j, t) ) \rt]  =  \lim_{\epsilon \to 0} \epsilon^{-1} (q^P m \epsilon) =  q^Pm, \quad \mbox{if } \blw_0(t)>0. 
\label{eq:dotbld0_1}
\end{align}
Now suppose that $\blw_0(t)=0$. In this case, we will exploit the properties that $\blw_0(\cdot)$ is differentiable at $t$, and is \emph{non-negative} over $[0,T]$. Since $\blw_0(t)=0$, the two properties together imply $\dot{\blw} _0(t)$ must be zero.  We have that $\dot{\bld}_0(t) = -(\dot \blw_0(t) - \dot\bla_0(t)) = n^P$, if $\blw_0(t) =0$.  This proves the case for $k=0$. 

We next consider the case of $k=1, \ldots,  \szK$. For $\dot \bld_k(t)$, the analysis is identical to the case of $k=0$, which we shall omit.\del{We obtain that $\bld_k(t) = r_k m  q^A$, if $\blw_k(t)>0$, and $\bld_k(t) =  \dot{\bla}_k(t)$, otherwise.} We now turn to the analysis for $\dot{\bla}_k(\cdot)$. In particular, we will show that 
\begin{equation}
\lt\{ \dot{\bla}_k(t) \rt\}_{k\in \calK} \leq  (1+\ln^{-1}(1/\delta))  \sum_{h\in \calH} \pi^P_h \calN^*_{h}(\blw(t)),
\label{eq:docak>0}
\end{equation}
where $\calN^*_h(\blw(t))$ was defined in Eq.~\eqref{eq:Nstarw(t)}. The proof of Eq.~\eqref{eq:docak>0} is more involved than for the other coordinates of the fluid solution, because our inspection policy determines the workload vector of a job entering the Adaptive stage using two sources of information:  $(1)$ the coarse estimate of the job's label from the Preparation stage, and $(2)$ the existing workload in the Adaptive stage. Our proof will also proceed in two steps: 
\begin{enumerate}[(a)]
\item We first show that the arrival process of jobs with the same coarse label estimate converges locally to its ``mean value'' under fluid scaling. 
\item We then show that the average workload vector among jobs with the same coarse estimator converges locally to a point that is dominated by the set $(1+\ln^{-1}(1/\delta))  \calN^*_{h}(\blw(t))$.
\end{enumerate}
These two steps together then yield Eq.~\eqref{eq:docak>0}. 

For $t\in (0,T)$, denote by $\calB^s(t)$ the set of all jobs that arrived to the Adaptive stage during the interval $[0, t]$ in the $s$-th system, and by $\calB_h^s(t)$ the subset of jobs in $\calB^s(t)$  whose ML estimators upon exiting the Preparation stage, $\what{H}^P_i$, are equal to $h$. Let $B^s_h(t)$ be the size of the set $\calB^s_h(t)$. The following lemma formalizes step $(a)$ above, whose proof is given in Appendix \ref{app:lem:Bhdiff}. As was mentioned in the beginning of Section \ref{sec:StabilityPrepAdaptive}, the proof of this lemma involves a careful analysis of the potential shuffling in the order in which jobs depart from the Preparation stage. 
\begin{lemma} 
\label{lem:Bhdiff}
Fix $h\in \calH$. For almost all $t\in (0,T)$, we have that
\begin{equation}
\dot{\blb}_h(t) \bydef \lim_{\epsilon\to 0}  \, \lim_{j \to \infty} \,  \frac{B^{\bar{s}_j}_h (\bar{z}_j, t+\epsilon) - B^{\bar{s}_j}_h (\bar{z}_j, t) }{\epsilon}= \left\{ \begin{array}{ll}
          \frac{m q^P}{n^P}  \pi^P_h , & \quad \mbox{if } \blw_0(t)>0,\\
          \\
          \pi^P_h , & \quad  \mbox{if } \blw_0(t)=0,\\
         \end{array} \right.
\label{eq:Bhdrift}
\end{equation}
where $\pi^P_h\bydef \pb\lt( \what{H}_1^P  = h\rt)$, as was defined in Eq.~\eqref{eq:piPdef}. 
\end{lemma}

Step $(b)$ is summarized in the following lemma. It states that, over a small time interval around $t$, the average workload among the jobs in $\calB^{\bar{s}_j}_h(\cdot)$ stays close to a point that is dominated by the set $\calN_h^*(\blw(t))$. The proof, which utilizes  a semi-continuity property of the solution set of the workload-generating optimization problem (Eq.~\eqref{eq:argminLinearWeight}),  is given in Appendix \ref{app:lem:uppersemicont_Nstar}. 
\begin{lemma}
\label{lem:uppersemicont_Nstar}
Fix $t\in (0,T)$, and $h\in \calH$. Let $\overline{\Lambda}_h(t, \epsilon, j)$ be the average workload across all jobs arriving to the Adaptive stage during $[\bar{z}_j t, \bar{z}_j (t+\epsilon)) $ whose ML estimator is $h$, i.e., 
\begin{equation}
\overline{\Lambda}_{h}(t, \epsilon, j) = \frac{1}{B^{\bar{s}_j}_h( \bar{z}_j (t+\epsilon)) -B^{\bar{s}_j}_h( \bar{z}_j  t)}\sum_{i \in 
 \calB^{\bar{s}_j}_h( \bar{z}_j (t+\epsilon))  \backslash  \calB^{\bar{s}_j}_h( \bar{z}_j t)} 
\Lambda_i. 
\end{equation}
We have that
\begin{equation}
\limsup_{\epsilon \downarrow 0}  \limsup_{j \to \infty} \inf_{y \leq \calN_h^*(\blw(t))}\|\overline{\Lambda}_h(t, \epsilon, j)- y\|_2 =0. 
\end{equation}
\end{lemma}

We are now ready to complete the proof of Proposition \ref{prop:limpointarefluid}.  Let $t\in (0,T)$ be a point where all coordinates of $\bla(\cdot)$ are differentiable. We have that
\begin{align}
(\dot{\bla}_k(t))_{k=1, \ldots,  \szK} = & \lim_{\epsilon \downarrow 0}\lim_{j \to \infty}\frac{\bua^{\bar{s}_j}(\bar{z}_j, t+\epsilon)-\bua^{\bar{s}_j}(\bar{z}_j, t)}{\epsilon} \sk{a}{=}  \lim_{\epsilon \downarrow 0}\lim_{j \to \infty} \frac{1}{\epsilon \bar{z}_j}\sum_{h \in \calH} \, \, \, \sum_{i \in
 \calB^{\bar{s}_j}_h( \bar{z}_j (t+\epsilon))  \backslash  \calB^{\bar{s}_j}_h( \bar{z}_j t) } \Lambda_i \nln
= &  \lim_{\epsilon \downarrow 0}\lim_{j \to \infty} \sum_{h \in \calH} \frac{B^{\bar{s}_j}_h(\bar{z}_j(t+\epsilon)) -B^{\bar{s}_j}_h(\bar{z}_j t)}{\epsilon \bar{z}_j} \cdot  \frac{\sum_{i \in 
  \calB^{\bar{s}_j}_h( \bar{z}_j (t+\epsilon))  \backslash  \calB^{\bar{s}_j}_h( \bar{z}_j t)} \Lambda_i}{B^{\bar{s}_j}_h(\bar{z}_j(t+\epsilon)) -B^{\bar{s}_j}_h(\bar{z}_j t)} \nln
\del{= &  \lim_{\epsilon \downarrow 0}\lim_{j \to \infty} \sum_{h \in \calH} \frac{B^{\bar{s}_j}_h(\bar{z}_j(t+\epsilon)) -B^{\bar{s}_j}_h(\bar{z}_j t)}{\epsilon \bar{z}_j} \, \overline{\Lambda}_{h}(t, \epsilon, j)  \nln}
= &  \lim_{\epsilon \downarrow 0}\lim_{j \to \infty} \sum_{h \in \calH} \frac{B^{\bar{s}_j}_h(\bar{z}_j,t+\epsilon) -B^{\bar{s}_j}_h(\bar{z}_j, t)}{\epsilon} \,  \overline{\Lambda}_{h}(t, \epsilon, j) \nln
\sk{b}{=} & \lim_{\epsilon \downarrow 0}\lim_{j \to \infty}   \sum_{h \in \calH}  \dot{\blb}_h(t) \, \overline{\Lambda}_{h}(t, \epsilon, j),
\label{eq:dotak1}
\end{align}
where step $(a)$ follows from the definition of the policy in Eq.~\eqref{eq:Lamb_i_k}: $\bua^s_k(t) = \sum_{h\in \calH} \sum_{i \in \calB^s_h(t)} \Lambda_{i,k}$, and step $(b)$ from Lemma \ref{lem:Bhdiff}. Applying Lemma \ref{lem:uppersemicont_Nstar} for every $h \in \calH$, we have that
\begin{align}
\limsup_{\epsilon \downarrow 0}\limsup_{j \to \infty}  \inf_{y \leq \sum_{h \in \calH}  \dot{\blb}_h(t) \,\calN^*_h (\blw(t))}  \left\| y- \sum_{h \in \calH}  \dot{\blb}_h(t) \, \overline{\Lambda}_{h}(t, \epsilon, j) \right\|_2 = 0. 
\label{eq:dotak2}
\end{align}
Since $\calN_h$, the feasible solutions of the linear program in Eq.~\eqref{eq:argminLinearWeight}, is a compact set, the set of optimal solutions, $\calN_h^*(\blw(t))$, is also compact. The compactness of $\calN_h^*(\blw(t))$ combined with Eq.~\eqref{eq:dotak2} implies 
\begin{equation}
 \lim_{\epsilon \downarrow 0}\lim_{j \to \infty}   \sum_{h \in \calH}  \dot{\blb}_h(t) \, \overline{\Lambda}_{h}(t, \epsilon, j)
\leq  \sum_{h \in \calH}  \dot{\blb}_h(t) \,\calN^*_h (\blw(t)) 
\leq   (1+\ln^{-1}(1/\delta)) \sum_{h \in \calH}  \pi_h^P \calN^*_h (\blw(t)),
\label{eq:dotak3}
\end{equation}
where the last inequality follows from the fact that $ \dot{\blb}_h(t)\leq \pi^P_h\frac{mq^P}{n^P}\leq \pi^P_h(1+\ln^{-1}(1/\delta))$ (Eq.~\eqref{eq:qPdef}). Eqs.~\eqref{eq:dotak1} and \eqref{eq:dotak3} together imply that 
\begin{equation}
(\dot{\bla}_k(t))_{k=1, \ldots,  \szK} \leq (1+\ln^{-1}(1/\delta)) \sum_{h \in \calH}  \pi_h^P \calN^*_h (\blw(t)).
\end{equation}
We have verified that the all conditions in Definition \ref{def:fluidSol} are met, and $\blw(\cdot)$ is a fluid solution. This completes the proof of Proposition \ref{prop:limpointarefluid}. 
\qed



\subsection{Drift Properties of Fluid Solutions}
\label{sec:driftFluid}
We show in this subsection that the fluid solutions exhibit a certain contraction property with respect to the Lyapunov function, $L: \rp^{ \szK+1}\to \rp$, defined as: 
\begin{equation}
L(\blw) = \|\blw \|_2 = \sqrt{\sum_{k = 0}^{\szK} \blw_k^2}, \quad \blw\in \rp^{ \szK+1}. 
\end{equation}
The following proposition is the main result of this subsection, which shows that when $m$ is sufficiently large, the value of the Lyapunov function always decreases by a constant amount starting from any initial condition with unit value. 
\begin{proposition}
\label{prop:drainFluid} Let $m_F^*$ be the optimal value of \flp. Suppose that $m q^P > n^P$, and
\begin{equation}
m q^A > \lt(1+\frac{\ln(2\szH)+ g_\delta}{\ln(1/\delta)}\rt)\lt(1+  \frac{2\szH^2  \old }{\uld\, \ulr \ln(1/\delta) }\rt)(1+\ln^{-1}(1/\delta))\,  m_F^*. 
\end{equation}
Then, there exist $\tau, \epsilon'>0$, such that, for any  $\blw^0\in \rp^{ \szK+1}$ with $L(\blw^0)=1$, 
\begin{equation}
L(\blw(\tau)) \leq 1-\epsilon', \quad \forall \blw \in \calW(\blw^0). 
\label{eq:contraction1}
\end{equation}
\end{proposition}
The proof of the proposition is given in Appendix \ref{app:prop:drainFluid}. A main step of the proof is to couple the drift properties of the fluid solutions, a result of the workload creation in the Adaptive stage, to the structure of the Fundamental Linear Program (\flp) in Definition \ref{def:FLP}, and show that the contraction property holds as long as the system size is approximately greater than the optimal solution of \flp. To this end, we leverage an upper bound on the error of the coarse label estimate produced by the Preparation stage, and show that its accuracy is sufficiently high so that the appropriation of resouces in the Adaptive stage resembles an optimal solution to  \flp.

\subsection{Proof of Stability of Preparation and Adaptive Stages}
\label{sec:completePStabiPrepandAdap}

We now complete the proof of Theorem \ref{thm:prepAndAdapStable} by establishing the joint stability of the Preparation and Adaptive stages. We will use the following version of the Foster-Lyapunov criterion. The main steps in this subsection are similar to those used in \cite{massoulie2007structural}. 
\begin{proposition}[Theorem 8.13, \cite{robert2003}]
\label{prop:roberts}
Let $\{X(t)\}_{t\in \rp}$ be a Markov jump process on a countable state space, $\overline{\calX}$. Suppose that there exists a function $L: \overline{\calX} \to \rp$, constants $C, \epsilon>0$, and an integrable stopping time $\hat \tau>0$, such that for all $x\in \overline{\calX}$ such that $L(x) > C$, we have that 
\begin{equation}
 \E\lt (L(X(\hat \tau))  \bbar X(0)=x\rt ) \leq L(x) - \epsilon \E\lt (\hat \tau \bbar X(0)=x\rt ). 
\label{eq:fosterDrift}
\end{equation}
Suppose, in addition, that the set $\{x: L(x)\leq C\}$ is finite, and that $ \E\lt (L(X(\hat \tau)) \bbar X(0)=x\rt ) < \infty$ for all $x\in \overline{\calX}$. Then, $\{X(t)\}_{t\in \rp}$ is positive recurrent. 
\end{proposition}

Denote by $\{\buu(t)\}_{t\in \rp}$ the Markov process that describes the system dynamics in the Preparation and Adaptive stages, where the states are defined in Section \ref{sec:stateRep}, and by $\calU$ the state space. We will denote by $\pb_\blu(\cdot)$ the probability distribution associated with the process $\buu(\cdot)$ with initial condition $\blu\in \calU$.  Define  the function $\hat{L}: \calU \to \rp$ as the extension of $L(\cdot)$ on $\calU$, i.e., $\hat{L}(\buu(t)) = L(\buw(t)) = \|\buw(t)\|_2$. Let $\tau$ and $\epsilon'$ be defined as in Proposition \ref{prop:drainFluid}. For every $\tilde \blu \in \calU$, let 
\begin{equation}
\hat \tau = \hat L(\tilde{\blu}) \tau = L(\tilde{\blw})\tau. 
\end{equation}
Consider the family of probability measures 
\begin{equation}
\big\{\pb_{\tilde \blu}\big (\buw_{k} (\hat{\tau}) / \hat L(\tilde \blu) \in \cdot \big) \big\}_{k=0, \ldots,  \szK, \tilde \blu\in \calU, \tilde{\blw} \neq 0}. 
\label{eq:familydis1}
\end{equation}
We first show that the above family is uniformly integrable. We have the dominance relation: 
\begin{equation}
\buw_k(t) \leq \buw_k(0)+ \max\{v_\delta, n^P\} \Xi_0(t), \quad t\in \rp, k = 0, \ldots,  \szK. 
\label{eq:wkdomatt}
\end{equation}
To see why this is true, recall that a job creates $n^P$ inspections in the Preparation stage ($k=0$) and at most $v_\delta$ total inspections in the Adaptive stage ($k=1, \ldots,  \szK$). Hence, the second term on the right-hand side of Eq.~\eqref{eq:wkdomatt} dominates the total number of inspections that could have been added to $\buw_k(\cdot)$ by time $t$. Scaling both sides of Eq.~\eqref{eq:wkdomatt} by $\hat L (\tilde \blu) $ and setting $t=\hat{\tau}$, we have that, when $\buu(0) = \tilde \blu$, 
\begin{align}
\buw_k(\hat \tau) / \hat L (\tilde \blu) \leq & \tilde \blw_k/ \hat L( \tilde \blu) + \max\{v_\delta, n^P\} \Xi_0(\hat \tau)/ \hat L( \tilde \blu) \nln
\del{=& \tilde \blw_k/ \hat L( \tilde \blu) + \max\{v_\delta, n^P\} \Xi_0(\hat L( \tilde \blu) \tau)/ \hat L( \tilde \blu) \nln}
{=} & \tilde \blw_k/  L( \tilde \blw)  +  \max\{v_\delta, n^P\} \Xi_0(\hat L( \tilde \blu) \tau)/ \hat L( \tilde \blu) \nln
\sk{a}{\leq} & {1}/{{\alpha_1}}  +  \max\{v_\delta, n^P\} \Xi_0(\hat L( \tilde \blu) \tau)/ \hat L( \tilde \blu),
\label{eq:wkboundafterscaling}
\end{align}
for some constant $\alpha_1>0$, where step $(a)$ follows from the first inequality in Eq.~\eqref{eq:L1homo} of Lemma \ref{lem:L1homo} in Appendix \ref{app:lem:L1homo}. Since $\Xi_0(\cdot)$ is a unit-rate Poisson process, we have $\E \lt( \Xi_0(\hat L( \tilde \blu) \tau)/ \hat L( \tilde \blu) \rt) = \var{\Xi_0(\hat L( \tilde \blu) \tau)/ \hat L( \tilde \blu)} = \tau.$ The above arguments show that the first and second moments of $\buw_k(\hat \tau)/ \hat L (\tilde \blu) $ are both bounded from above uniformly over all $\tilde \blu \in \calU$, so long as $\tilde{\blw} \neq 0$. This implies the uniform integrability of the family of distributions in Eq.~\eqref{eq:familydis1}. Using again Eq.~\eqref{eq:L1homo} of Lemma \ref{lem:L1homo}, we have that 
\begin{align}
\hat L (\buu (\hat \tau )) / \hat L (\tilde \blu ) = &  L (\buw (\hat \tau )) / \hat L (\tilde \blu ) \leq \alpha_2\|\buw(\hat \tau) \|_\infty / \hat L (\tilde \blu )  \leq \frac{\alpha_2}{\alpha_1}  +  \max\{v_\delta, n^P\} \Xi_0(\hat L( \tilde \blu) \tau)/ \hat L( \tilde \blu).
\label{eq:ukboundafterscaling}
\end{align}
Using the same arguments as those following Eq.~\eqref{eq:wkboundafterscaling}, we have that the set of probability measures
\begin{equation}
\big\{\pb_{\tilde \blu} \big ( \hat L(\buu(\hat{\tau}) )  / \hat L(\tilde \blu) \in \cdot \big) \big\}_{\tilde \blu\in \calU, \tilde{\blw} \neq 0}
\label{eq:familydis2}
\end{equation}
is also uniformly integrable. 

We now prove the validity of Eq.~\eqref{eq:fosterDrift} in our context. In particular, we will show that there exist $A, \epsilon>0$, 
such that for all  $\tilde{\blu} \in \calU$, if $\hat L(\tilde \blu) \geq A$, then we have that
\begin{equation}
\E_{\tilde \blu }\lt (\hat L (\buu (\hat L (\tilde \blu) \tau ) ) \rt ) \leq \hat L(\tilde \blu) (1-\epsilon \tau), 
\label{eq:driftcond1}
\end{equation}
where we have used the substitution $\hat \tau = \hat L (\tilde \blu) \tau$. By the definitions of $\hat L$ and $\hat \tau$, we have that
\begin{align}
\E_{\tilde \blu }\lt ({\hat L (\buu (\hat L(\tilde \blu) \tau) ) }/{\hat L (\blu) } \rt ) = & \E_{\tilde \blu }\big ( {L (\buw (L(\tilde \blw) \tau) )}/{L (\tilde \blw)} \big ) \sk{a}{=}  \E_{\tilde u} \big( L(L(\tilde \blw)^{-1} \buw (L(\tilde \blw) \tau))  \big)  \nln
= & \E_{\tilde u} \big( L( \buw (L(\tilde \blw), \tau))  \big),
\label{eq:Lusimp1}
\end{align}
where step $(a)$ is based on the second property in Lemma \ref{lem:L1homo} (one-homogeneity of $L(\cdot)$). In light of Eq.~\eqref{eq:Lusimp1}, the condition in Eq.~\eqref{eq:driftcond1} is equivalent to 
\begin{equation}
\E_{\tilde \blu} \lt( L( \buw (L(\tilde \blw), \tau))  \rt)\leq 1-\epsilon \tau. 
\label{eq:driftcond2}
\end{equation}
We will now prove Eq.~\eqref{eq:driftcond2} by contradiction. For the sake of contradiction, suppose that there exists a sequence of initial conditions, $\{\tilde \blu^s \}_{s \in \N}$, where $\hat L (\tilde \blu ^s) \to \infty $ as $s \to \infty$, such that 
\begin{equation}
\limsup_{s \to \infty} \E_{\tilde \blu ^s } \lt( L( \buw ^s (L(\tilde \blw^s  ), \tau))  \rt)\geq 1, 
\label{eq:contraassmEL}
\end{equation}
where we have re-introduced the superscript $s$ in $\buw$ to signify the initial condition $\buw^s(0) = \tilde \blw ^s$. In particular, after scaling, we can define the quantity $\blw^s$: 
\begin{equation}
\blw^s \bydef \buw ^s (L(\tilde \blw^s  ), 0) = \buw^s(0)/L(\tilde \blw^s  ) = \tilde \blw ^s / L(\tilde \blw^s  ). 
\end{equation}
Since $L(\blw) = \|\blw\|_2$, we have that $\|\blw^s\|_2 = L(\blw ^s / L(\tilde \blw^s  )) = 1$, i.e., $\blw^s$ belongs to the compact set $\{\blw\in \rp^{ \szK+1}: \|\blw\|_2=1\}$. The compactness implies that there exists $\blw^0\in \rp^{ \szK+1}$, with $L(\blw^0) = \|\blw^0\|_2=1$, and a sub-sequence of $\{\blw^s\}$, $\{\blw^{s_i}\}_{i\in \N}$, such that
\begin{equation}
\lim_{i \to \infty} \blw^{s_i} = \blw^0. 
\label{eq:wsitow0}
\end{equation}
Since $L(\tilde \blw ^s) = \hat L(\tilde \blw ^s) \to \infty$ as $s\to \infty$ by our assumption, by Proposition \ref{prop:convrgToFluid} and Eq.~\eqref{eq:wsitow0}, we have 
\begin{equation}
\lim_{i \to \infty} \inf_{\bly \in \calV(\blw^0, \tau)} \| \bly - \buw ^{s_i} (L(\tilde \blw^{s_i}  ), \tau) \|_2 =0, \quad \mbox{a.s.}, 
\label{eq:yclosetoWlws}
\end{equation}
where $\calV(\blw^0,t)$ is defined to be the set of all states of the fluid solution at time $t \in \rp$, starting with an initial condition $\blw^0$: $\calV(\blw^0, t) = \bigcup_{\blw(\cdot) \in \calW(\blw^0)} \blw(t)$. It can be verified from the definition of the fluid solutions that the set $\calV(\blw^0,t)$ is compact for all $t$ and $\blw^0$. 
The compactness of $\calV(\blw^0, \tau)$ along with Eq.~\eqref{eq:yclosetoWlws} implies that there exist $ \bly \in \calV(\blw^0, \tau)$ and  a sub-sequence of $\{\tilde \blw^{s_i}\}$, $\{\tilde \blw^{s_{i_j}}\}_{j \in \N}$, such that 
\begin{equation}
\lim_{j\to \infty} \buw ^{s_{i_j}} (L(\tilde \blw^{s_{i_j}}) , \tau)  = \bly, \quad \mbox{a.s.}
\end{equation}
Since $L(\cdot)$ is a continuous function, the above equation further implies that 
\begin{equation}
\lim_{j\to \infty} L( \buw ^{s_{i_j}} (L(\tilde \blw^{s_{i_j}}) , \tau) ) = L(\bly) \leq 1-\epsilon',  \quad \mbox{a.s.}, 
\label{eq:wsijA.S.converg}
\end{equation}
where the last inequality follows from the fact that $L(\blw^0)=1$, $\bly\in \calV(\blw^0, \tau)$, and Proposition \ref{prop:drainFluid}. As was shown earlier, the family of distributions, $\lt\{\pb_{\tilde \blu}\lt ( \hat L(\buu(\hat{\tau}) )  / \hat L(\tilde \blu) \in \cdot \rt) \rt\}_{\tilde \blu\in \calU, \tilde{\blw} \neq 0}$,  is uniformly integrable (Eq.~\eqref{eq:familydis2}), and hence one of its subsets
\begin{equation}
\lt\{\pb_{\tilde{\blu} ^{s_{i_j}}}(L( \buw ^{s_{i_j}} (L(\tilde \blw^{s_{i_j}}) , \tau) \in \cdot ) \rt\}_{i \in \N}
\end{equation}
is also uniformly integrable. Combining Eq.~\eqref{eq:wsijA.S.converg} with the above uniform integrability implies that 
\begin{equation}
\lim_{j\to \infty} \E_{\tilde{\blu}^{s_{i_j}}}(L( \buw ^{s_{i_j}} (L(\tilde \blw^{s_{i_j}}) , \tau) ))  = L(\bly) \leq 1-\epsilon'<1, 
\end{equation}
which contradicts with Eq.~\eqref{eq:contraassmEL}. This completes the proof of Eq.~\eqref{eq:driftcond1}. 

The fact that $\E_{\tilde \blu}(\tilde{L} (\buu (\hat\tau) ) < \infty$ for all $\tilde \blu$ is readily verifiable from Eq.~\eqref{eq:ukboundafterscaling}. Finally, for any finite $C$, the set $\{\blu \in \calU: L(\blu) < C\}$ contains states where the total number of jobs currently in system is bounded. Since each job can receive at most $n^P+v_\delta$ inspections, and the result of each inspection belongs to a finite set, it follows (Section \ref{sec:stateRep}) that the set $\{\blu \in \calU: L(\blu) < C\}$ is finite. We have thus verified all the conditions in Proposition \ref{prop:roberts}, and established the positive recurrence of $\buu(\cdot)$. This completes the proof of Theorem \ref{thm:prepAndAdapStable}.   \qed

\section{Stability of Residual Stage}
\label{sec:StabilityResidual}

The following result gives a sufficient condition on the stability of the Residual stage; the proof is given in Appendix \ref{app:prop:residualStable}. 
\vspace{-10pt}
\begin{proposition}  
\label{prop:residualStable}
Suppose that the Preparation and Adaptive stages are stable. Then, the Residual stage is stable whenever
\begin{equation}
m q^R >3\szH\zeta_0(1+\ln(4\szH) \ln^{-1}(1/\delta)).
\end{equation}
\end{proposition}

\section{Proof of Theorem \ref{thm:main}}
\label{sec:CompleteProof}

We now complete the proof of  our main result, Theorem \ref{thm:main}. We begin by showing that the policy proposed in Section \ref{sec:DefOfPolicy} is $\delta$-accurate: for \emph{every} job $i$ that departs from the system, the probability of it being misclassified is at most $\delta$. The proof is given in Appendix \ref{app:prop:classErr}. 

\begin{proposition}\label{prop:classErr} 
Fix $\delta \in (0,1)$. Under the inspection policy described in Section \ref{sec:DefOfPolicy}, 
\begin{equation}
\pb\lt(\what{H}_i \neq H_i \bbar H_i = h \rt) \leq \delta, \quad \forall h \in \calH, i \in \N. 
\label{eq:exitAccu}
\end{equation}
\end{proposition}

\paragraph{Completing the Proof of Theorem \ref{thm:main}.} We now show that the minimum value of $m$ required to stabilize the system under the proposed inspection policy satisfies the inequality in Eq.~\eqref{eq:mainIneq}. Combining Theorem \ref{thm:prepAndAdapStable} and Proposition \ref{prop:residualStable}, we have that the system is stable as long as all of the following hold: 
\begin{align}
\label{eq:mrequirement1}
m q^A > &  \lt(1+\frac{\ln(2\szH)+ g_\delta}{\ln(1/\delta)}\rt)\lt(1+  \frac{2\szH^2  \old }{\uld\, \ulr \ln(1/\delta) }\rt) m_F^*, \\
\label{eq:mrequirement2}
m q^P > & \zeta_0 \ln\ln(1/\delta),  \quad m q^R >3\szH\zeta_0(1+\ln(4\szH) \ln^{-1}(1/\delta)),
\end{align}
where $m^*_F$ is the optimal value of \flp (Definition \ref{def:FLP}).  Recall from the definition of our policy in Eqs.~\eqref{eq:qPdef} and \eqref{eq:qAdef} that 
$ q^P = \lt(\zeta_0 \ln\ln(1/\delta)+\ln^{-1}(1/\delta))\rt) /{m}$, $ q^R = \lt[3\szH\zeta_0(1+\ln(4\szH) \ln^{-1}(1/\delta)) +1 \rt]/{m}$, and $q^A =  1-q^P-q^R$. We have that Eq.~\eqref{eq:mrequirement2} is automatically satisfied  by construction, whenever $q^P, q^R \in (0,1)$. Therefore, the system is stable if Eq.~\eqref{eq:mrequirement1} is true, which, by inserting the expressions for $q^P$ and $q^R$, requires that
\begin{equation}
m >   \lt(1+\frac{\ln(2\szH)+ g_\delta}{\ln(1/\delta)}\rt)\lt(1+  \frac{2\szH^2  \old }{\uld\, \ulr \ln(1/\delta) }\rt) (1+\ln^{-1}(1/\delta)) m_F^* + \iota_\delta, 
\label{eq:msufficient1}
\end{equation}
where $\iota_\delta \bydef \zeta_0 \ln\ln(1/\delta)+3\szH\zeta_0(1+\ln(4\szH) \ln^{-1}(1/\delta))+ 1+ \ln^{-1}(1/\delta)$. By dividing both sides of Eq.~\eqref{eq:msufficient1} by $m_F^*$, the condition becomes
\begin{equation}
\frac{m}{m_F^*} >   \lt(1+\frac{\ln(2\szH)+ g_\delta}{\ln(1/\delta)}\rt)\lt(1+  \frac{2\szH^2  \old }{\uld\, \ulr \ln(1/\delta) }\rt)(1+\ln^{-1}(1/\delta)) + \frac{\iota_\delta}{ m_F^*}. 
\label{eq:mtom1}
\end{equation}
The two terms on the right-hand side of the above equation can be bounded as follows. For the first term, recall from the definition in Eq.~\eqref{eq:gdeltadef} that $
g_\delta  = 3\olz\uld^{-1/2} \sqrt{ \ln(1/\delta)\ln\ln(1/\delta)}$. It is not difficult to see that there exist  constants $\bar{c}_1$ and $\bar{\delta}_1>0$, independent of $\pi$, such that for all $\delta < \bar{\delta}_1$
\begin{align}
\lt(1+\frac{\ln(2\szH)+ g_\delta}{\ln(1/\delta)}\rt)\lt(1+  \frac{2\szH^2  \old }{\uld\, \ulr \ln(1/\delta) }\rt)(1+\ln^{-1}(1/\delta))  < & 1+ \frac{\bar{c}_1\sqrt{\ln(1/\delta)\ln\ln(1/\delta)}}{\ln(1/\delta)} =  1+ \bar{c}_1 \sqrt{ \frac{\ln\ln(1/\delta)}{\ln(1/\delta)} }. 
\label{eq:mtomm1-t1}
\end{align}
For the second term in Eq.~\eqref{eq:mtom1},  $\iota_\delta/m^*_F$, note that the dominating term in $\iota_\delta$ is of order $\ln\ln(1/\delta)$. Recall from Lemma \ref{lem:lp1structure} that $m_F^* \geq \old^{-1}\ln(1/\delta)$, and hence there exist constants $\bar{c}_2$ and $\bar{\delta}_2>0$, independent of $\pi$, such that for all $\delta < \bar{\delta}_2$, we have that $\frac{\iota_\delta}{ m_F^*} \leq \frac{\iota_\delta}{\old^{-1}\ln(1/\delta)} \leq \bar{c}_2 \frac{\ln\ln(1/\delta)}{\ln(1/\delta)}$. This, along with Eqs.~\eqref{eq:mtom1} and \eqref{eq:mtomm1-t1},  implies that, for all $\delta<\min\{\bar{\delta}_1, \bar{\delta}_2\}$, it suffices to have 
\begin{equation}
\frac{m}{m_F^*} > 1+ \bar{c}_1 \sqrt{ \frac{\ln\ln(1/\delta)}{\ln(1/\delta)} } + \bar{c}_2 \frac{\ln\ln(1/\delta)}{\ln(1/\delta)}. 
\label{eq:mtomm1-clean}
\end{equation}
Finally, the lower bound in Proposition \ref{prop:mstarLWB} shows that the optimal system size, $m^*(\delta,\pi)$, satisfies
\begin{equation}
m^*(\delta,\pi)  \geq (1 - \delta)\lt[ 1 - \lt(\ln\frac{1}{1-\delta}+e^{-1} \rt) \ln^{-1}(1/\delta) \rt] m_F^*. 
\end{equation}
In particular, there exist constants $\bar{c}_3$ and $\bar{\delta}_3>0$, independent of $\pi$, such that for all $\delta<\bar{\delta}_3$
\begin{equation}
 m_F^* \leq \lt(1+ \frac{\bar{c}_3}{\ln(1/\delta)}\rt) m^*(\delta,\pi) . 
\label{eq:m1tombar}
\end{equation}
Substituting Eq.~\eqref{eq:m1tombar} into Eq.~\eqref{eq:mtomm1-clean}, and noticing that the leading  term in $\lt( \frac{m}{m^*(\delta,\pi)}-1 \rt) $ is of order $\sqrt{ \frac{\ln\ln(1/\delta)}{\ln(1/\delta)} }$, we conclude that there exist $c_0, \delta_0>0$, independent of $\pi$, such that, for all $\delta \in (0,\delta_0)$, the system is stable under the proposed policy whenever 
\begin{equation}
\frac{m}{m^*(\delta,\pi)} > 1+ {c_0} \sqrt{ \frac{\ln\ln(1/\delta)}{\ln(1/\delta)} }. 
\end{equation}
Together with Proposition \ref{prop:classErr}, this completes the proof of Theorem \ref{thm:main}. \qed

\section{Concluding Remarks}

The main objective of this paper is to understand the design principles and fundamental limitations involved as one aims to efficiently classify a large set of items using a finite amount of processing resources. The main result demonstrates a prior-oblivious inspection architecture that asymptotically  uses the minimum number of experts, in the regime where the required classification error tends to zero. More broadly speaking, our result could be viewed as an attempt towards understanding how to build effective processing architectures and algorithms for large-scale statistical learning or information extraction tasks, given limited resources or processing power. We believe that there are many other problems in this domain, situated at the intersection between stochastic modeling and statistics, that may be of interest for future research. 

The present paper leaves open three main questions. First, we have thus far required the inspection policies to be stable, while the more refined metric of delays experienced by the jobs has not been investigated. Second, the resource efficiency property of the proposed policy applies only in the regime where the allowable error is close to zero. If significantly greater errors can be tolerated, however, it becomes less clear what inspection policy one should choose and whether a different criterion of resource efficiency should be adopted. Finally, it would be interesting to identify, and provide rigorous guarantees for, simpler heuristic policies, such as the one we propose in Appendix \ref{sec:heuristicPolicy}, which are more easily implementable in practical applications. 


\bibliographystyle{ormsv80}
\bibliography{EXPERT.bib} 


\begin{APPENDICES}

\normalsize

\section{Proofs}

\subsection{Proof of Lemma \ref{lem:EofShl}}
\label{app:lem:EofShl}

\bpf  Fix $i\in \N$. Denote by $N_i$ the total number of inspections job $i$ receives before departing from the system. Let $Y_{i,n}$ be the {inspection history} of job $i$ up till the $n$th inspection: 
\begin{equation}
Y_{i,n} = \{(K_{i,j}, X_{i,j})\}_{j = 1, \ldots, n}.
\label{eq:historyDef} 
\end{equation}
We will use $y_n = \{(k_{j},x_{j})\}_{j = 1, \ldots, n} $ to denote a specific realization of $Y_{i,n}$. Define $p_{h}(y_n)$ as the likelihood associated with an inspection history, $y_n=\{(k_{j},x_{j})\}_{j = 1, \ldots, n}$ under label $h$: 
\begin{equation}
p_{h}(y_n) = \prod_{j=1}^n p(h,k_j,x_j). 
\label{eq:historyLikelihoodDef0}
\end{equation}
Define $q_h(y_n)$ as the probability
\begin{align}
q_h(y_n) = & \pb( N_i= n, Y_{i,n}=y_n \bbar H_i = h ).
\label{eq:qhyndef}
\end{align}
We have that
\begin{align}
q_h(y_n) = & \pb( N_i = n, Y_{i,n}=y_n \bbar H_i = h ) \nln
= & \pb(N_i = n \bbar Y_{i,n}=y_n)\prod_{j=1}^n p(h,k_{j}, x_{j})\pb(K_{i,j}= k_{j}, N_i\geq j\bbar Y_{i,j-1}=y_{j-1}) \nln
=&  \Big(\prod_{j=1}^n p(h,k_{j}, x_{j})\Big)\Big(\pb(N_i = n \bbar Y_{i,n}= y_n)\prod_{j=1}^n \pb(K_{i,j}= k_{j},  N_i\geq j \bbar Y_{i,j-1}= y_{j-1})\Big)  \nln
= & p_h(y_n)v(y_n)
\label{eq:fhy0}
\end{align}
where $y_j$, $1\leq j\leq n$, denotes the first $j$ inspections in $y_n$, with $y_0 \bydef \emptyset$, and $v(y_n) \bydef \pb(N_i = n \bbar Y_{i,n}= y_n)\prod_{j=1}^n \pb(K_{i,j}= k_j,  N_i\geq j \bbar Y_{i,j-1}= y_{j-1})$.  Importantly, note that for a given policy, $v(y_n)$ only depends on the realization of the history, and not on the choice of $h$.   

We have that
\begin{align}
\E\lt(\frac{p_{l}(Y_{i, N_i})}{p_{h}(Y_{i, N_i})} \hbar \what{H}_i = h,  H_i = h\rt) = \E\lt( \frac{p_l(Y_{i, N_i})}{p_h(Y_{i, N_i})} \mathbb{I}(\what{H}_i = h)  \hbar H_i = h\rt)/ \pb(\what{H}_i = h \bbar H_i=h),
\label{eq:Eplph1}
\end{align}
where
\begin{align}
\E\lt( \frac{p_l(Y_{i, N_i})}{p_h(Y_{i, N_i})} \mathbb{I}(\what{H}_i = h)  \bbar H_i = h\rt) \sk{a}{=} & \sum_{n=1}^\infty \sum_{y_n} \frac{p_l(y_n)}{p_h(y_n) } \pb(N_i = n, Y_{i,n}=y_n, \what{H}_i = h \bbar H_i = h)  \nln
= & \sum_{n=1}^\infty \sum_{y_n} \frac{p_l(y_n)}{p_h(y_n) } q_h(y_n)\pb(\what{H}_i = h \bbar N_i = n, Y_{i,n}=y_n,  H_i = h)  \nln
\sk{b}{=} & \sum_{n=1}^\infty \sum_{y_n} p_l(y_n) v(y_n)\pb(\what{H}_i = h\bbar N_i = n, Y_{i,n}=y_n,  H_i = h)  \nln
\sk{c}{=} & \sum_{n=1}^\infty \sum_{y_n} p_l(y_n)v(y_n) \pb(\what{H}_i = h\bbar N_i = n, Y_{i,n}=y_n,  H_i = l)  \nln
\sk{d}{=} & \sum_{n=1}^\infty \sum_{y_n} q_l(y_n) \pb(\what{H}_i = h\bbar N_i = n, Y_{i,n}=y_n,  H_i = l)  \nln
=  & \sum_{n=1}^\infty \sum_{y_n} \pb(N_i = n, Y_{i,n}=y_n, \what{H}_i = h\bbar  H_i = l)  \nln
= & \pb(\what{H}_i = h \bbar H_i = l). 
\label{eq:Eplph2}
\end{align}
The second summation in step $(a)$ is over all possible realizations of $Y_{i,n}$ for which $p_h(y_n)$ is non-zero. Steps $(b)$ and $(d)$ follow from Eq.~\eqref{eq:fhy0}. For step $(c)$, we used the fact that  conditioning on the event $\{N_i = n, Y_{i,n}=y_n\}$, the event $\{\what{H}_i=h\}$ is independent from the true label of job $i$. Eqs.~\eqref{eq:Eplph1} and \eqref{eq:Eplph2} combined yield that
\begin{equation}
 \E\lt(\frac{p_l(Y_{i, N_i})}{p_h(Y_{i, N_i})} \hbar \what{H}_i = h,  H_i = h\rt) = \frac{\pb(\what{H}_i = h\bbar H_i = l)}{\pb(\what{H}_i=h \bbar H_i = h)} \leq \frac{\delta}{1-\delta}, 
 \label{eq:Eplphall}
 \end{equation}
 where the last equality follows from the assumption of the policy being $\delta$-accurate.

Define
\begin{equation}
\delta_h = \pb(\what{H}_i \neq h \bbar H_i = h).
\end{equation}
Using essentially the same steps of deduction as those leading to Eq.~\eqref{eq:Eplphall}, we have that
\begin{equation}
 \E\lt(\frac{p_l(Y_{i, N_i})}{p_h(Y_{i, N_i})} \hbar \what{H}_i \neq h,  H_i = h\rt) = \frac{\pb(\what{H}_i \neq h \bbar H_i = l)}{\pb(\what{H}_i \neq h \bbar H_i = h)} \leq \frac{1}{\delta_h}
\label{eq:Eplphall2} 
\end{equation}

We are now ready to prove the main claim of the lemma. We have that
\begin{align}
\E\lt(S_{i}(h,l) \bbar H_i = h\right) = &  \E\lt(\sum_{j=1}^{N_i}\ln\frac{p(h,K_{i,j}, X_{i,j})}{p(l,K_{i,j}, X_{i,j})}\hbar  H_i=h\right) \nln
= &  \E\lt(- \ln\frac{p_l(Y_{i, N_i})}{p_h(Y_{i, N_i})} \hbar  H_i=h \right) \nln
= & - \E\lt(\ln\frac{p_l(Y_{i, N_i})}{p_h(Y_{i, N_i})} \hbar \what{H}_i = h,  H_i=h \right) \pb(\what{H}_i =h \bbar H_i = h) \nln
& -\E\lt(\ln\frac{p_l(Y_{i, N_i})}{p_h(Y_{i, N_i})} \hbar \what{H}_i \neq h,  H_i=h\right) \pb(\what{H}_i \neq h \bbar H_i = h)  \nln
\sk{a}{\geq} & - \ln \lt(\E\lt(\frac{p_l(Y_{i, N_i})}{p_h(Y_{i, N_i})} \hbar \what{H}_i = h, H_i=h\right) \rt) \pb(\what{H}_i =h \bbar H_i = h) \nln
& -\ln \lt( \E\lt(\frac{p_l(Y_{i, N_i})}{p_h(Y_{i, N_i})} \hbar \what{H}_i \neq h,  H_i=h\right) \rt) \pb(\what{H}_i \neq h\bbar H_i = h ) \nln
\sk{b}{\geq} & (1-\delta) \ln\frac{1-\delta}{\delta} +  \delta_h \ln \delta_h \nln
\sk{c}{\geq} & (1-\delta) \ln\frac{1-\delta}{\delta} - e^{-1}. 
\end{align}
Step $(a)$ is based on the Jensen's inequality, by noting that that $-\ln(\cdot)$ is a convex function. Step $(b)$ follows from Eqs.~\eqref{eq:Eplphall} and \eqref{eq:Eplphall2}, and the policy being $\delta$-accurate. Step $(c)$ is based on the observation that $(x\ln x)$ is monotonically decreasing over $x\in (0,e^{-1})$ and monotonically increasing over $(e^{-1},1)$. This proves Lemma \ref{lem:EofShl}.  \qed

\subsection{Proof of Proposition \ref{prop:mstarLWB}}
\label{app:prop:mstarLWB}

\bpf Fix a stable inspection policy, $\psi$. Denote by $N_i$ the total number of inspections received by job $i$ before a classification is produced, and by $S_i(h,l)$ the sum of the log-likelihood ratios by the time when job $i$ departs from the system: 
\begin{equation}
S_{i}(h,l) =  \sum_{j=1}^{N_i}Z_{i,j}(h,l, K_{i,j}). 
\label{eq:Sifinaldef}
\end{equation}
Define: 
\begin{enumerate}
\item $\calI_h^n$, the set of the first $n \pi_h $ jobs of true label $h$ to depart from the system;
\item $i^j_h$, the index of the $j$th job of true  label $h$ to depart from the system; 
\item $T^n_h$, the first time when $n \pi_h $ jobs of true  label $h$ have departed form the system; 
\item $N_{i,k}$,  the number of times that job $i$ has been inspected by an expert of type $k$ by the time of its departure;  
\item $M^t_k$, the total number of inspections completed by  the type-$k$ experts by time $t$. 
\end{enumerate}

Fix $h\in \calH$. Let $T^n$ be the time by which there has been at least $\pi_h n$ departures from every job label: 
\begin{equation}
T^n = \max_{h\in \calH} T_h^n. 
\end{equation}
We have that, for all $k\in \calK$, 
\begin{align}
M_k^{T^n} \geq  \sum_{h\in \calH} \sum_{i \in \calI^n_h} N_{i,k}  =  \sum_{h\in \calH} \sum_{j=1}^{n\pi_h } N_{i^j_h, k},	
\label{eq:Mk1}
\end{align}
where the first inequality follows from the definition of $M^t_k$ and $T^n$. Take expectations of both sides of Eq.~\eqref{eq:Mk1}, and define the variables
\begin{equation}
{\upsilon}^n_{h,k} = \frac{1}{n\pi_h }\sum_{j=1}^{n\pi_h} \E\lt( N_{i^j_h, k} \rt). 
\label{eq:upsilonnhknonNeg}
\end{equation}
We have that
\begin{align}
\E(M_k^{T^n} ) \geq   \sum_{h\in \calH} \sum_{j=1}^{n\pi_h } \E\lt( N_{i^j_h, k} \rt)= n \sum_{h\in \calH}\pi_h {\upsilon}^n_{h,k}. 
\label{eq:Mk2}
\end{align}

Recall that jobs arrive to the system according to a unit rate Poisson process, and each job has label $h$ with probability $\pi_h$. Since we have assumed that the inspection policy is stable, the time when there are $n\pi_h$ jobs with label $h$  that have departed from the system, $T_h^n$,  cannot be much later than when the $(n\pi_h)$th job with label $h$ arrived to the system, which, by the strong law of large numbers applied to a Poisson process, is at most $n+o(n)$ almost surely, as $n\to \infty$. Formally, there exists a function $f: \rp \to \rp$, with $\lim_{x \to \infty}x^{-1}f(x) =0$,  such that, almost surely, 
\begin{equation}
\limsup_{n \to \infty} \lt[T^n_h - (n+ f(n)) \rt] \leq 0, \quad \forall h \in \calH. 
\label{eq:Tnalmostn}
\end{equation}
We thus have that
\begin{equation}
\limsup_{n \to \infty}\frac{1}{n}\E(M_k^{T^n} ) \sk{a}{\leq} \limsup_{n \to \infty} \frac{1}{n}\E(M_k^{n+f(n)} ) \sk{b}{\leq}  \limsup_{n \to \infty} m r_k(1+f(n)/n) = m r_k, 
\label{eq:EMTUPB}
\end{equation}
Step $(a)$ follows from Eq.~\eqref{eq:Tnalmostn} and the definition of $T^n$. Step $(b)$ is based on the fact that the total number of completed inspections by experts of type $k$ during an interval $[0,x]$ is stochastically dominated by a Poisson random variable with mean $mr_k x$, which corresponds to the case where all experts of type $k$ work without any idling during that interval. Define 
\begin{equation}
n_{h,k} = \limsup_{n \to \infty}\upsilon_{h,k}^n. 
\end{equation}
Combining Eqs.~ \eqref{eq:Mk2} and \eqref{eq:EMTUPB}, we have that
\begin{equation}
 \sum_{h \in \calH} \pi_h {n}_{h,k}=  \sum_{h \in \calH} \pi_h \limsup_{n \to \infty}\upsilon_{h,k}^n \leq \limsup_{n \to \infty}\frac{1}{n}\E(M_k^{T^n} ) \leq mr_k, \quad \forall k \in \calK. 
\label{eq:avgOverH}
\end{equation}

By Lemma \ref{lem:EofShl}, we have that, for any job $i$ whose true label is $h$, 
\begin{equation}
 \sum_{k \in \calK} \E(N_{i,k})D(h,l,k) \sk{a}{=} \E\Big (\sum_{ j=1}^{ N_{i}}Z_{i,j}(h,l, K_{i,j}) \Big) = \E(S_{i}(h,l)) \geq (1-\delta)\ln\frac{1-\delta}{\delta} - e^{-1},
 \label{eq:NikD}
\end{equation}
where step $(a)$ follows from Wald's identity by noting that $N_{i,k}$ is a stopping time with respect to the sequence $Z_{i,1}(h,l,K_{i,1}), Z_{i,2}(h,l,K_{i,2}), \ldots$. Summing both sides of Eq.~\eqref{eq:NikD} over the set $\calI_h^n$, we have that, for all $l \neq h$, 
\begin{align}
\sum_{k\in \calK} \upsilon_{h,k}^nD(h,l,k) = &  \sum_{k\in \calK} \Big( \frac{1}{n \pi_h}   \sum_{j=1}^{n\pi_h } \E \big( N_{i^j_h, k} \big) \Big)  D(h,l,k)  \nln
 =& \frac{1}{n \pi_h } \sum_{j=1}^{n \pi_h } \Big( \sum_{k \in \calK}\E\big( N_{i^j_h, k} \big)   D(h,l,k) \Big) \nln
\geq & (1-\delta)\ln\frac{1-\delta}{\delta} - e^{-1}.
\end{align}
Taking the limits as $n\to \infty$ on both sides of the above equation, we obtain that
\begin{equation}
\sum_{k\in \calK} n_{h,k} D(h,l,k) \geq (1-\delta)\ln\frac{1-\delta}{\delta} - e^{-1}. 
\label{eq:avgOverK}
\end{equation}

Note that the $\upsilon^n_{h,k}$'s are non-negative by definition (Eq.~\ref{eq:upsilonnhknonNeg}) and hence $n_{h,k}\geq 0$ for all $h$ and $k$.  From Eqs.~\eqref{eq:avgOverH} and \eqref{eq:avgOverK}, we conclude that the number of experts, $m$, under any $\delta$-accurate inspection policy must be no more than the optimal value of the following linear program: 
\begin{align}
\mbox{minimize}  &   \quad\quad\quad m \\
\mbox{s.t.} & \quad   \sum_{k\in \calK} {n}_{h,k}D(h,l,k)\geq  (1-\delta)\ln\frac{1-\delta}{\delta} - e^{-1} , \quad \forall h, l \in\calH, h\neq l, \label{eq:constr1}\\
&\quad \sum_{h \in \calH} {n}_{h,k} \pi_h  \leq r_k m, \quad \forall k\in \calK, \label{eq:constr2} \\
& \quad n_{h,k} \geq 0, \quad \forall h\in \calH, k\in \calK
\label{eq:LP2}
\end{align}
where Eqs.~\eqref{eq:avgOverH} and \eqref{eq:avgOverK} correspond to the constraints in Eqs.~\eqref{eq:constr2} and \eqref{eq:constr1}, respectively. Note that the above linear program differs from \flp (Definition \ref{def:FLP}) only through changing $\ln(1/\delta)$ to $(1-\delta)\ln\frac{1-\delta}{\delta} - e^{-1}$ in Eq.~\eqref{eq:constr2}. Hence, it is not difficult to verify that any value, $m$, associated with a feasible solution of \flp must satisfy
\begin{equation}
\frac{m}{\overline{m}(\delta, \pi)} \geq \frac{(1-\delta)\ln\frac{1-\delta}{\delta} - e^{-1} }{\ln(1/\delta)} = 1 - \delta-\frac{(1-\delta)\ln\frac{1}{1-\delta}+e^{-1}}{\ln(1/\delta)}.  
\end{equation}
The completes the proof of Proposition \ref{prop:mstarLWB}. \qed

\subsection{Proof of Lemma \ref{lem:lp1structure}} 
\label{app:lem:lp1structure}

\bpf  The first inequality in Eq.~\eqref{eq:sumnhk1} in fact holds for \emph{any} feasible solution of \flp, and it follows directly from   the constraints in Eq.~\eqref{eq:Lp1repconstr2} by noting that $D(h,l,k)$ is always no greater than $\old$. To show the second inequality in Eq.~\eqref{eq:sumnhk1}, note that since $D(h,l,k)$ is always no smaller than $\uld$, if we were to have $\frac{1}{\ln(1/\delta) }\sum_{k \in \calK}n^*_{h,k} >  \uld^{-1}  $, for some $h\in \calH$,  then we would be able to decrease one coordinate of $\{n^*_{h,k}\}_{k\in \calK}$ without violating any constraint, and hence the second inequality in Eq.~\eqref{eq:sumnhk1} must hold for some optimal solution of \flp.  Finally, to show Eq.~\eqref{eq:mstar}, we sum both sides of the the constraints in Eq.~\eqref{eq:Lp1repconstr1} over $\calK$, and obtain
\begin{align}
m^*_F = & \sum_{k\in \calK} r_k m^*_F \nln
\geq & \sum_{k\in \calK} \sum_{h\in \calH} n^*_{h,k}\pi_h  \nln
= &  \sum_{h\in \calH}  \pi_h \sum_{k\in \calK}  n^*_{h,k} \nln
\sk{a}{\geq} & \old^{-1} \ln(1/\delta) \sum_{h\in \calH} \pi_h \nln
=&  \old^{-1} \ln(1/\delta), 
\end{align}
where step $(a)$ follows from the first inequality in Eq.~\eqref{eq:sumnhk1}, which we have just shown. This proves Lemma \ref{lem:lp1structure}. 
\qed

\subsection{A Sufficient Condition for Accurate Classification}
\label{app:lem:SijSufficient2}

Fix $i, n \in \N$, and denote by $Y_{i,n} = \{(X_{i,j}, K_{i,j}) \}_{j = 1, \ldots, n}$, the {inspection history} of job $i$ up till the $n$th inspection it receives, where the inspections are ranked according to their time of initiation. Let $S_{i, [n]}(h,l)$ be the cumulative log-likelihood from the first $n$ inspections: $S_{i, [n]}(h,l) = \sum_{j=1}^n Z_{i,j}(h,l,K_{i,j})$, and similarly, let $\what{H}_{i,[n]}$ be the ML estimator of $H_i$ from the first $n$ inspections. 

\begin{lemma}
\label{lem:SijSufficient2}
Fix $i \in \N$ and $x>0$. Let $N^F_i$ be the total number of inspections job $i$ receives before its departure from the system. Let $N \in \N$, $N\leq N^F_i$, be a stopping time with respect to $Y_{i,N^F_i}$. Denote by $\calG_x$ the event:
 \begin{equation}
 \calG_x = \{\exists l \in \calH, \mbox{ s.t. } S_{i,[N]}(l, h') \geq x, \, \,  \forall h' \in \calH, h'\neq l \}.
 \end{equation}
 We have that
\begin{equation}
\pb(\what{H}_{i,[N]} \neq h \, , \, \calG_x\bbar H_i=h)\leq \szH\exp(-x), \quad \forall h \in \calH. 
\end{equation}
\end{lemma}
Note that Lemma \ref{lem:SijSufficient} follows from the above lemma by setting $N = N^F_i$. 

\bpf Fix $i \in \N$. Let $y_n = \{(x_j,k_j)\}_{j=1, \ldots, n}$ be a particular realization of $Y_{i,n}$. We will use the notation of $q_h(y_n)$, defined in Eq.~\eqref{eq:qhyndef} in Appendix \ref{app:lem:EofShl}, with a slight modification of changing $N_i$ to $N$, i.e.,
\begin{align}
q_h(y_n) = & \pb( N = n, Y_{i,n}=y_n \bbar H_i = h ).
\end{align}
Using a derivation that is identical to Eq.~\eqref{eq:fhy0} in Appendix \ref{app:lem:EofShl}, we have that
\begin{align}
q_h( y_n) =  v(y_n)\prod_{j=1}^n p(h,k_{j}, x_{j}). 
\label{eq:fhy}
\end{align}
where $v(y_n) \bydef \pb(N = n \bbar Y_{i,n}= y_n)\prod_{j=1}^n \pb(K_{i,j}= k_j,  N\geq j \bbar Y_{i,j-1}= y_{j-1})$, and it does not depend on $h$. We thus have  that
\begin{equation}
\frac{q_h(Y_{i,n})}{q_l(Y_{i,n})} = \frac{\prod_{j=1}^n p(h,K_{i,j}, X_{i,j})}{\prod_{j=1}^n p(l,K_{i,j}, X_{i,j})}= \exp(- S_{i, [n]}(l,h)). 
\label{eq:fynToSit}
\end{equation}

Fix $h, l \in \calH$, $h\neq l$. Denote by $\calC^l_n$ the event
\begin{equation}
\calC_n^l = \{N=n\, , \, \what{H}_{i,[n]} = l\}. 
\end{equation}
In particular,
\begin{equation}
\pb(\what{H}_{i,[N]} = l \, , \, \calG_x \bbar H_i = h) = \sum_{n=0}^\infty \pb(C^l_n \, , \, \calG_x \bbar H_i=h). 
\label{eq:pbwrongtoCnl}
\end{equation}
We now employ a change-of-measure argument. Denote by $\calY_n$ the set of realizations, $y_n$, for which $\pb(\calC^l_n , Y_{i,n} = y_n \bbar H_i = h) >0$, and 
\begin{equation}
S_{i,[n]}(l,h') \geq x, \quad \forall h' \in \calH, h' \neq l. 
\label{eq:calY2}
\end{equation}
We then obtain that
\begin{align}
\pb(C^l_n \, ,  \, \calG_x \bbar H_i =h)=& \sum_{y_n\in \calY_n} q_h( y_n) \nln
=& \sum_{ y_n \in \calY_n }\frac{q_h( y_n)}{q_l( y_n)}q_l( y_n) \nln
\stackrel{(a)}{=}  & \sum_{y_n\in \calY_n}\exp(-S_{i, [n]}(l,h))q_l( y_n)  \nln
\stackrel{(b)}{\leq}  & \exp(-x) \sum_{y_n \in \calY_n }q_l( y_n)  \nln
=  & \exp(-x) \pb(C^l_n \, ,  \, \calG_x\bbar H_i =l), 
\label{eq:pbCn^l}
\end{align}
where step $(a)$ is based on Eq.~\eqref{eq:fynToSit}, and step $(b)$ from Eq.~\eqref{eq:calY2}. Because Eq.~\eqref{eq:pbCn^l} holds for any $l \neq h$ and $n\in \N$, using a union bound, we have that
\begin{align}
\pb(\what{H}_i \neq h  \, ,  \, \calG_x \bbar H_i =h) = & \sum_{l \in \calH, l \neq h} \pb(\what{H}_i =l\, ,  \, \calG_x\bbar H_i =h ) \nln
\stackrel{(a)}{=} & \sum_{l \in \calH, l \neq h}\left(\sum_{n=0}^\infty \pb(C^l_n \, ,  \, \calG_x \bbar H_i =h)\right)  \nln
\stackrel{(b)}{\leq} &  \sum_{l \in \calH, l \neq h}\exp(-x)  \left(\sum_{n=0}^\infty \pb(C^l_n \, ,  \, \calG_x\bbar H_i =l)\right) \nln
=  & \exp(-x)  \sum_{l \in \calH, l \neq h}\pb(\what{H}_{i,[n]}=l \bbar H_i =l) \nln
\leq &  \szH \exp(-x) , 
\end{align}
where step $(a)$ follows from Eq.~\eqref{eq:pbwrongtoCnl}, and step $(b)$ from Eq.~\eqref{eq:pbCn^l}. This proves our claim.  \qed

\subsection{Proof of Proposition \ref{prop:erroAfterPrep}} 
\label{app:prop:erroAfterPrep}

\bpf Fix $h\in \calH$ and $i \in \N$. Denote by $S^P_i(h,l)$ the value of $S_{i,t}(h,l)$ at the time when job $i$ exits the Preparation stage. Define the event
\begin{equation}
\calB =  \lt\{S^P_i(h,l) \geq \ln\ln(1/\delta) , \, \forall l \in \calH \bslh \rt\}. 
\end{equation} 

The following lemma is the main technical result that we will use. Its proof makes use of the Azuma-Hoeffding ineuqality by noting that $S_{i,t}(h,l)$ is a martingle under proper conditioning. 

\begin{lemma}
\label{lem:reachSijPrep} 
Denote by $\overline{\calB}$ the complement of $\calB$. For all $h\in \calH$, 
\begin{equation}
\pb \lt(\overline{\calB} \bbar H_i = h \rt) \leq \szH \ln^{-1}(1/\delta). 
\end{equation}
\end{lemma}

\bpf We will index all inspections received by job $i$ in the Preparation stage in an arbitrary fashion, and denote by $K_j$ the type of the expert who performed the $j$th inspection for job $i$, and by $Z_{j}(\cdot,\cdot, K_j)$ the corresponding log-likelihood ratio. Define 
\begin{equation}
M^l_n = \sum_{s=1}^{n \vee n^P }  Z_{s}(h,l,K_s) - d(h,l), \quad n \in \N,
\end{equation}
where
\begin{equation}
d(h,l) = \sum_{k \in \calK} D(h, l, k) r_k. 
\end{equation}
Recall that the experts visit the three stages according to the randomized rule described in Section \ref{sec:expertVis},  and hence the probability of the $j$th expert to inspect job $i$ being type $k$ is equal to $r_k$. Conditional on the true label of job $i$ being $h$, we have that the $Z_{s}(h,l,K_s)$'s are i.i.d., with $•\E(Z_1(h,l,K_1) | H_i = h)=d(h,l)$.  It is hence not difficult to verify that $\{M^l_n \bbar H_i = h\}_{n \in \N}$ is a martingale. Recall that  $\dav \bydef \min_{h, l \in \calH, h\neq l}d(h,l)$. We have that for any $l \in \calH\bslh$, 
\begin{align}
\pb(S^P_i(h,l)) \leq \ln\ln(1/\delta)  \bbar H_i = h) =& \pb(S^P_i(h,l) - n^P d(h,l)\leq \ln\ln(1/\delta) - n^P d(h,l) \bbar H_i = h) \nln 
\leq& \pb(S^P_i(h,l) - n^P d(h,l)\leq \ln\ln(1/\delta) - n^P \dav \bbar H_i = h) \nln 
= & \pb(S^P_i(h,l) - n^P d(h,l)\leq  - (\dav- \ln\ln(1/\delta)/n^P) n^P \bbar H_i = h) \nln 
\sk{a}{\leq} & \exp\lt(- \frac{(\dav- \ln\ln(1/\delta)/n^P)^2 (n^P)^2 }{8 \olz^2 n^P}\rt) \nln
\leq & \exp\lt(- \frac{\dav^2 n^P- 2\dav  \ln\ln(1/\delta) }{8 \olz^2 }\rt) \nln
\sk{b}{=} & \exp\lt(- \ln\ln(1/\delta) \rt)  \nln
=& \ln^{-1}(1/\delta). 
\label{eq:sphl1}
\end{align}
Step $(a)$ follows from the Azuma-Hoeffding's Inequality (Lemma \ref{lem:asuma} in Appendix \ref{app:azuma}), and that $|M^l_n-M^l_{n-1}| \leq 2\olz$ for all $n\in \N$. Step $(b)$ follows from the fact that
\begin{equation}
n^P = \zeta_0  \ln\ln(1/\delta)= \frac{8\olz^2 + 2 \dav }{\dav^2}\ln\ln(1/\delta). 
\end{equation}
Because Eq.~\eqref{eq:sphl1} holds for all $l$, Lemma \ref{lem:reachSijPrep} follows by applying a union bound over $l \in \calH$. \qed

We now return to the proof of Proposition \ref{prop:erroAfterPrep}.  Combining Lemmas \ref{lem:SijSufficient} and \ref{lem:reachSijPrep}, we have that
\begin{align}
\pb(\what{H}_i^P \neq h \bbar H_i = h ) \leq&  \pb(\overline{\calB} \bbar H_i = h  )+ \pb(\what{H}_i^P \neq h \, ,\, \calB  \bbar H_i = h) \nln
\sk{a}{\leq} &   \szH\ln^{-1}(1/\delta) + \pb(\what{H}_i^P \neq h \, ,\, \calB  \bbar H_i = h) \nln
\sk{b}{\leq} &  \szH\ln^{-1}(1/\delta) + \szH\exp(-\ln\ln(1/\delta)) \nln
= & 2 \szH\ln^{-1}(1/\delta), 
\label{eq:Hpi}
\end{align}
where step $(a)$ follows from Lemma \ref{lem:reachSijPrep}, and  $(b)$ from Lemma \ref{lem:SijSufficient2} by setting $x = \ln\ln(1/\delta)$, $N = n^P$, and noting that $\calB \subset \calG_x$. Our claim that $\epsilon^P  \leq 2  \szH\ln^{-1}(1/\delta)$ follows by noting that the above equation holds for all $h\in \calH$. To show Eq.~\eqref{eq:piP}, note that, for all $h\in \calH$, 
\begin{align}
\pi^P_h \leq & \pi_h + \sum_{h' \neq h} \pi_{h'} \pb(\what{H}_1^P = h \bbar H_1 = h') \nln
\leq & \pi_h + \epsilon^P\sum_{h' \neq h} \pi_{h'} \nln
\leq& \pi_h + 2  \szH\ln^{-1}(1/\delta). 
\end{align}
where the last inequality follows from the upper bound on $\epsilon^P$, and that $\sum_{h' \neq h} \pi_{h'} \leq 1$. This completes the proof of Proposition \ref{prop:erroAfterPrep}.  \qed

\subsection{Proof of Proposition \ref{prop:subseqConvrg}} 
\label{app:prop:subseqConvrg}

\bpf We start by looking at the workload process for the Preparation stage, $\buw^s_0(t)$.  For every $a,b>0$, with $a>b$, $s\in \N$, and $z>0$, applying the triangle inequality to Eq.~\eqref{eq:Ws0v1}, we have that
\begin{align}
|\buw^s_0(z, a) - \buw^s_0(z, b)| \leq  & |\bua^s_0(z,a)- \bua^s_0(z,b)| + |\budta^s_0(z,a)-\budta^s_0(z,b)| \nln
\sk{a}{=} & n^P|\Xi_0(z,a)-\Xi_0(z,b)|+|\budta^s_0(z,a)-\budta^s_0(z,b)| \nln
\sk{b}{\leq} &  n^P|\Xi_0(z,a)-\Xi_0(z,b)| + \sum_{e\in \calE} |R_{e,P}(z,a)-R_{e,P}(z,b)|. 
\label{eq:wstriangle}
\end{align}
where  $\Xi_0(\cdot)$ and $R_{e,T}(\cdot)$ are defined in Eqs.~\eqref{eq:regXi} and \eqref{eq:regRAP}, respectively.  Step $(a)$ follows from the design of our inspection policy, where each job that arrives to the Preparation stage adds $n^P$ new inspections to $\buw_0(t)$. Step $(b)$ is due to the fact that the number of initiated inspections during an interval is no greater than the total number of times that an expert visits the Preparation stage during the same period. 

Recall that the average inspection rate, $\overline{\mu}$, is equal to $1$ (Eq.~\eqref{eq:avgRate}).  Because $\omega\in \calC$ and $\lim_{s\to \infty}z_s = \infty$, by Eqs.~\eqref{eq:regXi}, \eqref{eq:regRAP} and \eqref{eq:wstriangle}, we have that there exists a positive sequence $\{\epsilon_s\}_{s\in \N}$, with $\epsilon_s\to 0$ as $s\to \infty$, such that for all $s\in \N \mbox{ and } a,b \in [0,T], a>b$, 
\begin{align}
\label{eq:as0tight}
|\bua^s_0(z,a)- \bua^s_0(z,b)|  \leq & n^p(a-b)+ \epsilon_s, \\
|\budta^s_0(z,a)- \budta^s_0(z,b)|  \leq & q^P m (a-b)+ \epsilon_s, \\
|\buw^{s}_0(z_s,a) - \buw^{s}_0(z_s,a)| \leq & ( n^P+q^P m )(a-b) + \epsilon_{s}. 
\label{eq:ws0tight}
\end{align}
Recall that $\buw_0(0) = \blw^{(s)}$ in the $s$th system, and $\lim_{s\to \infty}{z_s}^{-1}\blw^{(s)} = \blw^0$. We invoke the following variation of the Arzel\`{a}-Ascoli theorem (Lemma 6.3 of \cite{ye2005stability}). 
\begin{lemma}
\label{lem:Arzela-AscoliVar}
Fix $T>0$. Let $\{f_s(\cdot)\}_{s \in \N}$ a sequence of functions that satisfies
\begin{enumerate}
\item $\sup_{s \in \N} |f_s(0)| < \infty $. 
\item There exists a positive sequence $\{\epsilon_s\}_{s\in \N}$ with $\epsilon_s \to 0$ as $ s \to \infty$, such that 
\begin{equation}
|f_s(a)-f_s(b)|\leq c|a-b| + \epsilon_s, \quad \forall a,b\in [0,T], \, s\in \N. 
\end{equation}
\end{enumerate}
Then, there exists a $c$-Lipschitz function $f^*(\cdot)$, and an increasing sequence $\{s_j\}_{j\in \N} \subset \N$, such that
\begin{equation}
\lim_{j \to \infty} \snorm{f_{s_j}(t) - f^*(t)}  = 0. 
\end{equation}
\end{lemma}
Using Lemma \ref{lem:Arzela-AscoliVar} and Eqs.~\eqref{eq:as0tight} through \eqref{eq:ws0tight}, we have that there exist $(n^P+q^P m )$-Lipschitz functions, $\bla_0(\cdot)$, $\bld_0(\cdot)$, and  $\blw_0(\cdot)$, and an increasing sequence $\{i^0_j\}_{j\in \N}$, such that
\begin{equation}
\lim_{j \to \infty} \snorm{ (\buw_0, \bua_0, \budta_0)^{s_{i^0_j}} (z_{s_{i^0_j}}, t) - (\blw_0, \bla_0, \bld_0)(t)} = 0. 
\label{eq:w0convg}
\end{equation}
This shows that the sequence $\{(\buw_0, \bua_0, \budta_0)^{s_i}(z_{s_i}, \cdot)\}_{i \in \N}$ admits Lipschitz-continuous limit points. 

We now repeat essentially the same argument for the workload in the Adaptive stage, and use a diagonal argument to identify a family of nested sub-sequences along which the workload converges coordinate-wise to a Lipschitz function. Fix $k\in \{1, \ldots, \calK\}$. We claim that
\begin{align}
\lt|\bua^s_k(a)-\bua^s_k(b) \rt |\leq& v_\delta\left( m+ \sum_{e\in \calE}|R_{e,P}(a)-R_{e,P}(b)| \rt), \quad \forall a,b \in [0,T], a>b. 
\label{eq:atoR1}
\end{align}
To justify the above inequality, note that by Eq.~\eqref{eq:calNh3}, each job that arrives to the Adaptive stage could add at most $v_\delta$ inspections to the aggregate workload in the Adpative stage. Also, an arrival to the Adaptive stage corresponds to a departure to the Preparation stage, and the number of such departures during an interval is at most the number of times an expert visits the Preparation stage, $\sum_{e\in \calE}|R_{e,P}(a)-R_{e,P}(b)| $, plus the number of experts who were already inspecting a job in the Preparation stage prior to time $t=a$, which is at most $m$. This yields Eq.~\eqref{eq:atoR1}. Analogous to Eq.~\eqref{eq:wstriangle}, we have that
\begin{align}
|\buw^s_k(z,a) - \buw^s_0(z,b)| \leq& n^P|\bua^s_k(z,a)-\bua^s_k(z,b)|+|\budta^s_k(z,a)-\budta^s_k(z,b)| \nln
\leq&  \gamma_\delta\left( \frac{m}{z}+ \sum_{e\in \calE}|R_{e,P}(z,a)-R_{e,P}(z,b)| \rt) + \sum_{e\in \calE} |R_{e,A}(z,a)-R_{e,A}(z,b)|. 
\label{eq:wstriangle2}
\end{align}
Using the fact that $\omega \in \calC$, $\lim_{s\to \infty}z_s = \infty$, and Eqs.~\eqref{eq:regRAP} and \eqref{eq:wstriangle2}, we have that there exists a positive sequence $\{\epsilon_s\}_{s\in \N}$, with $\lim_{s\to \infty}\epsilon_s=0$, such that, for all $s\in \N \mbox{ and } a,b \in [0,T], a>b$, 
\begin{align}
|\bua^s_k(z,a)- \bua^s_k(z,b)|  \leq & v_\delta q^P m (a-b)+ \epsilon_s, \\
|\budta^s_k(z,a)- \budta^s_k(z,b)|  \leq & q^A m (a-b)+ \epsilon_s, \\
|\buw^{s}_k(z_s,a) - \buw^{s}_k(z_s,b)| \leq& (v_\delta q^P+q^A )m(a-b) + \epsilon_{s}.  
\label{eq:wsotherktight}
\end{align}
Using Eq.~\eqref{eq:wsotherktight} and  the same argument involving Lemma \ref{lem:Arzela-AscoliVar} as before, we have that there exist 
\begin{enumerate}
\item a collection of coordinate-wise $(v_\delta q^P+q^A )m$-Lipschitz functions, $\{(\blw_k, \bla_k, \bld_k)\}_{k \in \calK}$, and 
\item a family of sequences, $\{ \{i^k_{j}\}_{j\in \N}\}_{k \in \calK }$, where $\{i^k_{j}\}$ is a sub-sequence of $\{i^{k-1}_{j}\}$ for all $k\in \calK$, and $\{i^0_j\}$ was defined in Eq.~\eqref{eq:w0convg}, 
\end{enumerate}
such that
\begin{equation}
\lim_{j\to \infty} \snorm{ (\buw_k, \bua_k, \budta_k) ^{s_{i^k_j}}(z_{s_{i^k_j}}, t) - (\blw_k, \bla_k, \bld_k)(t)} = 0, \quad \forall k \in \{1,\ldots,  \szK\}. 
\end{equation}
Because the sequences, $\{i^{k}_{j}\}_{k\in \calK}$, are nested, by setting $k= \szK$, the above equation further implies that
\begin{equation}
\lim_{j\to \infty} \snorm{(\buw, \bua, \budta)^{s_{i^{ \szK}_j}}(z_{s_{i^k_j}}, t) - (\blw, \bla, \bld) (t)} = 0,
\end{equation}
Note that, by construction, all coordinates of $(\blw, \bla, \bld)$ are Lipschitz-continuous, with a Lipschitz constant of $c =  [q^A+(v_\delta+1)q^P]m+ n^P$. We complete the proof of Proposition \ref{prop:subseqConvrg} by setting $\{i_j\}_{j\in \N} = \{i^{ \szK}_j \}_{j \in \N}$. \qed

\subsection{Proof of Lemma \ref{lem:Bhdiff}}
\label{app:lem:Bhdiff}

\bpf The general strategy of our proof is to relate the behavior of $B^s_h(\cdot)$ to the number of initiated inspections in the Preparation stage, $\Delta^s_0(\cdot)$, while being cautious of the fact that the lengths of inspections are random and can cause the order of job departures from the Preparation stage to be different from that of the arrivals. Let $B^s(t) = |\calB^s(t)|$. Fix $t\in (0,T)$ and $t'\in (t,T)$. We first show that the following property holds: 
\begin{equation}
\left|  (B^s(t') - B^s(t)) - (\Delta^s_0(t')- \Delta^s_0(t))/n^P \right| \leq 2m.
\label{eq:BtoDelta}
\end{equation}
In words, this means that the number of jobs that enter the Adaptive stage during the interval $(t,t']$ multiplied by $n^P$ deviates from the number of initiated inspections during the same period by at most $2mn^P$. 

To prove Eq.~\eqref{eq:BtoDelta}, recall that each job must receive $n^P$ inspections in the Preparation stage before it moves to the Adaptive stage, and hence having  $(\Delta^s_0(t')- \Delta^s_0(t))$ initiated inspections means that there could be at most $\lceil (\Delta^s_0(t')- \Delta^s_0(t))/n^P \rceil +m$ jobs arriving at the Adaptive stage during $(t, t']$, where the addition of $m$ captures the possibility that all $m$ servers were processing a job from the Preparation stage at the start of the interval, $t$. This shows that
\begin{equation}
 (B^s(t') - B^s(t)) - (\Delta^s_0(t')- \Delta^s_0(t))/n^P \leq  m+1 \leq 2m.
 \label{eq:bdifftodelta1}
\end{equation}
\del{Similarly, since there are $m$ experts and hence at most $m$ on-going inspections at any moment in time, the number of inspections \emph{completed} during the interval $(t,t']$ is at least $(\Delta^s_0(t')- \Delta^s_0(t)) - m$.}
We next invoke the following fact: 
\begin{lemma}
\label{fac:ballbinsm}
Fix $x\in \N$ and an interval $J \subset \rp$. If there are $x$ inspections initiated in the Preparation stage during $J$, then there are at least $(x/n^P  -m)$ jobs departing from the Preparation stage during $J$. 
\end{lemma}

\bpf Denote by ${\calI}'$ the set of jobs who have had  \emph{any} inspection initiated during the interval $J$. We partition ${\calI}'$ into ${\calI}'_0 \cup {\calI}'_1$, where $\calI'_0$ corresponds to those jobs who have departed from the Preparation stage by the end of $J$, and $\calI'_1$ to those who have not. Because the experts perform inspections in a first-come-first-server manner, a job cannot start receiving inspections until all $n^P$ inspections for the previous job have been initiated. Since there are $m$ experts, this implies that at any point in time, there can be at most $m$ jobs in the Preparation stage who have initiated any inspection, and therefore
 \begin{equation}
|\calI'_1| \leq m. 
\label{eq:calI1big}
\end{equation}
Recall that $x$ is the total number of initiated inspections in the Preparation stage during $J$. We have that
\begin{equation}
x \sk{a}{\leq} n^P|\calI'| =  n^P(|\calI'_0|+|\calI'_1|) \sk{b}{\leq}  n^P(|\calI'_0|+m),
\end{equation}
where step $(a)$ follows from the fact that each job can lead to at most $n^P$ inspections, and step $(b)$ from Eq.~\eqref{eq:calI1big}. This yields
\begin{equation}
|\calI'_0| \geq x/n^P - m. 
\end{equation}
That is, there are at least $(x/n^P - m)$ jobs departing from the Preparation stage during $J$. \qed

We now apply Lemma \ref{fac:ballbinsm} with $x=(\Delta^s_0(t')- \Delta^s_0(t))$, and obtain

\begin{equation}
B^s(t') - B^s(t) \geq (\Delta^s_0(t')- \Delta^s_0(t))/n^P  - m, 
\end{equation}
or
\begin{equation}
 (B^s(t') - B^s(t)) - (\Delta^s_0(t')- \Delta^s_0(t))/n^P \geq -m, 
 \end{equation} 
which, combined with Eq.~\eqref{eq:bdifftodelta1} leads to Eq.~\eqref{eq:BtoDelta}. 


To prove the main claim of Lemma \ref{lem:Bhdiff}, we will invoke the functional law of large number on the sequence of ML estimators,  $\{\what{H}^P_i\}_{i \in \N}$, as expressed in Eq.~\eqref{eq:regHhatP}. However, we have to be careful in doing so: the jobs that depart from the Preparation stage do not necessarily preserve the order in which they arrived to the system, since the lengths of inspections in the Preparation stage are random and can potentially cause a job to depart earlier than another with a smaller index. Thus, we  begin by showing that such ``shuffling'' is relatively minor and can essentially be ignored when studying the local behavior of $B^s_h(\cdot)$. This is formalized in the next lemma, which states that the set $\calB^s(t') \backslash \calB^s(t)  $ can be ``sandwiched'' by two intervals of consecutive integers that differ by no more than $2m$ elements from one another. The use of consecutive integers will allow us to invoke the law of large numbers with greater ease. 

We will adopt the following short-hand notation. For a set-valued process $\{\calS(t)\}_{t\in \rp}$, and $t,t' \in \rp$, we will denote by $\calS(t' \backslash t) $ the difference between $\calS(t')$ and $\calS(t)$: 
\begin{equation}
\calS(t' \backslash t) \bydef \calS(t')\backslash \calS(t). 
\end{equation}
For $a,b \in \N$, we will use $\{a\to b\}$ to denote the set of consecutive integers $\{a, a+1, \ldots, b\}$, if $a\leq b$, and $\emptyset$, otherwise. We have the following lemma. 
\begin{lemma}
\label{lem:orderpresB}
Fix $t,t'\in \rp$, $t'>t$. There exist $a, a', b$ and $b' \in \N$, with $\max\{|a-a'|, |b-b'|\}\leq m$,  such that
\begin{equation}
\{a \to b\} \subset \calB^s(t' \backslash t)  \subset \{a' \to b'\}. 
\label{eq:orderpresB1}
\end{equation}
It follows that
\begin{equation}
\big| (b-a) - |\calB^s(t' \backslash t)| \big| \leq 2m, \mbox{ and }  \big| (b'-a')  - |\calB^s(t' \backslash t)| \big| \leq 2m. 
\label{eq:orderpresB2}
\end{equation}
\end{lemma}

\bpf Denote by $\calD^s(t)$ the set of jobs in the Preparation stage for whom at least one inspection has been initiated by time $t$. Because a job can depart from the Preparation stage only after completing $n^P$ inspections, we have
\begin{equation}
\calB^s(t) \subset \calD^s(t), \quad \forall t\in \rp. 
\label{eq:calBinD}
\end{equation}
On the other hand, the initiations of inspections in the Preparation stage are performed in a first-come-first-serve (FCFS) manner. Recall also that at any given point in time there can be at most $m$ experts performing inspections. Combining the above two facts, we know that an inspection for job $i+m$ can be initiated only after all jobs prior to, and including, job $i$ have departed from the Preparation stage. We thus have that
\begin{equation}
\{1\to \max \calD^s(t) - m\} \subset \calB^s(t),\quad \forall t\in \rp. 
\label{eq:calDinB}
\end{equation}
where $\max \calD^s(t)$ is the index of the last job in $\calD^s(t)$. Combining Eqs.~\eqref{eq:calBinD} and \eqref{eq:calDinB}, we have 
\begin{align}
\label{eq:dtoBtod1}
\calB^s(t' \backslash t)   \supset & \ \{1\to \max \calD^s(t') - m\}  \backslash \calD^s(t), \\
\calB^s(t' \backslash t)   \subset & \ \calD^s(t') \backslash \{1\to \max \calD^s(t) - m\} 
\label{eq:dtoBtod2}
\end{align}
Finally, the above-mentioned FCFS property of the inspection initiation rule further implies that  $\calD^s(t)$  is a set of consecutive positive integers. Therefore, it follows that the right-hand sides of Eqs.~\eqref{eq:dtoBtod1} and \eqref{eq:dtoBtod2} are both intervals of consecutive integers. We can prove Eq.~\eqref{eq:orderpresB1} by setting
\begin{align}
\{a\to b\} = & \{1\to \max \calD^s(t') - m\}   \, \backslash \,  \calD^s(t), \\
\{a' \to b'\} = & \calD^s(t') \, \backslash \, \{1\to \max \calD^s(t) - m\}, 
\end{align}
or, equivalently, that
\begin{align}
a =& \max{\calD^s(t)}+1, \quad a' = \max{\calD^s(t)}+1-m  \nln
b = & \max{\calD^s(t')}-m+1, \quad b' = \max{\calD^s(t')}+1
\end{align}
and it is clear from the above equations that $\max\{|a-a'|,|b-b'|\}\leq m$.  For Eq.~\eqref{eq:orderpresB2}, note that by  the ordering in Eq.~\eqref{eq:orderpresB1} and the triangle inequality, we have
\begin{equation}
\big| (b-a) - |\calB^s(t' \backslash t)| |  \leq \big| (b'-a') - (b-a) |  \leq 2\max\{|a-a'|, |b-b'|\}\leq 2m, 
\end{equation}
and a similar argument shows $\big| (b'-a') - |\calB^s(t' \backslash t)| | \leq 2m$. This proves Lemma \ref{lem:orderpresB}. 
\qed

We are now ready to prove the main claim of Lemma \ref{lem:Bhdiff}. Fix $t \in (0,T)$ and $\epsilon >0$. Let $a,a',b$ and $b'$ be defined as in Lemma \ref{lem:orderpresB}, corresponding to the set $ \calB^{\overline{s}_j}( \overline{z}_j (t+\epsilon) \backslash \overline{z}_j t)$, all of which depend on $j$ and $\epsilon$, though we will suppress the dependencies in our notation for simplicity. We have that, for every sufficiently small $\delta'>0$,  there exists   $j^*>0$, such that for all $j>j^*$, 
\begin{align}
\frac{B_h^{\bar{s}_j}(\bar{z}_j, t+\epsilon)- B_h^{\bar{s}_j}(z_{\bar{s}_j}, t)}{\epsilon} =& \frac{B_h^{\bar{s}_j}(\bar{z}_j(t+\epsilon)) - B_h^{\bar{s}_j}(z_{\bar{s}_j}t)}{\epsilon \overline{z}_j}  \nln
= & \frac{1}{\epsilon \overline{z}_j}\sum_{i \in\calB^{\overline{s}_j}( \overline{z}_j  (t+\epsilon)\backslash \overline{z}_j t)  } \mathbb{I}(\what{H}^P_i = h) \sk{a}{\geq} \frac{1}{\epsilon \overline{z}_j}\sum_{i =a }^b \mathbb{I}(\what{H}^P_i = h)  \nln
\sk{b}{\geq} & \frac{1}{\epsilon \overline{z}_j}\sum_{i =a }^{a + (B^{\overline{s}_j}(\overline{z}_j (t+\epsilon)- B^{\overline{s}_j}(\overline{z}_j t )) -2m} \mathbb{I}(\what{H}^P_i = h) \nln
\sk{c}{\geq} & \frac{1}{\epsilon \overline{z}_j}\sum_{i =a }^{a + (\budta^{\overline{s}_j}_0(\overline{z}_j (t+\epsilon)- \budta^{\overline{s}_j}_0(\overline{z}_j t ))/n^P -4m} \mathbb{I}(\what{H}^P_i = h) \nln
\sk{d}{\geq} & \frac{1}{\epsilon \overline{z}_j} \sum_{i =a }^{a + (\bld_0(t+\epsilon)-\bld_0(t) -\delta)\overline{z}_j/n^P -4m} \mathbb{I}(\what{H}^P_i = h) \nln
\sk{e}{\geq} & \frac{\pi^P_h}{\epsilon n^P }  (\bld_0(t+\epsilon)-\bld_0(t) -2\delta') \nln
=&\frac{\pi^P_h}{ n^P }  \cdot \frac{\bld_0(t+\epsilon)-\bld_0(t) -2\delta'}{\epsilon}.
\label{eq:limitOfBh1} 
\end{align}
Step $(a)$ follows from Eq.~\eqref{eq:orderpresB1} of Lemma \ref{lem:orderpresB}, and step $(b)$ from the first inequality in Eq.~\eqref{eq:orderpresB2} in the same lemma. Step $(c)$ is a consequence of $\left|  (B^s(t') - B^s(t)) - (\budta^s_0(t')- \budta^s_0(t))/n^P \right| \leq 2m $ (Eq.~\eqref{eq:BtoDelta}). Step $(d)$ follows from the uniform convergence of $\budta^{\overline{s}_j}_0(\cdot) $ to $\bld_0(\cdot)$ over $[0,T]$. Finally, step $(e)$ is based on the functional law of large number for the sequence $\{\what{H}^P_i\}_{i \in \N}$ (Eq.~\eqref{eq:regHhatP}). 

Since Eq.~\eqref{eq:limitOfBh1} holds for all sufficiently small $\delta'>0$, we have that
\begin{align}
\lim_{\epsilon \downarrow 0} \lim_{j \to \infty} \frac{B_h^{\bar{s}_j}(\bar{z}_j, t+\epsilon)- B_h^{\bar{s}_j}(z_{\bar{s}_j}, t)}{\epsilon}\geq & \lim_{\epsilon \downarrow 0} \lim_{j \to \infty} \lim_{\delta' \downarrow 0} \frac{\pi^P_h}{ n^P }  \cdot \frac{\bld_0(t+\epsilon)-\bld_0(t) -2\delta'}{\epsilon} \nln
= & \frac{\pi^P_h}{ n^P }  \dot{\bld}_0(t) \nln
= &  \left\{ \begin{array}{ll}
          \frac{m q^P}{n^P}  \pi^P_h , & \quad \mbox{if } \blw_0(t)>0,\\
          \\
          \pi^P_h , & \quad  \mbox{if } \blw_0(t)=0,\\
         \end{array} \right.
         \label{eq:limitOfBh2}
\end{align}
A line of argument identical to Eqs.~\eqref{eq:limitOfBh1} through \eqref{eq:limitOfBh2} shows the other direction of the inequality, namely, that $\lim_{\epsilon \downarrow 0} \lim_{j \to \infty} \frac{B_h^{\bar{s}_j}(\bar{z}_j, t+\epsilon)- B_h^{\bar{s}_j}(z_{\bar{s}_j}, t)}{\epsilon} \leq  \frac{\pi^P_h}{ n^P }  \dot{\bld}_0(t)$. This completes the proof of Lemma \ref{lem:Bhdiff}. \qed

\subsection{Proof of Lemma \ref{lem:uppersemicont_Nstar}}
\label{app:lem:uppersemicont_Nstar}

\bpf We first state a useful technical lemma. The proof makes use of Berge's maximum theorem concerning the continuity of optimal solutions to a convex optimization problem. 
\begin{lemma}
\label{lem:uppersimconNh}
Fix $h\in \calH$ and $\blw\in \rp^{ \szK+1}$. The set-valued function, $\calN_h^*(\cdot) $, satisfies the following semi-continuity property. For every $\delta>0$, there exists $\epsilon>0$ such that, for all $\blw'\in \rp^{ \szK+1}$,  $\|\blw'-\blw\|_2 \leq \epsilon$, 
\begin{equation}
\sup_{x\in \calN_h^*(\blw')} \ \inf_{y \in \calN_h^*(\blw)} \|x-y\|_2 \leq \delta. 
\end{equation}
\end{lemma}

\bpf The set $\calN_h(\cdot)$ represents the set of feasible solutions of the linear optimization problem defined in Eq.~\eqref{eq:argminLinearWeight}. We first invoke Berge's maximum theorem (pp.~116 of \cite{berge1963topological}), which roughly states that $\calN_h(\cdot)$ depends semi-continuously in its argument. More specifically, it is easy to verify that the set,  $\calN_h$, is compact (Eq.~\eqref{eq:calNh1} to \eqref{eq:calNh3}), and the objective function is continuous. Berge's maximum theorem thus implies that $\calN_h^*(\cdot)$ is upper-hemicontinuous at $\blw$, in the following sense: let  $\{\blw_n\}_{n\in \N}$ be a sequence, where $\blw_n \to \blw$ as $n\to \infty$, and $\{\lambda_n\}_{n \in \N}$ be a sequence where $\lambda_n \in \calN^*_h(\blw^n)$ for all $n\in \N$. Then, if $\lambda_n\to \lambda$ as $n\to \infty$, we have that $\lambda\in \calN^*_h(\blw)$. 

We now use the above upper hemicontinuity property, along with the boundedness of the constraint set, $\calN_h$, to prove our claim.  Suppose, for the sake of contradiction,  that there exists a sequence $\{\blw_n\}_{n\in \N}$ with $\lim_{n\to \infty}\blw_n= \blw$, such that
\begin{equation}
\liminf_{n \to \infty} \sup_{x\in \calN_h^*(\blw_n)} \ \inf_{y \in \calN_h^*(\blw)} \|x-y\|_2 > 0. 
\end{equation}
This implies the existence of a sequence $\{\lambda_n\}_{n \in \N}$, with $\lambda_n \in \calN^*_h(\blw^n)$ for all $n\in \N$, such that
\begin{equation}
\liminf_{n \to \infty} \inf_{y \in \calN_h^*(\blw)} \|\lambda_n-y\|_2 > 0. 
\label{eq:lamnotconv1}
\end{equation}
By definition, for every $\blw$,  the coordinates of the  elements of $\calN^*_h(\blw)$ are non-negative and bounded from above by $v_\delta$. Therefore, by sequential compactness, there exist $\lambda \in \rp^{ \szK+1}$ and a sub-sequence of $\{\lambda_{n}\}$, $\{\lambda_{n_k}\}_{k\in \N}$, such that $\lim_{k \to \infty} \lambda_{n_k} = \lambda$. By Eq.~\eqref{eq:lamnotconv1}, we have that
\begin{equation}
\inf_{y \in \calN_h^*(\blw)} \|\lambda - y\|_2 > 0.
\end{equation}
This contradicts with the hemicontinuity of $\calN^*_h(\cdot)$, which would imply that $\lambda \in \calN^*_h(\blw)$, and thus proves Lemma \ref{lem:uppersimconNh}.  \qed

We now return to the proof of Lemma \ref{lem:uppersemicont_Nstar}. Fix $\epsilon \in (0, T-t)$. 
We first observe that the set $\calN^*_h(\cdot)$ is scale-invariant, in the sense that for every $\blw \in \rp^{ \szK+1}$
\begin{equation}
\calN^*_h ( a \blw) = \calN^*_h(\blw), \quad a >0. 
\end{equation}
Therefore, 
\begin{equation}
\calN^*_h (\buw^s(zt) )= \calN^*_h (z^{-1}\buw^s(zt) )= \calN^*_h (\buw^s(z,t) ), \quad  
\label{eq:Nhscalinvar}
\end{equation}

Recall that $\buw^{s_{j}}( \bar{z}_j , \cdot)$ converges uniformly over $[0,T]$ to $\blw(\cdot)$ as $j\to \infty$, and that $\blw(\cdot)$ is a continuous function. These two facts, along with the semi-continuity property of Lemma \ref{lem:uppersimconNh}, imply that for every $\delta>0$, there exist $\epsilon, j^* >0$, such that for all $j\geq j^*$
\begin{equation}
\inf_{y \in \calN_h^*(\blw(t))} \| \lambda - y\|_2 \leq \delta, \quad \forall t' \in (t, t+\epsilon] , \, \lambda\in \calN_h^*(\buw^{\bar{s}_j}( \bar{z}_j , t')). 
\label{eq:infyinNh1}
\end{equation}
Note that all jobs in the set $
 \calB^{\bar{s}_j}_h( \bar{z}_j (t+\epsilon) \backslash \bar{z}_j t)$ arrived during the interval $( \bar{z}_j t,  \bar{z}_j (t+\epsilon)]$. Therefore, for all $i\in 
 \calB^{\bar{s}_j}_h( \bar{z}_j (t+\epsilon) \backslash \bar{z}_j t)$, there exists $t' \in (t,t+\epsilon]$, such that 
\begin{equation}
\Lambda_i \leq \calN^*_h(\buw^{\bar{s}_j}( \bar{z}_j , t')),
\label{eq:laminNcal}
\end{equation}
where the inequality follows from round-down procedure in Eq.~\eqref{eq:Lamb_i_k}. 

Combining Eqs.~\eqref{eq:infyinNh1} and \eqref{eq:laminNcal}, we conclude that for every $\delta>0$,  there exist $\epsilon^*, j^*>0$, such that for all $j\geq j^*, \epsilon<\epsilon^*$, 
\begin{equation}
\inf_{y \in  \calN^*_h(\blw(t))} 
\max_{k \in \calK} \, \lt( \Lambda_{i,k} - y_k \rt)  \leq \delta, \quad \forall i \in 
 \calB^{\bar{s}_j}_h( \bar{z}_j (t+\epsilon) \backslash \bar{z}_j t). 
\end{equation}
Since the above inequality holds for \emph{all} jobs in $ \calB^{\bar{s}_j}_h( \bar{z}_j (t+\epsilon) \backslash \bar{z}_j t)$, it further implies that the average workload among the jobs in $ \calB^{\bar{s}_j}_h( \bar{z}_j (t+\epsilon) \backslash \bar{z}_j t)$ satisfies: 
\begin{equation}
\inf_{y \in \calN^*_h(\blw(t))}  \, \max_{k \in \calK} \, \lt( \overline{\Lambda}_{h,k}(t, \epsilon, j) - y_k \rt)   \leq \delta. 
\end{equation}
Because the above inequalities hold for all $\delta>0$, we conclude that
\begin{equation}
\limsup_{\epsilon \downarrow 0}  \limsup_{j \to \infty} \inf_{y \leq \calN_h^*(\blw(t))}\|\overline{\Lambda}_h(t, \epsilon, j)- y\|_2 =0. 
\end{equation}
This completes the proof of Lemma \ref{lem:uppersemicont_Nstar}. \qed

\subsection{Proof of Proposition \ref{prop:drainFluid}}
\label{app:prop:drainFluid}
\bpf  Fix $\blw^0\in \rp^{ \szK+1}$ such that $L(\blw^0)=1$, and a fluid solution $\blw\in \calW(\blw^0)$. Fix $t\in (0,T)$ to be a point where all coordinates of $\blw(\cdot)$, $\bla(\cdot)$, and $\bld(\cdot)$ are differentiable. We have, by the chain rule of differentiation,
\begin{align}
\frac{d}{dt} L(\blw(t)) = &\frac{d}{dt} \sqrt{ \sum_{k = 0}^{ \szK} \blw_k(t)^2 } \nln
 = & \|\blw(t)\|_2^{-1/2} \lt(  \sum_{k = 0}^{ \szK} \blw_k(t)\dot{\blw}_k(t) \rt) \nln
 = & \|\blw(t)\|_2^{-1/2}  \lt(\blw_0(t)\dot \blw_0(t) + \sum_{k \in \calK} \blw_k(t)\dot{\blw}_k(t) \rt). 
 \label{eq:dLdt1}
\end{align}
We next inspect separately at the two terms in the parentheses in Eq.~\eqref{eq:dLdt1}. For the first term, corresponding to the workload in the Preparation stage, we have
\begin{equation}
\blw_0(t) \dot \blw _0(t) = \blw_0(t) (\dot{\bla}_0(t) - \dot{\bld}_0(t)) \sk{a}{=} \blw_0(t)(n^P-m q^P) = - \blw_0(t)\tilde{c},  
 \label{eq:dLdt-w0term}
\end{equation}
where we define $\tilde{c} \bydef mq^P-n^P >0$, and step $(a)$ follows from the definition of a fluid solution.  

We now analyze the second term, which corresponds to the workloads in the Adaptive stage. Recall that $\calN_h$ is the set of vector satisfying  Eqs.~\eqref{eq:calNh1} through \eqref{eq:calNh3}. 
Fix $c_2>0$. We will define the following linear program, which we refer to as {\lptwo}: 
\begin{align}
\mbox{minimize}  &   \quad\quad\quad m \\
\label{eq:LP2c1}
\mbox{s.t.}  \quad &  (1+\ln^{-1}(1/\delta)) \sum_{h\in \calH} n_{h, k} \pi^P_h \leq  r_k m-c_2, \quad \forall k \in \calK,  \\
& \{n_{h,k}\}_{k\in \calK} \in \calN_{h}, \quad \forall h\in \calH, 
\label{eq:LP2c}
\end{align}
We will denote by $m^*_2$ the optimal value of \lptwo. For the remainder of the proof, we will assume that
\begin{equation}
m q^A > m^*_2, 
\end{equation}
and demonstrate that $(1)$ Eq.~\eqref{eq:contraction1} holds whenever the above inequality is true, and $(2)$ the value of $m^*_2$ is not far from the optimal solution to $\flp$. 

From the definition of a fluid solution, we know that 
\begin{equation}
\dot{\bld}_k(t) = r_kmq^A, \quad \mbox{if }  \blw_k(t)>0, 
\end{equation}
and 
\begin{equation}
  \{\dot{\bla}_k(t)\}_{k \in \calK} \leq (1+\ln^{-1}(1/\delta))  \sum_{h\in \calH}\pi_{h}^P \calN^*_h(\blw(t)). 
\end{equation}  
 For concreteness, we can write  
\begin{equation}
  \{\dot{\bla}_k(t)\}_{k \in \calK} = (1+\ln^{-1}(1/\delta))  \sum_{h \in \calH}\pi_n^P n^*_{h,k}, 
\end{equation}
for some $\{n^*_{h,k}\}_{h\in \calH, k\in \calK}$, where $\{n^*_{h,k}\}_{k\in \calK} \leq \calN^*_h(\blw(t))$ for all $h\in \calH$. We have that
\begin{align}
\sum_{k\in \calK}\blw_k(t) \dot{\blw}_k(t) =  & \sum_{k\in \calK}\blw_k(t) (\dot{\bla}_k(t)-\dot{\bld}_k(t)) \nln
= & \sum_{k\in \calK }\blw_k(t) \lt[ (1+\ln^{-1}(1/\delta)) \sum_{h \in \calH}\pi_h^Pn^*_{h,k}- r_k m q^A \rt] \nln
= & \sum_{h \in \calH} \pi_h^P \lt[ (1+\ln^{-1}(1/\delta)) \sum_{k \in \calK} \blw_k(t)n^*_{h,k} - \sum_{k \in \calK} \blw_k(t)r_kmq^A \rt], 
\label{eq:dotL1}
\end{align}
where the last step follows from the fact that $\sum_{h}\pi^P_h = 1$.  Because $m q^A > m^*_2$, the value $m q^A$ must belong to some feasible solution of \lptwo. By the first constraint of \lptwo in Eq.~\eqref{eq:LP2c1}, this implies that there exists $\{n_{h,k}\}_{h\in \calH, k\in \calK}$, where $\{n_{h,k}\}_{k\in \calK}\in \calN_h$ for all $h\in \calH$, such that 
\begin{equation}
r_k m q^A \geq c_2+ (1+\ln^{-1}(1/\delta)) \sum_{h\in \calH} \pi^P_h n_{h,k}, \quad \forall k\in \calK. 
\end{equation}
We thus have that
\begin{align}
\sum_{k \in \calK} \blw_k(t) r_k m q^A \geq & \sum_{k \in \calK} \blw_k(t) \big(c_2+ (1+\ln^{-1}(1/\delta)) \sum_{h\in \calH} \pi^P_h n_{h,k} \big) \nln
= & c_2 \sum_{k \in \calK}\blw_k(t)+ (1+\ln^{-1}(1/\delta)) \sum_{h \in \calH}\pi^P_h \lt(\sum_{k \in \calK} \blw_k(t)n_{h,k} \rt). 
\label{eq:rkmlwb}
\end{align}
Substituting the above inequality into Eq.~\eqref{eq:dotL1}, we have that
\begin{align}
\sum_{k\in \calK}\blw_k(t) \dot{\blw}_k(t) \leq & - c_2 \sum_{k \in \calK}\blw_k(t)+(1+\ln^{-1}(1/\delta)) \sum_{h \in \calH} \pi^P_h\lt( \sum_{k \in \calK} \blw_k(t)n^*_{h,k} - \sum_{k \in \calK} \blw_k(t)n_{h,k}  \rt) \nln
{\leq} & - c_2 \sum_{k \in \calK}\blw_k(t), 
\label{eq:dotL2}
\end{align}
where the last inequality follows from the fact that $\{n^*_{h,k}\}_{k\in \calK} \leq \calN^*_h(\blw(t))$, and hence 
\begin{equation}
\sum_{k\in \calK} \blw_k(t) n^*_{h,k} \leq \min_{\{n_{k}\}_{k\in \calK}} \sum_{k\in \calK}  \blw_k(t) n_k. 
\end{equation}
Substituting Eqs.~\eqref{eq:dLdt-w0term} and \eqref{eq:dotL2} into Eq.~\eqref{eq:dLdt1}, we have that
\begin{align}
\frac{d}{dt} L(\blw(t)) {=} & \|\blw(t)\|_2^{-1/2}  \lt(\blw_0(t)\dot \blw_0(t) + \sum_{k \in \calK} \blw_k(t)\dot{\blw}_k(t) \rt) \nln
\leq & - \min\{\tilde{c}, c_2\} \|\blw(t)\|_2^{-1/2} \lt( \sum_{k=0}^{ \szK+1} \blw_k(t) \rt) \nln
\sk{a}{\leq} & - \min\{\tilde{c}, c_2\} \|\blw(t)\|_2^{-1/2} \|\blw(t)\|_2 \nln
= & - \min\{\tilde{c}, c_2\} \sqrt{L(\blw(t))}
\label{eq:dLdtnew}
\end{align}
where in step $(a)$ we used the elementary inequality:  $\sum_{i=1}^n x_i \geq  \sqrt{\sum_{i=1}^n x_i^2}, \, \forall (x_1, \ldots, x_n)\in \rp^n$.

Recall that all coordinates of $\blw(\cdot), \bla(\cdot)$ and $\bld(\cdot)$ are {differentiable} for almost all $t\in (0,T)$. Eq.~\eqref{eq:dLdtnew} thus implies that $L(\blw(t))$ is strictly decreasing whenever $\blw(t) \neq 0$. It follows that if $L(\blw(0))=1$, then, for every $t \in (0,T)$, either $L(\blw(t))\leq 1/4$, or $L(\blw(s))> 1/4$ for all $s\in [0,t]$. In the latter case, we have
\begin{equation}
L(\blw(t))\leq 1-  \min\{\tilde{c}, c_2\} \int_{1}^t \sqrt{L(\blw(s))}ds \leq 1-\frac{\min\{\tilde{c}, c_2\} }{2} \, t. 
\end{equation}
Setting the right-hand side of the above equation to $1/4$, we conclude that if  \begin{equation}
\epsilon' =\frac{3}{4} \quad \mbox{and} \quad \tau = \frac{3}{2\min\{\tilde{c}, c_2\} }, 
\end{equation}
then 
 \begin{equation}
L(\blw(\tau)) \leq 1-\epsilon', \quad \mbox{ whenever} \quad L(\blw(0))=1. 
\label{eq:Lgood1}
\end{equation}
Therefore, we have shown that the contraction property of Eq.~\eqref{eq:contraction1} holds,  whenever 
\begin{equation}
m q^P > n^P,  \mbox{ and } m q^A > m_2^*. 
\end{equation}
To complete the proof, it remains to relate the optimal value of \lptwo, $m_2^*$, to that of the the optimal value of \flp. This is accomplished in the following lemma, whose proof is given in Appendix \ref{app:lem:LP2toOrig}.

\begin{lemma} \label{lem:LP2toOrig}
Let $m_F^*$ and $m^*_2$ be the optimal values of \flp and \lptwo, respectively, and $c_2$ be defined as in \lptwo. For every $\epsilon>0$, there exist $c'>0$, such that if $c_2<c'$, then
\begin{equation}
m^*_2 \leq \lt(1+\frac{\ln(2\szH)+ g_\delta}{\ln(1/\delta)}\rt)\lt(1+ \epsilon^P\frac{\szH \old }{\uld\, \ulr}\rt)(1+\ln^{-1}(1/\delta)) m_F^* + \epsilon. 
\end{equation}
\end{lemma}%

The proof of Proposition \ref{prop:drainFluid} is completed by recalling from Proposition \ref{prop:erroAfterPrep} that 
\begin{equation}
\epsilon^P \leq 2\szH/ \ln(1/\delta), 
\end{equation}
and by choosing a sufficiently small $\epsilon$ in Lemma \ref{lem:LP2toOrig}.  \qed

\subsection{Proof of Proposition \ref{prop:residualStable}}
\label{app:prop:residualStable}

\bpf We first look at the macro-level arrival primitive associated with the Residual stage. The arrivals to the Residual stage are jobs leaving the Adaptive stage who have failed to meet the criterion for exiting the system in Eq.~\eqref{eq:adaptiveExitCriteria}. If the Preparation and Adaptive stages are stable, the job arrival process to the Residual stage is generated according to a Markov modulated Poisson process (MMPP), modulated by a stable countable-state Markov process that corresponds to the dynamics of Preparation and Adaptive stages, $\{\calI(t),\{Y_i(t)\}_{i\in \calI(t)}\}_{t \in \rp}$, defined in Section \ref{sec:stateRep}. Using an elementary application of the Foster-Lyapunov criteria, it is thus not difficult to show that the Residual stage is stable whenever the following rate condition holds: the average rate of inspections from the MMPP should be \emph{less} than the service rate of the Residual stage. In the remainder of the proof, we verify that this is indeed the case. 

Recall that $\what{H}^P_i$ is the ML estimator job $i$ obtains from the inspections in the Preparation stage. We say that job $i$ is ``good'', if $\what{H}^P_i$ happens to be equal to the true label, $H_i$, and ``bad'', otherwise. The next result states that most of the good jobs will be able to depart from the system after being processed in the Adaptive stage, and hence not enter the Residual stage. Denote by $S_i^A(h,l)$ the value of $S_{i,t}(h,l)$ when job $i$ exits the Adaptive stage, and define the event 
\begin{equation}
\calB^A = \{\exists h' \in \calH, \mbox{ s.t. } S^A_i(h',l) \geq \ln(2\szH/\delta), \, \forall l \in \calH, l\neq h'\}.
\label{eq:calBA}
\end{equation}
We have the following lemma, whose proof is given in  Appendix \ref{app:lem:ivisitResi}. 
\begin{lemma} 
\label{lem:ivisitResi}
Fix $i\in \N$. There exists $\delta_0 >0$, independent of $\pi$, such that
\begin{equation}
\pb \lt( \calB^A \bbar  \what{H}_i^P = H_i = h\rt) \geq  1- \szH \ln^{-1}(1/\delta), \quad \forall h \in  \calH, 
\end{equation}
for all $\delta \in (0,\delta_0)$. 
\end{lemma}

By Proposition \ref{prop:erroAfterPrep}, we have that at least a $(1-2\szH\ln^{-1}(1/\delta))$ fraction of all jobs are good when exiting the Preparation stage. By Lemma \ref{lem:ivisitResi}, each one of these good jobs has a probability of at most $\szH\ln^{-1}(1/\delta)$ for entering the Residual stage. By assuming that \emph{all} bad jobs eventually enter the Residual stage, the above reasoning shows that the arrival rate of jobs to the Residual stage, $\lambda^R$, satisfies
\begin{equation}
\lambda^R \leq 2\szH\ln^{-1}(1/\delta)+ \szH\ln^{-1}(1/\delta) = 3\szH\ln^{-1}(1/\delta). 
\label{eq:lamRupperB}
\end{equation}
In terms of service speed, recall that each job requires $n^R$ inspections in the Residual stage, and the experts visit the Residual stage at the rate of $m q^R$. Hence, by Eq.~\eqref{eq:lamRupperB}, it suffices to have 
\begin{align}
mq^R >&  3\szH\ln^{-1}(1/\delta) n^R  \nln
= & 3\szH\ln^{-1}(1/\delta) \zeta_0 \ln(4\szH/\delta) \nln
= & 3\szH\zeta_0(1+\ln(4\szH) \ln^{-1}(1/\delta)). 
\end{align}
This proves Proposition \ref{prop:residualStable}. 
\qed

\subsection{Proof of Proposition \ref{prop:classErr}}
\label{app:prop:classErr}

\bpf Recall that there are only two ways through which a job could depart from the system: either from the Adaptive stage, or the Residual stage. Fix $i\in \N$ and $h\in \calH$. Denote by $S^A_i(h,l)$ the value of $S_{i,t}(h,l)$ at the time job $i$ completes all inspections in the Adaptive stage, and by $\what{H}^R_i$ the maximum likelihood estimator of $H_i$ using only the inspections in the Residual stage.  

Recall the event $\calB^A$ defined in Eq.~\eqref{eq:calBA} in Appendix \ref{app:prop:residualStable}, i.e., whether job $i$ will go through the Residual stage, and denote by $\overline{\calB^A}$ its compliment.  We have that: 
\begin{align}
\pb(\what{H}_i \neq h \bbar H_i = h) = & \pb(\what{H}_i \neq h, \calB^A \bbar H_i = h) + \pb(\what{H}_i \neq h, \overline{\calB^A} \bbar H_i = h) \nln
= & \pb(\what{H}_i \neq h, \calB^A \bbar H_i = h) + \pb(\what{H}_i \neq h \bbar \overline{\calB^A}, H_i = h)\pb( \overline{\calB^A} \bbar H_i = h) \nln
\sk{a}{=} & \pb(\what{H}_i \neq h, \calB^A \bbar H_i = h) + \pb(\what{H}^R_i \neq h \bbar  \overline{\calB^A},  H_i = h)\pb( \overline{\calB^A} \bbar H_i = h) \nln
\leq & \pb(\what{H}_i \neq h, \calB^A \bbar H_i = h) + \pb(\what{H}^R_i \neq h \bbar  \overline{\calB^A},  H_i = h), 
\label{eq:whatHbound1}
\end{align}
where step $(a)$ follows from the definition of the inspection policy, where a job will be sent to the Residual stage if and only if the event $\calB^A$ does not occur. The two terms on the right-hand-side of Eq.~\eqref{eq:whatHbound1} correspond to the classifications made by the Adaptive and Residual stages, respectively. By noting that the total number of inspections by job $i$ by the time it exits the Adaptive stage is a stopping time with respect to the past inspections, applying Lemma \ref{lem:SijSufficient2}, with $x = \ln(2\szH/\delta)$, we have that
\begin{equation}
\pb(\what{H}_i \neq h, \calB^A \bbar H_i = h) \leq \delta/2. 
\label{eq:whatHadap}
\end{equation}
Now suppose that job $i$ exits the system after having gone through the Residual stage. Recall that by construction, the Residual stage employs an inspection procedure that is the same as the Preparation stage, except for having a larger number of inspections per job. Therefore, using an essentially identical line of arguments to that of Proposition \ref{prop:erroAfterPrep}, and, in particular, to the portion of the proof leading to Eq.~\eqref{eq:Hpi}, by replacing the quantity $\ln\ln(1/\delta)$ with $\ln(2\szH/\delta)$, we can show that
\begin{equation}
\pb(\what{H}^R_i \neq h \bbar  \overline{\calB^A},  H_i = h) \leq \delta/2. 
\label{eq:whatHResi}
\end{equation}
Substituting Eqs.~\eqref{eq:whatHadap} and \eqref{eq:whatHResi} into Eq.~\eqref{eq:whatHbound1} completes the proof of Proposition \ref{prop:classErr}. 
\qed

\subsection{Proof of Lemma \ref{lem:LP2toOrig}} 
\label{app:lem:LP2toOrig} 

\bpf 


Let $(m_F^*, \{n^*_{h,k}\}) $ be an optimal solution to \flp that satisfies the conditions stated in Lemma \ref{lem:lp1structure}. Define 
\begin{equation}
\phi_\delta = 1+(\ln (2\szH) + g_\delta) {\ln^{-1}(1/\delta)},
\end{equation}
and let
\begin{equation}
\tilde{n}_{h,k} = n^*_{h,k}\phi_\delta. 
\end{equation}
Recall from the definition in Eq.~\eqref{eq:vdeltadef} that 
\begin{equation}
v_\delta = 2\uld^{-1}\ln(1/\delta) [1+(\ln(2\szH) + g_\delta) {\ln^{-1}(1/\delta)}] =  2\uld^{-1}\ln(1/\delta) \phi_\delta, 
\end{equation}
which implies, by Eq.~\eqref{eq:sumnhk1} of Lemma \ref{lem:lp1structure}, that $\sum_{k \in \calK } \tilde{n}_{h,k} \leq v_\delta$ for all $h \in \calH$. We thus conclude that the variables $\{\tilde{n}_{h,k}\}$ satisfy the second set of constraints of \lptwo, i.e., 
\begin{equation}
 \{\tilde{n}_{h,k}\}\in \calN_{h},\quad  \forall h\in \calH. 
\end{equation}

We now derive a sufficient condition on $\tilde{m}$ so that $(\tilde{m}, \{\tilde{n}_{h,k}\})$ is a feasible solution to \lptwo. Fixing $k\in \calK$, we have that
\begin{align}
\sum_{h \in \calH} \tilde{n}_{h,k} \pi^P_h = & \phi_\delta\sum_{h \in \calH} {n}^*_{h,k} \pi^P_h  \nln
\sk{a}{\leq} & \phi_\delta\lt(\sum_{h \in \calH} {n}^*_{h,k} \pi_h  + \epsilon^P\sum_{h \in \calH} {n}^*_{h,k} \rt) \nln
\sk{b}{\leq} & \phi_\delta\lt(\sum_{h \in \calH} {n}^*_{h,k} \pi_h  + \epsilon^P  \szH  \uld^{-1} \ln(1/\delta)\rt) \nln 
\sk{c}{\leq} & \phi_\delta\lt(\sum_{h \in \calH} {n}^*_{h,k} \pi_h  +  \epsilon^P  \szH  \uld^{-1} \old m_F^*  \rt) \nln
\sk{d}{\leq} & \phi_\delta\lt(m_F^* r_k +  \epsilon^P  \szH  \uld^{-1} \old m_F^*  \rt). 
\label{eq:tilnhkpi}
\end{align}
Step $(a)$ follows from Proposition \ref{prop:erroAfterPrep}, $(b)$ and  $(c)$ from Lemma \ref{lem:lp1structure}, and $(d)$ from the constraint of \flp in Eq.~\eqref{eq:Lp1repconstr1}. To satisfy the first constraint of \lptwo, it suffices to let the right-hand side of Eq.~\eqref{eq:tilnhkpi} satisfy
\begin{equation}
\phi_\delta\lt(m_F^* r_k +  \epsilon^P  \szH  \uld^{-1} \old m_F^*  \rt)\leq \frac{r_kq^A \tilde{m}-c_2}{1+\ln^{-1}(1/\delta)}. 
\end{equation}
That is, if 
\begin{equation}
\tilde{m}q^A  \geq m_F^* \phi_\delta \lt(1+ \epsilon^P\frac{ \szH  \old }{\uld\, \ulr}\rt)(1+\ln^{-1}(1/\delta)) + \frac{c_2}{\ulr}, 
\end{equation}
then
\begin{equation}
 \sum_{h\in \calH} \tilde{n}_{h, k} \pi^P_h \leq  r_k q^A \tilde{m}-c_2, \quad \forall k \in \calK, 
\end{equation}
which, in light of the fact that $\{\tilde{n}_{h,k}\}_{k\in \calK}\in \calN_h$ for all $h\in \calH$, further implies that $(\tilde{m}, \{\tilde{n}_{h,k}\})$ is a feasible solution of \lptwo. Therefore, we conclude that, for all $c_2>0$, 
\begin{align}
m^*_2 \leq & m_F^*\phi_\delta\lt(1+ \epsilon^P\frac{ \szH  \old }{\uld\, \ulr}\rt)(1+\ln^{-1}(1/\delta))  + \frac{c_2}{\ulr} \nln
= & m_F^*\lt(1+\frac{\ln(2\szH) + g_\delta}{\ln(1/\delta)}\rt)\lt(1+ \epsilon^P\frac{ \szH  \old }{\uld\, \ulr}\rt)(1+\ln^{-1}(1/\delta)) + \frac{c_2}{\ulr}, 
\end{align}
We complete the proof of Lemma \ref{lem:LP2toOrig} by letting $c' = {\epsilon \ulr}$. 
\qed

\subsection{Proof of Lemma \ref{lem:ivisitResi}}
\label{app:lem:ivisitResi}

\bpf  The proof follows similar steps  as those in the proof of Lemma \ref{lem:reachSijPrep}. Fix $h\in \calH$ and $\{\lambda_{i, k}\}_{k \in \calK} \in \zp^{\szK}$. Denote by $\calB$ the event
\begin{equation}
\calB = \{\what{H}_i^P = H_i=h, \, \{\Lambda_{i, k}\}_{k \in \calK} = \{\lambda_{i, k}\}_{k \in \calK} \}
\end{equation}
For the remainder of the proof, we will assume that $h\in \calH$ and $\{\lambda_{i, k}\}_{k \in \calK}$ are such that $\pb(\calB) >0$.  We will index all inspections received by job $i$ in the Adaptive stage in an arbitrary fashion, and denote by $K_n$ the type of the expert who performed the $n$th inspection for job $i$, and by $Z_{n}(\cdot,\cdot, K_n)$ the corresponding log-likelihood ratio. Define $\Lambda_i = \sum_{k\in \calK}\Lambda_{i,k}$. For all $l \in \calH, l\neq h$,  let
\begin{equation}
M^l_n = \sum_{s=1}^{n \vee \Lambda_i }  Z_{s}(h,l,K_s) - D(h, l, K_s), \quad n \in \N, 
\end{equation}
and $M^l_0=0$. It is not difficult to see that,  conditional on the true label of job $i$ being $h$, the summands in the above equation have zero mean and are independent. Therefore, one can verify that $\{M^l_n \bbar  \calB\}_{n \in \zp}$ is a martingale. 

Define 
\begin{equation}
\phi_\delta = 1+(\ln(2\szH) + g_\delta) {\ln^{-1}(1/\delta)},
\end{equation}
By the construction of $\{\Lambda_{i,k}\}_{k\in \calK}$ (Eq.~\eqref{eq:Lamb_i_k}), we know that 
\begin{equation}
\Lambda_i \leq \sum_{k\in \calK}n_k \leq v_\delta = 2 \phi_\delta \uld^{-1}\ln(1/\delta). 
\label{eq:Lambdaiupperb}
\end{equation}
By Eq.~\eqref{eq:calNh1}, we know that if 
\begin{equation}
M^l_{\Lambda_i} - M^l_{0} >  - g_\delta + \szK \old, 
\end{equation}
then 
\begin{equation}
S_i^A(h,l) > \ln(2\szH/\delta), 
\end{equation}
where the term $\szK \old$ captures the potential discrepency induced by the rounding in Eq.~\eqref{eq:Lamb_i_k}.  We have that 
\begin{align}
\pb\lt(S_i^A(h,l) \leq \ln(2\szH/\delta) \bbar \calB \rt) =& \pb\lt ( M^l_{\Lambda_i} - M^l_{0} \leq - (g^A_\delta - \szK\old )    \bbar \calB \rt) \nln
\sk{a}{\leq} &\exp\lt(- \frac{(g_\delta - \szK \old)^2}{4\olz^2 \Lambda_i }\rt) \nln
\sk{b}{\leq} &\exp\lt(- \frac{  (g_\delta)^2 - 2g_\delta \szK \old }{8\olz^2 \uld^{-1}  \phi_\delta \ln(1/\delta)}   \rt). 
\label{eq:SAilowerbound}
\end{align}
Step $(a)$ follows from the Azuma-Hoeffding's Inequality (Lemma \ref{lem:asuma} in Appendix \ref{app:azuma}), and noting that $|M^l_n-M^l_{n-1}| \leq 2\olz$ for all $n\in \N$. Step $(b)$ follows from Eq.~\eqref{eq:Lambdaiupperb}. 

Recall that $g_\delta = 3\olz\uld^{-1/2} \sqrt{ \ln(1/\delta)\ln\ln(1/\delta) }$. The exponent in Eq.~\eqref{eq:SAilowerbound} can be further expanded as: 
\begin{align}
- \frac{  (g_\delta)^2 - 2g_\delta \szK \old }{8\olz^2 \uld^{-1}  \phi_\delta \ln(1/\delta)} {=}&  - \frac{9\olz^2\uld^{-1}\ln(1/\delta)\ln\ln(1/\delta)}{8\olz^2 \uld^{-1}  \phi_\delta \ln(1/\delta)} + \frac{ 6\olz\uld^{-1/2} \szK \old \sqrt{\ln(1/\delta)\ln\ln(1/\delta)} }{ 8\olz^2 \uld^{-1}  \phi_\delta \ln(1/\delta)} \nln
= & - \frac{9 \ln\ln(1/\delta) }{8\phi_\delta} +\frac{3\uld^{-3/2} \szK \old }{4\olz } \sqrt{\frac{\ln\ln(1/\delta)}{\ln(1/\delta)}}. 
\end{align}
As $\delta \downarrow 0$, we have that $\phi_\delta=1+(\ln(2\szH) + g_\delta) {\ln^{-1}(1/\delta)} \to 1$, and $\sqrt{\frac{\ln\ln(1/\delta)}{\ln(1/\delta)}} \to 0$. Therefore, the above expression yields\footnote{The notation $f(x)\sim g(x)$ means that $\lim_{x\to 0}\frac{f(x)}{g(x)} = 1$. }
\begin{equation}
- \frac{  (g_\delta)^2 - 2g_\delta \szK \old }{8\olz^2 \uld^{-1}  \phi_\delta \ln(1/\delta)} \sim - \frac{9 \ln\ln(1/\delta) }{8}, \quad \mbox{as $\delta \to 0$.}
\label{eq:exponenMarting} 
\end{equation}
Substituting  Eq.~\eqref{eq:exponenMarting} into Eq.~\eqref{eq:SAilowerbound}, we know that there exists $\delta_0>0$, independent of the prior distribution, $\pi$, and the choice of $l$ and $h$, such that for all $\delta \in (0,\delta_0)$, 
\begin{align}
\pb\lt(S_i^A(h,l) \leq \ln(2\szH/\delta) \bbar \calB \rt) \leq & \exp\lt( -\frac{  (g_\delta)^2 - 2g_\delta \szK \old }{8\olz^2 \uld^{-1}  \phi_\delta \ln(1/\delta)} \rt) \nln
= & \exp\lt( - \frac{9 \ln\ln(1/\delta) }{8\phi_\delta} +\frac{3\uld^{-3/2} \szK \old }{4\olz } \sqrt{\frac{\ln\ln(1/\delta)}{\ln(1/\delta)}}  \rt) \nln
\leq& \exp(- \ln\ln(1/\delta))\nln
= & \ln^{-1}(1/\delta). 
\end{align}
Since the above inequalities hold for all $h,l \in \calH$, $l \neq h$, and $\{\lambda_{i, k}\}_{k \in \calK}$, we can apply a union bound over $l$ and conclude that, for all $h\in \calH$, 
\begin{align}
\pb(\calB^A \bbar \what{H}_i^P = H_i = h) \geq & 
\pb \lt( S_i^A(h,l) \geq \ln(2\szH/\delta),  \, \forall l \neq h\bbar \what{H}_i^P = H_i = h\rt) \geq  1- \szH \ln^{-1}(1/\delta).
\end{align}
This completes the proof of Lemma \ref{lem:ivisitResi}. 
\qed

\section{Some Properties of $L(\cdot)$}
\label{app:lem:L1homo}

The following lemma states some basic properties of the Lyapunov function, $L(\cdot)$. 
\begin{lemma}
\label{lem:L1homo}
The function $L(\blw) = \|\blw\|_2$ admits the following properties: $(1)$ $L$ is continuous on $\rp^{ \szK +1 }$; $(2)$ For all $c \in \rp$ and $\blw\in \rp^{ \szK+1}$, we have  $L(c \blw) = cL(\blw)$; $(3)$ There exists positive constants, $\alpha_1, \alpha_2>0$, such that for all $\blw\in \rp^{ \szK+1}$, 
\begin{equation}
\alpha_1\|\blw\|_{\infty} \leq L(\blw) \leq \alpha_2 \|\blw\|_{\infty},
\label{eq:L1homo}
\end{equation}
where $\|\cdot\|_\infty$ is the $l_\infty$ norm, $\|\blw\|_\infty =\max_{k = 0, \ldots,  \szK+1}|\blw_k|$. 
\end{lemma}

\bpf The first two properties follow directly from the definition of $L(\cdot)$. For Eq.~\eqref{eq:L1homo}, we note that for all $\blw\in \rp^{ \szK+1}$, 
\begin{align}
L(\blw) \geq & \sqrt{\|\blw\|_\infty^2} = \|\blw\|_\infty, \nln
L(\blw) \leq & \sqrt{( \szK+1)\|\blw\|_\infty^2} = \sqrt{ \szK+1} \|\blw\|_\infty. \nonumber
\end{align}
Hence, Eq.~\eqref{eq:L1homo} holds for $\alpha_1 = 1$ and $\alpha_2 =  \sqrt{ \szK+1}$. \qed

\section{The Azuma-Hoeffding Inequality}
\label{app:azuma}
We will make use of the standard Azuma-Hoeffding inequality (Section 12.2 of \cite{grimmett2001probability}), which we state below for completeness. 
\begin{lemma}[Azuma-Hoeffding Inequality]
\label{lem:asuma} Let $\{M_n\}_{n \in \zp}$ be a martingale, and suppose that there exists $a>0$, such that $|M_{n}-M_{n-1}|\leq a$ for all $n\in \zp$. Then
\begin{equation}
\pb(M_n - M_0 \leq -x ) \leq \exp\lt(- \frac{x^2}{2n a^2}\rt), \quad \forall x \in \rp. 
\end{equation}
\end{lemma}

\section{A Heuristic Policy}
\label{sec:heuristicPolicy}
We discuss in this section a simple heuristic inspection policy, which essentially condenses the three-stage architecture of Figure \ref{fig:blockDiag} into only one that resembles the Adaptive stage. We expect the policy to be substantially easier  to implement than the one presented in Section \ref{sec:DefOfPolicy}. Nevertheless, we have not been able to establish its resource efficiency rigorously, though we discuss some indication that it may be resource efficient as well. 

For simplicity, we assume that $D(h,l,k)>0$ for all $k\in \calK, h\neq l$, i.e., every expert is able to distinguish any two job labels to some degree.  Denote by $\what{H}_{i,t}$ the ML estimator of job $i$'s type at time $t$, and we assume that the value of $\what{H}_{i,t}$ at the time when job $i$ has just arrived to the system is set to a value drawn uniformly at random from $\calH$. Denote by $\calI_t$ the indices of all jobs that are currently in the system at time $t$, and denote by $\calI_{h,t}$ the subset containing those jobs whose ML estimator is equal to $h$:
\begin{equation}
\calI_{h,t} = \{i \in \calI_t: \what{H}_{i,t}=h\}.
\end{equation}
Define $W_{i,t}(h,l)$ as the \emph{residual workload} associated with job $i$: 
\begin{equation}
W_{i,t}(h,l) = \lt(\ln( \szH/ \delta) - S_{i,t}(h,l) \rt)^+, \quad h,l\in \calH, 
\end{equation}
and define the \emph{aggregate residual workload} associated with job type $h$:  
\begin{equation}
\overline{W}_{t} (h,l) = \sum_{i \in \calI_{h,t} } W_{i,t}(h,l) , \quad h, l \in \calH. 
\end{equation}

The heuristic policy operates as follows.

\emph{Departure rule.} Job $i$ departs from the system as soon as there exists $h\in \calH$, such that 
\begin{equation}
\max_{l\in \calH, l\neq h}  W_{i,t}(h,l) = 0.  
\label{eq:indiviWitdef}
\end{equation}

\emph{Expert actions.} Suppose that an expert of type $k$ becomes available at time $t$. If there is no job in the system, she goes on a vacation. Otherwise, she inspects the oldest job from the set $\calI_{h^*,t}$, where 
\begin{equation}
h^* \in \argmax_{h \in \calH} \sum_{l \in \calH, l\neq h} D(h,l,k)\overline{W}_{t} (h,l), 
\label{eq:serAdap}
\end{equation}
with ties broken arbitrarily. 

This heuristic policy has a number of nice features. As is evident from the above description, the policy is much simpler than the one given in Section \ref{sec:DefOfPolicy}, and it does not require solving a linear optimization problem as a sub-routine. Also, the departure criterion ensures that the event $\calG_x$ in Lemma \ref{lem:SijSufficient} always occurs for every job that leaves the system, so the policy is $\delta$-accuracy whenever it is stable. There is also some indication that this heuristic policy might be resource efficient. For a moment, let us suppose that all ML estimators,  $\what{H}_{i,t}$, are correct. In this case, $\calI_{h,t}$ is the set of jobs whose true label is $h$, and Eq.~\eqref{eq:serAdap} can be thought of as a max-weight-like procedure, where $D(h,l,k)$ corresponds to the expected decrease a single inspection can incur on the aggregate residual workload $\overline{W}_t(h,l)$. Since max-weight scheduling policies (e.g., \cite{tassiulas1992stability}) are known to achieve the maximum stability region in many queueing systems, one may expect that this max-weight-like property would make our heuristic policy resource efficient as well. 

Unfortunately, there appear to be two intrinsic characteristics  of the policy that make it difficult to rigorously establish resource efficiency. Firstly, the ML estimators $\what{H}_{i,t}$ can be \emph{incorrect}, especially when a job has received only a small number of inspections, and therefore the vector $\{\overline{W}_{t} (h,l)\}_{h,l\in \calH}$ does not exactly capture the ``true'' workload of the system. Secondly, even assuming that the ML estimators are correct, a {synchronization} issue can prevent a type-$k$ expert's inspection from actually contributing an expected $D(h,l,k)$ units of decrease in $\overline{W}_t(h,l)$ for all $l\neq h$. This is because, for some job $i$, ${W}_{i,t}(h,l)$ along some coordinate $l$ could reach zero much earlier than other coordinates $l'\neq l$, during which period the expected decrease in $\overline{W}_t(h,l)$ from inspecting job $i$ would be zero, instead of $D(h,l,k)$, due to the capping at zero of Eq.~\eqref{eq:indiviWitdef}.  Both of the above issues make it difficult to bound the drift in a Lyapunov function in the course of proving stability. However, there is hope that they could be mitigated as the accuracy parameter, $\delta$, tends to $0$. In this regime, the number of inspections needed for each job tends to infinity, and hence most jobs in the system have already received a large number of inspections, implying that they likely have the correct ML estimators as well. Similarly, since the distortion of drift only occurs at the ``boundary'' where the residual workload of a job is small, it will likely have  much less impact as the number of inspections per job grows.

\end{APPENDICES}

\end{document}